\documentclass[journal]{IEEEtran}
\usepackage{amsmath,amssymb,amsfonts}
\usepackage{cases}
\usepackage{color,xcolor}
\usepackage{graphicx}
\usepackage{cite}
\usepackage{subfigure}
\usepackage{epstopdf}
\usepackage{multirow}
\usepackage{rotating}
\usepackage{booktabs}
\usepackage{stmaryrd}
\usepackage{tikz}
\usepackage{acronym}
\usepackage{pifont}
\usepackage{subfigure}
\usepackage{makecell}         
\usepackage{newtxmath}
\usepackage{array}
\usepackage{threeparttable}
\usepackage{pdflscape}
\usepackage{varwidth}
\usepackage{hyperref}
\usepackage[ruled]{algorithm2e}  

\newtheorem{Theorem}{Theorem}[section]
\newtheorem{Lemma}[Theorem]{Lemma}

\newtheorem{Definition}[Theorem]{Definition}

\begin{document}
	
	\title{Sparse Tensor CCA via Manifold Optimization\\ for Multi-View Learning}
	
	\author{Yanjiao Zhu, Wanquan Liu, \IEEEmembership{Senior~Member,~IEEE},
		Xianchao Xiu, \IEEEmembership{Member,~IEEE}, and Jianqin Sun

		\thanks{This work was supported in part by the National Natural Science Foundation of China under Grants 12371306 and 12271309, and in part by the Department of Science and Technology of Guangdong Province under Grant 2021CX02G450. \textit{(Corresponding author: Xianchao Xiu.)}}
        
\thanks{Y. Zhu and W. Liu are with the School of Intelligent Systems Engineering, Sun Yat-sen University, Guangzhou 510275, China (e-mail: \{zhuyj87; liuwq63\}@mail.sysu.edu.cn).}%
\thanks{X. Xiu is with the School of Mechatronic Engineering and Automation, Shanghai University, Shanghai 200444, China
        (e-mail: xcxiu@shu.edu.cn).}
\thanks{J. Sun is with the School of Mathematics and Statistics, Beijing Jiaotong University, Beijing 100044, China (e-mail: 21121631@bjtu.edu.cn).}
	}
	
	\IEEEpubid{\begin{minipage}{\textwidth}\ \\[30pt] \centering
			Copyright \copyright 2025 IEEE. Personal use of this material is permitted. 
			However, permission to use this material for any other purposes must \\ be obtained 
			from the IEEE by sending an email to pubs-permissions@ieee.org.
	\end{minipage}}
	
	\maketitle
	
	\begin{abstract}
		Tensor canonical correlation analysis (TCCA) has garnered significant attention due to its effectiveness in capturing high-order correlations in multi-view learning. However, existing TCCA methods often underemphasize the characterization of individual structures and lack algorithmic convergence guarantees. In order to deal with these challenges, we propose a novel sparse TCCA model called STCCA-L, which integrates sparse regularization of canonical matrices and Laplacian regularization of multi-order graphs into the TCCA framework, thereby effectively exploiting the geometric structure of individual views. To solve this non-convex model, we develop an efficient alternating manifold proximal gradient algorithm based on manifold optimization, which avoids computationally expensive full tensor decomposition and leverages a semi-smooth Newton method for resolving the subproblem. Furthermore, we rigorously prove the convergence of the algorithm and analyze its complexity. Experimental results on eight benchmark datasets demonstrate the superior classification performance of the proposed method. Notably, on the 3Sources dataset, it achieves improvements of at least 4.50\% in accuracy and 6.77\% in F1 score over competitors.
Our code is available at \href{https://github.com/zhudafa/STCCA-L}{https://github.com/zhudafa/STCCA-L}.
	\end{abstract}

	\begin{IEEEkeywords}
		Multi-view learning, tensor canonical correlation analysis (TCCA), sparse regularization, multi-order graph, manifold optimization.	
	\end{IEEEkeywords}
	
	\section{Introduction}
	\IEEEPARstart{W}{ulti-view} learning aims to address the challenge of data heterogeneity, which often arises when the same phenomenon is observed through different modalities or sources, such as images, audio, and textual metadata \cite{Li2019}. This paradigm provides a principled framework for integrating such complementary views to improve data representation, analysis, and interpretation \cite{Ke2024, Wang2025Comprehensive}. Generally, multi-view learning methods fall into three main categories: co-training based methods that iteratively refine classifiers across views via mutual agreement \cite{liu2021multiview,li2023multi}, multi-kernel learning methods that integrate heterogeneous information through composite kernels \cite{celikkanat2022multiple,li2022consensus}, and subspace learning methods that seek a common latent representation shared by multiple views \cite{liang2024robust,teng2025consensus}.
Among them, multi-view subspace learning \cite{Mandal2023, Cheng2025} has attracted growing interest owing to its effectiveness in capturing consensus structures across modalities while preserving complementary information, thereby facilitating robust performance in downstream tasks such as classification \cite{liu2025reliable}, clustering \cite{fang2023comprehensive}, and retrieval \cite{Qilin2023}.

 Canonical correlation analysis (CCA) is a foundational method in multi-view subspace learning \cite{weenink2003canonical,Yang2021}. By seeking projections of each view that maximize the correlation in the lower-dimensional space, CCA effectively captures the most significant and discriminative features shared across modalities \cite{Hu2014Multiview,shu2022d}. Matrix CCA \cite{xu2019canonical} is a popular multi-view CCA method that generalizes pairwise correlation analysis by employing a matrix formulation. This includes various extensions such as CCA \cite{witten2009penalized}, sparse CCA (SCCA) \cite{zhang2014mining}, and structured generalized CCA (SGCCA) \cite{LV2024}. As a nonlinear extension of matrix CCA, deep CCA \cite{Kumar2024,Xiu2022} mines more complex data associations by passing observations to deep neural networks \cite{chen2021comparative,wang2015unsupervised}, autoencoders \cite{kaloga2021variational}, and convolutional networks \cite{yang2017canonical,zhang2021canonical}. It has outstanding feature extraction capabilities in the era of big data. However, deep models often suffer from poor interpretability and a high dependency on large labeled datasets for effective training, limiting their applicability. In contrast, tensor CCA (TCCA) \cite{yang2019survey} leverages high-order covariance structures to model complex relationships among multiple views more effectively. Specifically, Luo \textit{et al.} \cite{luo2015tensor} was the first to propose TCCA, which captures complex, high-order dependencies that are often missed by matrix methods, leading to improved performance in multi-view learning tasks. However, a key limitation of TCCA is its failure to enforce the orthogonality of the regularization variables, which may result in redundant or highly correlated canonical vectors. To address this issue, Sun \textit{et al.} \cite{sun2023learning} integrated TCCA with orthogonality (TCCA-O) to ensure the irrelevance among canonical vectors. For multi-view tensor data, methods like multi-view graph CCA (TMCCA) \cite{Reddy2025}, trial selection TCCA \cite{Yang2025Trial}, and TCCA across multiple groups \cite{Zhou2024MG-TCCA} were widely used, especially in the biomedical field. However, due to the correlation within the constructed matrices, these methods exhibit suboptimal performance in capturing the complexity of data relationships.

\begin{figure*}[t]
\centering
\includegraphics[scale=0.5]{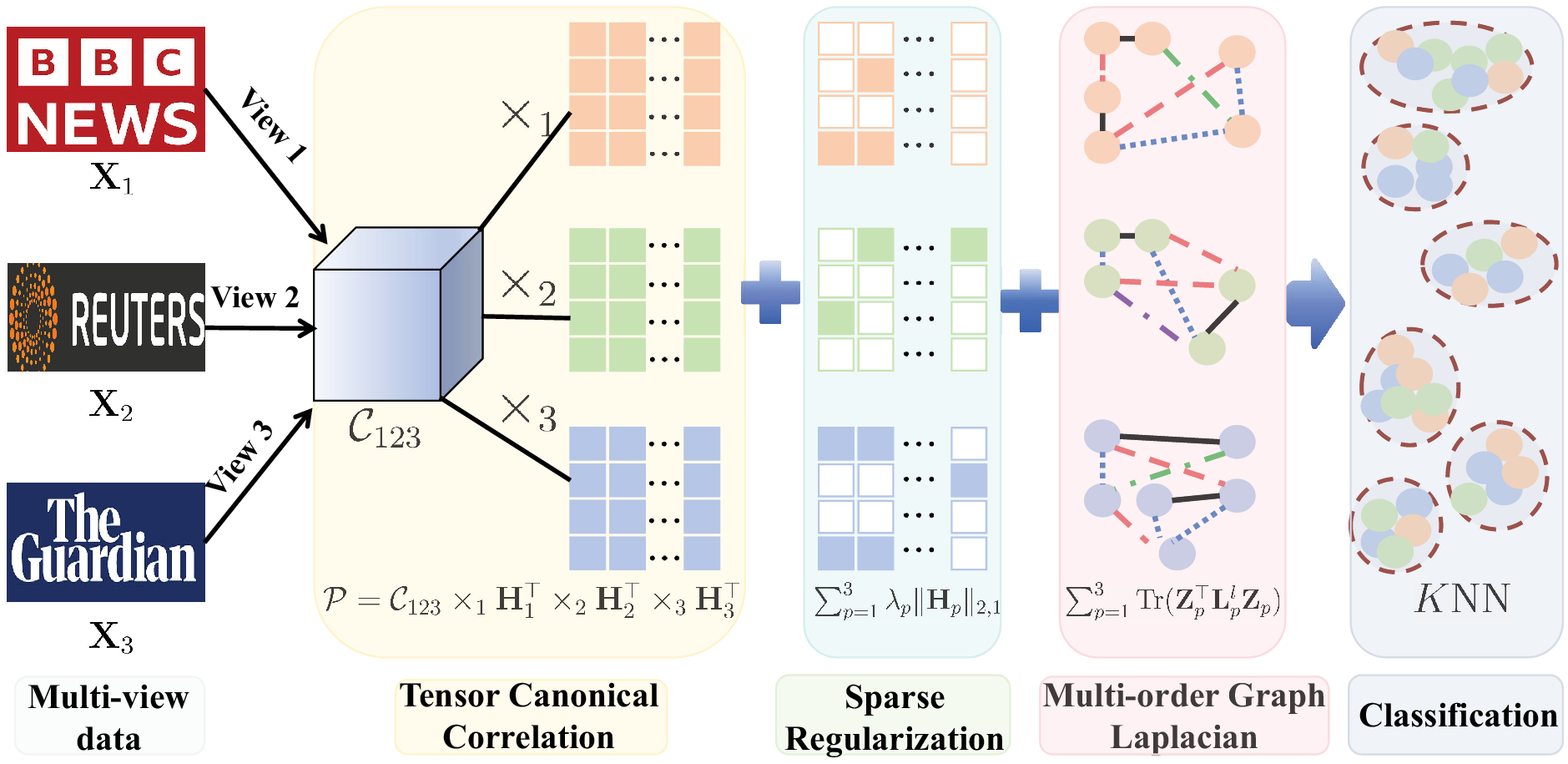}
\caption{Framework of our proposed method. Take the 3Sources dataset as an example. Given the data $\mathbf{X}=[\mathbf{X}_1,\mathbf{X}_2,\mathbf{X}_3] $ from three views, STCCA-L learns projection matrices $\mathbf{H}_1,\mathbf{H}_2,\mathbf{H}_3$ that maximize high-order canonical correlation via a core tensor $\mathcal{P}=\mathcal{C}_{123}\times_{1} \mathbf{H}_{1}^\top \times_{2} \mathbf{H}_{2}^\top\times_{3} \mathbf{H}_{3}^\top$. Moreover, the $\ell_{2,1}$-norm regularization is introduced to promote sparsity, and the Laplacian regularization of multi-order graphs is employed to preserve the intrinsic local structure within each view. Finally, the learned subspace is used to perform the classification task.}\label{ill}
\end{figure*}

Despite their effectiveness in capturing high-order correlations, TCCA methods often suffer from feature redundancy, where many components in canonical projection matrices contribute little to the representation and are difficult to interpret \cite{luo2018semi,wang2024one}. Sparse learning is an effective strategy to address this issue by promoting compact and meaningful representations \cite{li2020survey,zhu2025joint}. Recognizing its importance, Du \textit{et al.} \cite{du2020mining} incorporated a sparse regularization term into the TCCA objective function, enabling feature selection while analyzing complex high-order relationships in multi-modal brain imaging data. On the other hand, Sun \textit{et al.} \cite{sun2023learning} further introduced a structural sparse regularization term in TCCA-O to develop the TCCA-OS method. Due to the sparsity effect on the sample representation after projection, the focus is on selecting key samples.  However, its sparse structure complicates precise control over feature selection, limiting its operability and analytical tractability in the feature space. This limitation hinders the ability of the model to fully leverage high-dimensional data.

Furthermore, most existing TCCA methods lay particular emphasis on capturing structural relationships between views, but often underemphasizing the inherent characteristics of each view \cite{chen2019graph}. This inadequately addressing may lead to the loss of key features \cite{Zheng2024Flexible}. Graph learning offers a complementary solution by representing data as graphs, where the nodes correspond to the data points and the edges encode pairwise relationships \cite{liang2020multi}. Recent studies suggest that high-order graphs provide richer representations by capturing multi-point interactions \cite{Lee2025Hypergraph}. Multi-order graph learning further enhances flexibility by adaptively integrating graphs of varying orders using weighted schemes, and has shown promising results in multi-view tasks \cite{Wang2024Multi-View,liu2024multi}. 

It is worth noting that theoretical analyses of algorithm convergence have received limited attention in existing TCCA methods. Among the few exceptions, Du \textit{et al.} \cite{du2020mining} proved the monotonic increase of the objective function for their proposed method. However, most existing methods lack rigorous guarantees on the convergence of the optimization algorithm, which raises concerns about the stability and reliability of the learned representations in practice. This gap highlights the urgent need for theoretical guarantees of the TCCA methods. 

Motivated by these insights, we propose a novel method called sparse TCCA with the Laplacian regularization of multi-order graph (STCCA-L). 
One goal is to automatically select important features and eliminate redundant dimensions by leveraging the structural sparsity constraints on the projection matrix. The other goal is to effectively capture the intrinsic information of each view of data using the multi-order graph Laplacian regularization. Obviously, the introduction of regularization terms increases the computational load. To handle this limitation, we develop an alternating manifold proximal gradient algorithm based on Stiefel manifold optimization, which ultimately enables our proposed method to achieve good accuracy with acceptable computational efficiency. Taking the 3Sources dataset as an example, the framework of STCCA-L is illustrated in Fig. \ref{ill}.

Compared with the existing work, the main contributions of this paper can be summarized in the following three aspects.
	\begin{enumerate}
    \item \textit{(New Model)} We construct a new multi-view subspace learning model that not only introduces structural sparse regularization to effectively alleviate feature redundancy, but also enhances the exploration of the underlying data structure of each view through multi-order graph Laplacian regularization. To our knowledge, this is the first study to integrate multi-order graphs with TCCA.
    \item \textit{(Convergent Algorithm)} We develop an efficient alternating manifold proximal gradient algorithm on the Stiefel manifold by leveraging the semi-smooth Newton method (SSN). Mathematically, it is rigorously proved that our proposed algorithm converges to a stationary point. 
    \item \textit{(Empirical Superiority)} We validate the effectiveness, robustness, and stability of the proposed method by comparing it with some state-of-the-art CCA methods in the classification task on eight multi-view datasets.
	\end{enumerate}

The structure of this paper is as follows.
Section \ref{sec2} introduces notations and related basics.
Section \ref{sec3} formulates our model and develops the optimization algorithm.
Section \ref{sec4} validates the superiority of our proposed method.
Section \ref{sec5} concludes this paper.
	
	\section{Preliminaries} \label{sec2}
	\subsection{Notations}
	For clarification, the following notation conventions are used: Calligraphic letters for tensors, say $\mathcal{X}$; Bold capital letters for matrices, say $\mathbf{X}$; Bold lowercase letters for vectors, say $\mathbf{x}$; Lowercase letters for scalars, say $x$. For a matrix $\mathbf{X} \in \mathbb{R}^{n\times r}$, the $\ell_{2,1}$-norm is defined by $\|\mathbf{X}\|_{2,1}=\sum_{i=1}^n \|\mathbf{x}_{i}\|_2,$ where $\mathbf{x}_{i}$ is the $i$th row of the matrix $\mathbf{X}$.  
The operator $\operatorname{vec}(\mathbf{X}) \in \mathbb{R}^{nr}$ denotes the vector obtained by stacking the column vectors of $\mathbf{X}$.  
When $\mathbf{X} \in \mathbb{R}^{r\times r}$ is symmetric, let $\overline{\operatorname{vec}}(\mathbf{X}) \in \mathbb{R}^{\frac{1}{2}r(r+1)}$ denote the vector obtained from $\operatorname{vec}(\mathbf{X})$ by eliminating all super-diagonal elements of $\mathbf{X}$.
 For tensors $\mathcal{A},\mathcal{B} \in \mathbb{R}^{I_{1}\times \cdots\times I_{N}}$, their inner product is \begin{equation}
			\langle\mathit{\mathcal{A}, \mathcal{B}}\rangle=\sum_{i_1 =1}^{I_1}\cdots\sum_{ i_N =1}^{I_N} \mathit{a_{i_1, \cdots, i_N} b_{i_1, \cdots, i_N}},
		\end{equation}
and their outer product is $
		\mathcal{A} \circ \mathcal{B}\in \mathbb{R}^{I_1 \times \cdots \times I_N \times I_1 \times \cdots \times I_N},$
whose entries are composed by
 \begin{equation}
			(\mathit{\mathcal{A} \circ \mathcal{B})}_{i_1, \cdots, i_N, i_1, \cdots, i_N}=\mathit{a_{i_1, \cdots, i_N} b_{i_1, \cdots, i_N}}.			
		\end{equation}
For a tensor $\mathcal{A}$, if $\mathbf{V} \in \mathbb{R}^{r_n \times I_n}$ is a matrix, the $n$-mode product of $\mathcal{A}$ with $\mathbf{V}$ is denoted as
\begin{equation*}\mathcal{A} \times_{n} \mathbf{V} \in \mathbb{R}^{I_{1} \times \cdots \times I_{n-1} \times r_n \times I_{n+1}\times \cdots \times I_{N}}.\end{equation*} 
If $\mathbf{v} \in \mathbb{R}^{I_n}$ is a vector, the $n$-mode product of $\mathcal{A}$ with $\mathbf{v}$ is
\begin{equation*}
\mathcal{A} \times_{n} \mathbf{v} \in \mathbb{R}^{I_{1} \times \cdots \times I_{n-1} \times I_{n+1} \times \cdots \times I_{N}}.
\end{equation*} 

In what follows, denote $p=\{1,\cdots,N\}$ as $p\in [N]$.
Given a set of matrices $\{\mathbf{V}_p\}$ where $\mathbf{V}_p\in \mathbb{R}^{r_p \times I_p}$ and $p\in[N]$, the contracted tensor product of $\mathcal{A}$ with ${\mathbf{V}_p}$ is expressed as
		\begin{equation}
			\mathcal{B}=\mathcal{A} \times_1 \mathbf{V}_1 \times_2 \cdots \times_N \mathbf{V}_N \in \mathbb{R}^{r_1 \times \cdots \times r_N}.			
		\end{equation}
 Accordingly, the mode-$p$ unfolding matrix of the tensor $\mathcal{B}$ can be given by
	\begin{equation}
		\mathcal{B}_{(p)}=\mathbf{V}_p \mathcal{A}_{(p)}(\mathbf{V}_{N-1} \otimes \cdots \otimes \mathbf{V}_{p+1} \otimes \mathbf{V}_{p-1} \otimes \cdots \otimes \mathbf{V}_{1})^\top,
	\end{equation}
	where $\otimes$ is the Kronecker product.
		
Below, some definitions related to manifold optimization are introduced.
	
	\begin{Definition}[Stiefel Manifold]\label{ste}
		For a matrix $\mathbf{X}\in \mathbb{R}^{n \times r}$, the Stiefel manifold is
		\begin{equation}
			\textrm{St($n,r$)}=\{ \mathbf{X} \in \mathbb{R}^{n \times r} \mid  \mathbf{X}^\top\mathbf{X}=\mathbf{I}_r\},
		\end{equation}
		where $\mathbf{I}_{r}$ denotes a $r \times r$ identity matrix. \end{Definition}
        It is an orthogonality constraint on the mapping matrix $\mathbf{X}$.
		Its tangent space at a point $\mathbf{X} \in \textrm{St($n,r$)}$ can be expressed as
		\begin{equation}
	\textrm{T}_\mathbf{X}\textrm{St($n,r$)}=\{ \mathbf{U} \in \mathbb{R}^{n \times r} \mid\mathbf{X}^\top\mathbf{U}+\mathbf{U}^\top\mathbf{X}=0\}.
		\end{equation}

	\begin{Definition}[Retraction]
Let $\mathcal{M}$ be a Riemannian manifold, and $\textrm{T}_{\textbf{X}}\mathcal{M}$ be the tangent space of $\mathcal{M}$ at the point $\mathbf{X} \in \mathcal{M}$. A retraction is a smooth mapping defined as
\begin{equation}
\operatorname{Retr}_{\mathbf{X}} : \textrm{T}_{\mathbf{X}}\mathcal{M} \rightarrow \mathcal{M}.
\end{equation}
\end{Definition}
The retraction onto the Euclidean space is simply the identity mapping, \textit{i.e.}, $\operatorname{Retr}_{\mathbf{X}}(\mathbf{y})=\mathbf{X}+\mathbf{y}$. When a point $\mathbf{X} \in \textrm{St($n,r$)}$, QR-based retraction is a common approach for handling the Stiefel manifold, \textit{i.e.}, 
\begin{equation}
			\operatorname{Retr}_{\mathbf{X}}^{\textrm{QR}}(\mathbf{y})=\textrm{qf}(\mathbf{X}+\mathbf{y}),		
		\end{equation}
where $\textrm{qf}(\mathbf{A})$ is the $Q$ factor of the QR factorization of $\mathbf{A}$.

	\begin{Definition}[Proximal Operator]
		For two matrices $ \mathbf{X}, \mathbf{Y}\in \mathbb{R}^{n \times r}$ and a parameter $\beta>0$, the proximal operator is defined as
		\begin{equation}
			\operatorname{prox}_{2,1}( \mathbf{X},\beta)=\underset{ \mathbf{Y}}{\operatorname{argmin}}\{\| \mathbf{Y}\|_{2,1}+\frac{1}{2 \beta}\| \mathbf{Y}- \mathbf{X}\|_\textrm{F}^2\},
		\end{equation}
		whose $i$th row admits the closed-form expression
		\begin{equation}
			\mathbf{y}_i=\frac{\mathbf{x}_i}{\|\mathbf{x}_i\|_2}\max\{0, \|\mathbf{x}_i\|_2-\beta\},
		\end{equation}
		where $\mathbf{x}_i$ and $\mathbf{y}_i$ are the $i$th row of $\mathbf{X}$ and $\mathbf{Y}$, respectively. More details can be found in \cite{nie2010efficient}.
	\end{Definition}

	\subsection{Tensor CCA}
	 For the multi-view data $\mathbf{X} = [\mathbf{X}_1,\cdots,\mathbf{X}_m]$, where each view $\mathbf{X}_p \in \mathbb{R}^{d_p \times N}$  corresponds to a different feature representation of the same instances. TCCA \cite{luo2015tensor} is a representative method that models the high-order correlations among all views by forming a covariance tensor \begin{equation}
\mathcal{C}_{12 \cdots m}=\frac{1}{N} \sum_{n=1}^N \mathbf{x}_{1n}  \circ \mathbf{x}_{2n}\circ \cdots \circ \mathbf{x}_{mn} \in \mathbb{R}^{d_1 \times \cdots \times d_m},
\end{equation} where $\circ$ represents the outer product. The goal is to find a set of projection vectors $\{\mathbf{h}_{p}\},~p\in[m]$ that maximally correlate the projected features. The optimization problem of TCCA can be formulated as \begin{equation}\label{tcca}
		\begin{aligned}
			\underset{\{\mathbf{h}_{p}\}}{\operatorname{max}} &\quad \mathcal{C}_{12 \cdots m} \times_{1} \mathbf{h}_{1}^\top \times_{2}\cdots \times_{m} \mathbf{h}_{m}^\top \\
			\text { s.t. } &~~~ \mathbf{h}_{p}^\top \mathbf{C}_{p p}\mathbf{h}_{p}=1, p\in [m],
		\end{aligned}
	\end{equation}
where  $\mathbf{C}_{pp}=\mathbf{X}_p\mathbf{X}_p^\top$ denotes the variance matrix of the $p$th view. The formulation \eqref{tcca} provides a compact way to explore common latent directions. However, it is limited to discovering only one-dimensional projections for each view, and it does not guarantee the irrelevance among the canonical vectors, which may lead to highly correlated representations.

To overcome this drawback, TCCA-O \cite{sun2023learning} extended TCCA to learn multiple projection directions simultaneously by seeking a set of orthogonality projection matrices $\mathbf{H}_{p}=[\mathbf{h}_{p1},\cdots,\mathbf{h}_{pr}] \in \mathbb{R}^{d_p \times r}, p\in [m]$. The model is formulated as follows
	\begin{equation}\label{tcca-o}
		\begin{aligned}
			\underset{\left\{\mathbf{H}_{p}\right\}}{\operatorname{max}} &\quad \frac{1}{2}\|\mathcal{C}_{12 \cdots m} \times_{1} \mathbf{H}_{1}^\top \times_{2}\cdots \times_{m} \mathbf{H}_{m}^\top\|_\textrm{F}^{2} \\
			\text { s.t. } &~~~ \mathbf{H}_{p}^\top \mathbf{C}_{p p}\mathbf{H}_{p}=\mathbf{I}, p\in [m].
		\end{aligned}
	\end{equation}
   It is verified that TCCA-O not only retains the capability of capturing high-order correlations through tensor modeling, but also ensures that the learned projections form orthogonality bases within each view, thereby avoiding redundancy and enhancing representation capacity. This extension is particularly beneficial in downstream tasks.

\section{The Proposed Method}\label{sec3}

	\subsection{Problem Formulation}
    In graph learning, $\mathbf{X}_p$ corresponds to a weighted undirected graph $\mathbf{G}_p=(\mathbf{V}_p, \mathbf{W}_p)$, where $\mathbf{V}_p$ is the vertex set with the node set of $\mathbf{X}_p$ and $\mathbf{W}_p \in \mathbb{R}^{N\times N}$ is the first-order weight matrix, focusing on the pairwise relationship. High-order graphs adhere to the principle that the neighbor of a neighbor is also a neighbor, which excavates important structural information that is not easy to observe in a first-order graph. Given the first-order graph $\mathbf{W}_p$, the $h$th-order graph is defined as
\begin{equation}\label{high}
		\begin{aligned}
\mathbf{W}_p^h=
\begin{cases}
\mathbf{W}_p, &h=1,\\
\mathbf{W}_p^{h-1}\mathbf{W}_p, &h>1.
\end{cases} 	\end{aligned}
	\end{equation}
Moreover, to address the dilemma of order selection, multi-order graphs introduce weights to construct the most consistent graph, which is specifically defined as
\begin{equation}\label{high2}
		\begin{aligned}
\mathbf{W}_p^l=\sum_{i=1}^{l}q^{i}\mathbf{W}_p^{i},	\end{aligned}
	\end{equation}
where $q^{i}\in[0,1)$, $\sum_{i=1}^{l}q^{i}=1$, and $l$ is the maximum order.

  Next, we denote the correlation tensor as
		\begin{equation}
			\begin{aligned}
				\mathcal{P}&=\mathcal{C}_{12 \cdots m} \times_{1} \mathbf{H}_{1}^\top \times_{2}  \cdots \times_{m} \mathbf{H}_{m}^\top \in \mathbb{R}^{r\times \cdots \times r}.
			\end{aligned}
	\end{equation}
	Building upon this, to jointly learn the unique characteristics of individual views and their shared representation, we propose 
	\begin{equation}\label{21otcca}
		\begin{aligned}
			\underset{\{\mathbf{H}_{p}\}}{\operatorname{min}} &\quad -\frac{1}{2}\|\mathcal{P}\|_\textrm{F}^2 +\sum_{p=1}^{m}\lambda_{p}\|\mathbf{H}_p\|_{2,1}+\sum_{p=1}^{m}\textrm{Tr}(\mathbf{Z}_{p}^\top\mathbf{L}_{p}^{l}\mathbf{Z}_{p})\\
			\text { s.t. } &~~~ \mathbf{H}_{p}^\top  \mathbf{X}_p \mathbf{X}_p^\top\mathbf{H}_{p}=\mathbf{I}_{r}, p\in [m],
		\end{aligned}
	\end{equation}
	where $\mathbf{L}_{p}^{l}=\mathbf{S}_{p}^{l}-\mathbf{W}_{p}^{l}$, ${\mathbf{S}_{p}^{l}}_{ii}=\sum_{j}{\mathbf{S}_{p}^{l}}_{ij}$, and $\mathbf{W}_{p}^{l}$ represents the multi-order graph. 

In this paper, we refer to \eqref{21otcca} as STCCA-L. It can be seen that compared to \eqref{tcca-o}, STCCA-L directly applies structural sparse regularization to the projection matrix $\mathbf{H}_p$ to reduce feature redundancy. This design not only reduces the redundancy of the learning subspace but also promotes feature selection, thereby leading to a more interpretable and compact representation. Furthermore, STCCA-L incorporates the graph Laplacian regularization term $\textrm{Tr}(\mathbf{Z}_p^\top \mathbf{L}_p^l \mathbf{Z}_p)$ into the objective. This enables STCCA-L to explicitly retain the inherent local geometry of each view during the feature extraction.

  Following Definition \ref{ste}, the proposed model in \eqref{21otcca} can be rewritten as the manifold form
		\begin{equation}\label{stcca-l}
		\begin{aligned}
			\underset{\{\mathbf{H}_{p}\}}{\operatorname{min}} &\quad -\frac{1}{2}\|\mathcal{P}\|_\textrm{F}^2 +\sum_{p=1}^{m}\lambda_{p}\|\mathbf{H}_p\|_{2,1}+\sum_{p=1}^{m}\textrm{Tr}(\mathbf{Z}_{p}^\top\mathbf{L}_{p}^{l}\mathbf{Z}_{p})\\
			\text { s.t. } &~~~  \mathbf{X}_p^\top\mathbf{H}_{p} \in\textrm{St($r,N$)},  p\in [m].
		\end{aligned}
	\end{equation}
	
    Rewriting problem \eqref{21otcca} in the manifold form in \eqref{stcca-l} is essential, as it explicitly encodes the orthogonality constraints within the Stiefel manifold structure. This reformulation not only provides a geometrically consistent framework for modeling the projection matrices but also facilitates theoretical analysis and enables the use of manifold optimization tools in a principled manner.

	\subsection{Optimization}
    To solve problem \eqref{stcca-l}, tensor decomposition techniques can be incorporated into optimization strategies such as the alternating direction method of multipliers (ADMM) \cite{Xu2025} and gradient descent \cite{Yi2020}. However, the inclusion of regularization terms often leads to increased computational complexity, making most existing algorithms time-consuming. Considering the presence of the non-smooth $\ell_{2,1}$-norm, we employ the proximal gradient (PG) method for solving the problem. Meanwhile, the constraint involving the Stiefel manifold poses additional challenges for direct optimization. Therefore, for the problem defined in \eqref{stcca-l}, we design an algorithm based on alternating manifold PG, as detailed below.

    \subsubsection{Reformulation with Auxiliary Variables}
    For the convenience of notation, denote
    \begin{equation}		\begin{aligned}F(\{\mathbf{H}_p\})=-\frac{1}{2}\|\mathcal{P}\|_\textrm{F}^2+\sum_{p=1}^{m}\textrm{Tr}(\mathbf{Z}_{p}^\top\mathbf{L}_{p}^{l}\mathbf{Z}_{p}), \end{aligned}
    	\end{equation} and introduce the auxiliary variable $\mathbf{Y}_{p}=\mathbf{H}_p$, then problem \eqref{stcca-l} can be reformulated as
    \begin{equation}\label{stcca-l1}
		\begin{aligned}
			\underset{\{\mathbf{H}_{p},\mathbf{Y}_p\}}{\operatorname{min}}& \quad F(\{\mathbf{H}_p\}) +\sum_{p=1}^{m}\lambda_{p}\|\mathbf{Y}_p\|_{2,1}\\
			\text { s.t. } ~~&~~~  \mathbf{X}_p^\top\mathbf{H}_{p} \in\textrm{St($r,N$)},  p\in [m].
		\end{aligned}
	\end{equation}

\subsubsection{Proximal Gradient Step on the Stiefel Manifold}
For the Stiefel manifold, it needs to ensure that the descent direction lies in the tangent space. For brevity, it only focuses on the subproblems of the $p$th view $\mathbf{H}_{p}$. This motivates the following subproblem for finding the descent direction $\mathbf{D}_p^k$ in the $k$th iteration \cite{chen2020alternating}, which is 
\begin{equation}\label{smtcca0}
		\begin{aligned}
	    	\underset{\mathbf{D}_{p}}{\operatorname{min}} \quad &\langle \textrm{grad } F(\mathbf{H}_p^{k}), \mathbf{D}_p\rangle+\frac{1}{2t}\|\mathbf{D}_p\|_{\textrm{F}}^2\\
            & +\lambda_{p}\|\mathbf{H}_p^k+\mathbf{D}_p\|_{2,1}\\
			\text { s.t. }~~& \mathbf{D}_p \in \textrm{T}_{\mathbf{H}_p^k}\textrm{St($r,N$)},
		\end{aligned}
	\end{equation}
where $\textrm{T}_{\mathbf{H}_p^k} \textrm{St($r,N$)}=\{\mathbf{D}_p\mid \mathbf{D}_p^\top\mathbf{X}_p\mathbf{X}_p^\top \mathbf{H}_p+\mathbf{H}_p^\top\mathbf{X}_p\mathbf{X}_p^\top \mathbf{D}_p=0\}$ is the  tangent space of the Stiefel manifold $\textrm{St($r,N$)}$.
	According to the definition of Riemannian gradient, for any $\mathbf{D}_p \in \textrm{T}_{\mathbf{H}_p^k}\textrm{St($r,N$)}$, it has
	\begin{equation}\label{gradf}
		\langle \textrm{grad }F(\mathbf{H}_p^k),\mathbf{D}_p\rangle=\langle\nabla F(\mathbf{H}_p^k),\mathbf{D}_p\rangle.
	\end{equation}
Recall that $F$ consists of the tensor norm and a trace regularization, then $\nabla F(\mathbf{H}_p^k)$ can be decomposed accordingly. For the tensor part, it fixes all projection matrices except for the $p$th and computes the derivative along mode-$p$ as
\begin{equation}\begin{aligned}-\mathcal{C}_{12\cdots m}&\times_{1}\mathbf{H}_1^{k\top}\times_{2}\cdots\times_{p-1}\mathbf{H}_{p-1}^{k\top}\times_{p+1}\cdots\\ &\times_{m}\mathbf{H}_m^{k\top} \in \mathbb{R}^{r \times \cdots \times d_p\cdots \times r}, \end{aligned}\end{equation}
and reshaped into $m-1$ matrices $\mathbf{C_p}_i \in \mathbb{R}^{d_p \times r},  i\in [m-1]$. Combining both terms, $\nabla F$  is given by
\begin{equation}\nabla F(\mathbf{H}_p)=\sum_{i=1}^{m-1}(\mathbf{C_p}_i)+\mathbf{X}_p\mathbf{L}_p^{l}\mathbf{Z}_p^k.\end{equation}

Define the linear operator \begin{equation}\begin{aligned}A^k(\mathbf{D}_p)=\mathbf{D}_p^\top\mathbf{X}_p\mathbf{X}_p^\top \mathbf{H}_p+\mathbf{H}_p^\top\mathbf{X}_p\mathbf{X}_p^\top \mathbf{D}_p,\end{aligned}\end{equation} 
then problem \eqref{smtcca0} can be reformulated as
		\begin{equation}\label{smtcca}
		\begin{aligned}
	    	\underset{\mathbf{D}_{p}}{\operatorname{min}} \quad &\langle \nabla F, \mathbf{D}_p\rangle+\frac{1}{2t}\|\mathbf{D}_p\|_{\textrm{F}}^2+\lambda_{p}\|\mathbf{H}_p^k+\mathbf{D}_p\|_{2,1}\\
			\text { s.t. }~~& A^k(\mathbf{D}_p)=0,
		\end{aligned}
	\end{equation}
    where the PG step is restricted to the tangent space of the Stiefel manifold. Once the descent direction $\mathbf{D}_p^k$ is obtained by solving \eqref{smtcca}, an Armijo-type line search is employed to determine the step size $\alpha^k$. The update is then projected back onto the Stiefel manifold via
    \begin{equation}
\begin{aligned}\mathbf{H}_p^{k+1}=\operatorname{Retr}_{\mathbf{H}_p^k}(\alpha^k\mathbf{D}_p^k),\end{aligned}
	\end{equation}
ensuring feasibility under the manifold constraint.

\subsubsection{Efficient Solution via the Semi-Smooth Newton Method}

    The next important question is how to solve problem \eqref{smtcca} quickly? The semi-smooth Newton (SSN) method \cite{xiao2018regularized} has recently attracted considerable attention for its efficiency and accuracy in solving structured convex problems. It has been successfully applied across a variety of domains, including LASSO \cite{li2018highly} and sparse principal component analysis \cite{chen2020alternating}. In this regard, we attempt to develop an efficient SSN method.

   The Lagrangian function of problem \eqref{smtcca}  can be written as
   	\begin{equation}\label{larall}
   	\begin{aligned}
   		\mathcal{L} (\mathbf{D}_p;\mathbf{\Lambda}_p)&=\langle \nabla F(\mathbf{H}_p^k),\mathbf{D}_p \rangle+\lambda_{p}\|\mathbf{H}_p^k+\mathbf{D}_p\|_{2,1}\\
   		&+\frac{1}{2t}\|\mathbf{D}_p\|_\textrm{F}^{2}-\langle A^k(\mathbf{D}_p),\mathbf{\Lambda}_p\rangle,\\
   	\end{aligned}
   \end{equation}
   where $\mathbf{\Lambda}_p,  p\in [m]$ are the Lagrangian multipliers. We analyze the solution in four steps.

   Firstly, it constructs the Karush-Kuhn-Tucker (KKT) condition of problem \eqref{smtcca} as
 	\begin{equation}\label{kkt1}
   0\in \partial_{\mathbf{D}_p}\mathcal{L}(\mathbf{D}_p;\mathbf{\Lambda}_p),~A^k(\mathbf{D}_p)=0.
    \end{equation}
   The first condition leads to the proximal mapping
 \begin{equation}\label{skkt1}
   \mathbf{D}_p=\operatorname{prox}_{2,1}( \mathbf{B}(\mathbf{\Lambda}_{p}),t)-\mathbf{H}_p^k,
   \end{equation}
   where $\mathbf{B}(\mathbf{\Lambda}_{p})=\mathbf{H}_p^k-t(\nabla F(\mathbf{H}_p^k)-2\mathbf{X}_p\mathbf{X}_p^\top\mathbf{H}_p^k\mathbf{\Lambda}_p)$. Substituting \eqref{skkt1} into the second condition of \eqref{kkt1} derives
    \begin{equation}\label{skkt2}
Q(\mathbf{\Lambda}_p)=\mathbf{D}_p^\top\mathbf{X}_p\mathbf{X}_p^\top\mathbf{H}_p^k+\mathbf{H}_p^{k\top}\mathbf{X}_p\mathbf{X}_p^\top\mathbf{D}_p=0.
   \end{equation}

Secondly, the operator $Q$ is monotone and Lipschitz continuous \cite{chen2020proximal}, which makes it suitable for the SSN method. To proceed, we compute the generalized Jacobian of $Q$.
The vectorization of $Q(\mathbf{\Lambda}_{p})$ can be  showed as
   \begin{equation}
   	 \begin{aligned}
   	 	&\operatorname{vec}(Q(\mathbf{\Lambda}_p))=(\mathbf{K}_{r r}+\mathbf{I}_{r^2})(\mathbf{H}_p^{k\top}\mathbf{X}_p\mathbf{X}_p^\top\otimes\mathbf{I}_r)\\
   &[\operatorname{prox}_{2,1}( \operatorname{vec}( \mathbf{H}_p^{k\top}\mathbf{X}_p\mathbf{X}_p^\top)-t \nabla F(\mathbf{H}_p^k),t)]\\
   	 	&+2 t(\mathbf{X}_p\mathbf{X}_p^\top\mathbf{H}_p^k\otimes\mathbf{I}_r) \operatorname{vec}(\mathbf{\Lambda}_p)-\operatorname{vec}(\mathbf{H}_p^{k\top}),
   	 \end{aligned}
   \end{equation}
   where $\mathbf{K}_{r d_p}$ and $\mathbf{K}_{rr}$ are the commutation matrices. Define
 \begin{equation}
   	\mathbf{\Xi_p}_j=\begin{cases}
   		 \mathbf{I}_r-\frac{\tau_1 t}{\|\mathbf{b}_{j}\|_2}\mathbf{R}, & \text { if }\|\mathbf{b}_{j}\|_2>t \tau_1, \\ \gamma \frac{\mathbf{b}_{j} \mathbf{b}_{j}^{\top}}{(t \tau_1)^2}, & \text { if }\|\mathbf{b}_{j}\|_2=t \tau_1, \\ 0, & \text { otherwise},
   	\end{cases}
   \end{equation}
where $p\in [m]$, $j\in [d_p$], $\mathbf{R}=( \mathbf{I}_r-\frac{\mathbf{b}_{j} \mathbf{b}_{j}^{\top}}{\|\mathbf{b}_{j}\|_2^2}),$ $\gamma \in[0,1],$ and 
$\mathbf{b}_{j}$ is the $j$th column of $\mathbf{B}( \mathbf{\Lambda}_p)^\top$.
Let the generalized Jacobian be \begin{equation}
   	 \begin{aligned}\mathcal{J}(\mathbf{y})|_{\mathbf{y}=\operatorname{vec}(\mathbf{B}( \mathbf{\Lambda}_{p})^\top)}=\textrm{Diag}( \mathbf{\Xi_p}_1,\cdots, \mathbf{\Xi_p}_{d_p}). \end{aligned}
   \end{equation} Then the generalized Jacobian matrix $\mathbf{V}$ of $\operatorname{vec}(Q( \mathbf{\Lambda}_p))$ is
   \begin{equation}
    \begin{aligned} \mathbf{V}&=2t(\mathbf{K}_{rr}+\mathbf{I}_{r^2})(\mathbf{H}_p^{k\top}\mathbf{X}_p\mathbf{X}_p^\top\otimes\mathbf{I}_r)\\
    &\quad\mathcal{J}(\mathbf{y})(\mathbf{X}_p\mathbf{X}_p^\top\mathbf{H}_p^{k}\otimes\mathbf{I}_r).
    \end{aligned}
   \end{equation}
   By monotonicity of $Q$, it is seen that $\mathbf{V}$ is positive semi-definite \cite{xiao2018regularized} and serves as a valid surrogate of the true Jacobian. For $\text{any}~\sigma \in \mathbb{R}^{r^2}$, it has 
\begin{equation}
\begin{aligned}\mathbf{V}\sigma=\nabla(\operatorname{vec}(Q(\operatorname{vec}(\mathbf{\Lambda}_p))))\sigma.\end{aligned}
   \end{equation}

 Thirdly, as $\mathbf{\Lambda}_{p}$ is symmetric, it uses 
$\overline{\operatorname{vec}}(\mathbf{\Lambda}_p)$ to denote the 
$\frac{1}{2} r (r+1)$-dimensional vector obtained from $\operatorname{vec}(\mathbf{\mathbf{\Lambda}_{p}})$ by eliminating all superdiagonal elements of $\mathbf{\Lambda}$. Using the duplication matrix $\mathbf{U}_p \in \mathbb{R}^{r^2\times\frac{1}{2}r(r+1)}$ and its Moore-Penrose inverse $\mathbf{U}_p^+$, it has
  \begin{equation}
 \mathbf{U}_p\overline{\operatorname{vec}}(\mathbf{\Lambda}_{p})=\operatorname{vec}(\mathbf{\Lambda}_{p}),
  \end{equation}
 and the generalized Jacobian in the reduced space is
 \begin{equation}
 	{V}(\overline{\operatorname{vec}}(\mathbf{\Lambda}_{p}))=t\mathbf{U}_p^+\mathbf{V}\mathbf{U}_p.
 \end{equation}
 Then, the SSN update direction $\mathbf{d}_k$ is computed by solving the linear system
 \begin{equation}
 	(\mathbf{V}+\eta \mathbf{I}_{r^2})\mathbf{d}=-\overline{\operatorname{vec}}(Q(\overline{\operatorname{vec}}(\mathbf{\Lambda}_{p}^k))),
 \end{equation}
where $\eta > 0$.

Finally, the update rule of $\mathbf{\Lambda}_{p}^k$ is
\begin{equation}
	\overline{\operatorname{vec}}(\mathbf{\Lambda}_{p}^{k+1})=\overline{\operatorname{vec}}(\mathbf{\Lambda}_{p}^k)+\mathbf{d}_k.
\end{equation}

In summary, the full implementation details are provided in Algorithm~\ref{al1}.

\begin{algorithm}[t]
	\caption{Optimization algorithm for solving \eqref{stcca-l}}\label{al1}
	\KwIn{Multi-view data $\mathbf{X}=[\mathbf{X}_1,\cdots,\mathbf{X}_m]$, where $\mathbf{X}_p \in\mathbb{R}^{d_{p} \times N}$, $p\in[m]$, step-size $t$, maximum number of iterations $T$, and $\gamma$ $\in$ (0,1). Calculate covariance tensor $\mathcal{C}_{12 \cdots m}$, and initialize $\mathbf{H}_p^0 \in \textrm{St}(n,r)$}
	\KwOut{$\{\mathbf{H}_p^k\}$ }
	\For{$p\in[m]$}{
		\If{$k<T$}
		{
			Obtain $\mathbf{D}_p^k$ via \eqref{smtcca} using the SSN method\\
			\While{$F(\operatorname{Retr}_{\mathbf{H}_p^k}(\alpha\mathbf{D}_p^k))\geq F(\mathbf{H}_p^k)-\frac{\alpha\|\mathbf{D}_p^k\|_{\textrm{F}}^2}{2t}$}{
				$\alpha=\gamma \alpha$}
			Set $\mathbf{H}_p^{k+1}=\operatorname{Retr}_{\alpha\mathbf{H}_p^k}(\alpha\mathbf{D}_p^k)$ }
	}
\end{algorithm}

\subsection{Convergence Analysis}

Despite the empirical success of existing TCCA methods \cite{luo2015tensor,sun2023learning}, they lack rigorous convergence guarantees.  In what follows, we provide a detailed convergence analysis of the proposed algorithm to ensure its theoretical soundness.

It denotes by $\mathbf{H}_{[p]}^{k}(\alpha)$ the collection of projection matrices at iteration $k$, \textit{i.e.},
\begin{equation}
\mathbf{H}_{(p)}^{k}(\alpha) = \{\mathbf{H}_1^{k}, \cdots, \mathbf{H}_{p-1}^{k}, \mathbf{H}_p^{k} + \alpha \mathbf{D}_p^{k}, \mathbf{H}_{p+1}^{k}, \cdots, \mathbf{H}_m^{k}\}.
\end{equation}
Define the objective function of problem \eqref{smtcca} as \begin{equation}g(\mathbf{D}_p)=\langle \nabla F, \mathbf{D}_p\rangle+\frac{1}{2t}\|\mathbf{D}_p\|_{\textrm{F}}^2+\lambda_{p}\|\mathbf{H}_p^k+\mathbf{D}_p\|_{2,1} . \end{equation}

Now, we prove that $\mathbf{D}_p^k$ is a descending direction in the tangent space.
\begin{Lemma}\label{lem1}
 For any $\alpha\in [0,1],$ if $t\le \frac{1}{L_{p}}$, where $L_{p}$ is the Lipschitz constant of $\nabla_{\mathbf{H}_p} F$, the following inequality holds
\begin{equation}F(\mathbf{H}_{(p)}^{k}(\alpha))+\|\mathbf{H}_p^{k} +\alpha \mathbf{D}_p^{k}\|_{2,1}\le F(\mathbf{H}_{(p)}^{k}(0))+ \|\mathbf{H}_p^{k}\|_{2,1}.\end{equation}
\end{Lemma}
\begin{IEEEproof}
 Since the
objective function $g(\mathbf{D}_p)$ is $\frac{1}{t}$-strongly convex, for $\widehat{\mathbf{D}}_p,\mathbf{D}_p$, it has
\begin{equation}\label{c5}
			g(\widehat{\mathbf{D}}_p)\ge g(\mathbf{D}_p)+\langle\partial g(\mathbf{D}_p),\widehat{\mathbf{D}}_p-\mathbf{D}_p\rangle+\frac{\alpha}{2}\|\widehat{\mathbf{D}}_p-\mathbf{D}_p\|_{\textrm{F}}^{2}.
		\end{equation}
Specifically, if $ \widehat{\mathbf{D}}_p, \mathbf{D}_p \in \textrm{T}_\mathbf{\mathbf{H}_{p}^k}\textrm{St($N,r$)},$ then it has
\begin{equation}
\langle \partial g(\mathbf{D}_p), \widehat{\mathbf{D}}_p - \mathbf{D}_p \rangle = \langle \operatorname{proj}_{\textrm{T}_{\mathbf{H}_{p}^k}} ( \partial g(\mathbf{D}_p) ), \widehat{\mathbf{D}}_p - \mathbf{D}_p \rangle.
\end{equation}
From the Riemannian optimality condition, it follows
\begin{equation}
0 \in \operatorname{proj}_{\textrm{T}_{\mathbf{H}_{p}^k}}  ( \partial g(\mathbf{D}_p^k)  ).
\end{equation}
Letting $ \mathbf{D}_p = \mathbf{D}_p^k $, $ \widehat{\mathbf{D}}_p = \alpha \mathbf{D}_p^k $, and $ \alpha \in [0, 1] $ in \eqref{c5}, it yields
\begin{equation}
    	 g (\alpha \mathbf{D}_p^k )-g(\mathbf{D}_p^k) \ge \frac{(1-\alpha)^2}{2 t}\|\mathbf{D}_p^k\|_\textrm{F}^2.
    \end{equation}
This, together with the definition of $g$ and the convexity of the $\ell_{2,1}$-norm, implies that
\begin{equation}\begin{aligned}
    	&(1-\alpha)\langle \nabla F(\mathbf{H}_{(p)}^{k}(0)),\mathbf{D}_p^k\rangle +\frac{1-\alpha}{t}\|\mathbf{D}_p^k\|_\textrm{F}^2\\
    &+(1-\alpha)(\|\mathbf{H}_p^{k}+\mathbf{D}_p^{k}\|_{2,1}-\|\mathbf{H}_p^{k}\|_{2,1})\le 0.
    \end{aligned}\end{equation}
 Combining the convexity of $\ell_{2,1}$-norm and the Lipschitz continuity of $g$, it has
 \begin{equation}\begin{aligned}
 &F(\mathbf{H}_{(p)}^{k}(\alpha))-F(\mathbf{H}_{(p)}^{k}(0))\\
  &+\|\mathbf{H}_p^{k}+\alpha\mathbf{D}_p^{k}\|_{2,1}-\|\mathbf{H}_p^{k}\|_{2,1}\\
 &\le\alpha \langle \nabla F(\mathbf{H}_{(p)}^{k}(0)),\mathbf{D}_p^k\rangle+\frac{\alpha^2}{2t}\|\mathbf{D}_p^k\|_\textrm{F}^2\\
 &+\alpha(\|\mathbf{H}_p^{k}+\mathbf{D}_p^{k}\|_{2,1}-\|\mathbf{H}_p^{k}\|_{2,1})\\
 &\le-\frac{\alpha}{2t}\|\mathbf{D}_p^k\|_\textrm{F}^2.
    \end{aligned}\end{equation}
Thus, the proof is completed.
\end{IEEEproof}

    Furthermore, the following lemma shows that when $\mathbf{D}_p^k= 0,  p\in [m]$, then a stationary point is found. Specifically, a point $\mathbf{H} \in \textrm{St($r,N$)}$ is referred to as a stationary point of problem \eqref{stcca-l} if it satisfies the first-order optimization condition.
   \begin{Lemma}\label{lem2}
   	If the sequence $\{\mathbf{D}_p^k\}$ satisfies $\mathbf{D}_p^k = 0$ for all $ p\in [m]$, then the sequence $\{\mathbf{H}_p^k\}$ is a stationary point of problem~\eqref{stcca-l}.
   \end{Lemma}
   \begin{IEEEproof}
   	For any $p\in [m],$ the optimality conditions of problem \eqref{smtcca} can be written as
   	\begin{equation}
   	0 \in \frac{1}{t} \mathbf{D}_p^k+\nabla F(\mathbf{H}_p^k)+\operatorname{proj}_{\textrm{T}_{\mathbf{H}_{p}^k}} \partial \|\mathbf{H}_p^k+\mathbf{D}_p^k\|_{2,1},
   	\end{equation}where $\mathbf{D}_p^k \in \textrm{T}_{\mathbf{H}_{p}^k}\textrm{St($r,N$)}$. If $\mathbf{D}_p^k=0$, then we have \begin{equation}0 \in \nabla F(\mathbf{H}_p^k)+\operatorname{proj}_{\mathrm{T}_{\mathbf{H}_p^k}} \partial \|\mathbf{H}_p^k+\mathbf{D}_p^k\|_{2,1}.
   \end{equation} 
   It is the first-order necessary condition of problem \eqref{stcca-l}.
   	\end{IEEEproof}
    
   Define the objective function of \eqref{stcca-l} as
\begin{equation}G(\mathbf{H}_p)=F(\mathbf{H}_p)+\|\mathbf{H}_p+\mathbf{D}_p\|_{2,1}.\end{equation}

\begin{Lemma}
Assume that $\{\mathbf{H}_p^{k}\}$ is generated by Algorithm \ref{al1}, then $\{G(\mathbf{H}_p^{k})\}$ is monotonically decreasing. And it satisfies the following inequality
   \begin{equation}\label{Dp}
   	G(\operatorname{Retr}_{\mathbf{H}_p^k}(\alpha \mathbf{D}_p^k))-G(\mathbf{H}_p^k) \leq-\frac{\alpha}{2 t}\|\mathbf{D}_p^k\|_{\textrm{F}}^2,~p\in [m].
   \end{equation}
   \end{Lemma}
   \begin{IEEEproof}
   Let $\mathbf{H}_p^{k+}=\mathbf{H}_p^k+\alpha \mathbf{D}_p^k$. Following \cite{boumal2019global} and the $L$-Lipschitz continuity of $\nabla G$, for any $\alpha>0$, we have
  \begin{equation}
   	  \begin{aligned}
   		&G(\operatorname{Retr}_{ \mathbf{H}_p^k}(\alpha  \mathbf{D}_p^k))-G( \mathbf{H}_p^k)\\
   &\leq \langle \nabla G( \mathbf{H}_p^k),\operatorname{Retr}_{ \mathbf{H}_p^k}(\alpha  \mathbf{D}_p^k)-\mathbf{H}_p^{k+}+\mathbf{H}_p^{k+}-\mathbf{H}_p^{k}\rangle\\
   &+\frac{L}{2}\|\operatorname{Retr}_{ \mathbf{H}_p^k}(\alpha  \mathbf{D}_p^k)-\mathbf{H}_p^k\|_{\mathrm{F}}^2\\
   &\leq \zeta_2\|\nabla G( \mathbf{H}_p^k)\|_{\mathrm{F}}\|\alpha  \mathbf{D}_p^k\|_{\mathrm{F}}^2\\
   &+\alpha\langle\nabla G( \mathbf{H}_p^k),  \mathbf{D}_p^k\rangle+\frac{{L} \zeta_1^2}{2}\|\alpha  \mathbf{D}_p^k\|_{\mathrm{F}}^2.
   	\end{aligned}
 \end{equation}
   Since $\nabla G$ is continuous on the compact manifold $\textrm{St($r,N$)}$, there exists a constant $\mu>0$ such that $\|\nabla G(\mathbf{H}_p^k)\|_{\mathrm{F}} \leq \mu$, any $\mathbf{H}_p \in \textrm{St($r,N$)}$. It has
   	\begin{equation}
     \begin{aligned}
   	&G(\operatorname{Retr}_{\mathbf{H}_p^k}(\alpha \mathbf{D}_p^k))-G(\mathbf{H}_p^k) \\
    &\leq \alpha\langle\nabla G(\mathbf{H}_p^k), \mathbf{D}_p^k\rangle+{c}_0 \alpha^2\|\mathbf{D}_p^k\|_{\mathrm{F}}^2
    \end{aligned}
   \end{equation}
   	where ${c}_0=\zeta_2 \mu+L \zeta_1^2 / 2$. This implies that
   \begin{equation}
   	\begin{aligned}
   		&G(\operatorname{Retr}_{\mathbf{H}_p^k}(\alpha  \mathbf{D}_p^k))-G( \mathbf{H}_p^k) \\
       &\leq({c}_0+\delta \zeta_2-\frac{1}{\alpha t})\|\alpha  \mathbf{D}_p^k\|_{\mathrm{F}}^2,
   	\end{aligned}
   	\end{equation}
   	where $\delta$ is the Lipschitz continuity of $\ell_{2,1}$-norm. Upon setting $\bar{\alpha}=1 /(2({c}_0+\delta \zeta_2) t)$, for any $0<\alpha \leq \min \{\bar{\alpha}, 1\}$, it holds
   	\begin{equation}
    \begin{aligned}
   	&G(\operatorname{Retr}_{\mathbf{H}_p^k}(\alpha  \mathbf{D}_p^k))-G( \mathbf{H}_p^k) \\
      &\leq-\frac{1}{2 \alpha t}\|\alpha  \mathbf{D}_p^k\|_{\mathrm{F}}^2=-\frac{\alpha}{2 t}\| \mathbf{D}_p^k\|_{\mathrm{F}}^2.
    \end{aligned}
   	\end{equation}
    Therefore, after applying a retraction to $\mathbf{D}_p^k$, it is also a descending direction of the objective function in \eqref{stcca-l}. The proof is completed.
   \end{IEEEproof}
   
To end this section, we theoretically establish the global convergence of Algorithm \ref{al1} to a stationary point.  
 
   \begin{Theorem}\label{them}
 The sequence $\{\mathbf{H}_p^k\}$ generated by Algorithm \ref{al1} converges to
a stationary point of problem \eqref{stcca-l}.
   \end{Theorem}
  \begin{IEEEproof}
  	Since $G$ is bounded below on $\textrm{St($r,N$)}$, by \eqref{Dp}, it is not hard to obtain 
    \begin{equation}\underset{k \rightarrow \infty}{\operatorname{\lim}}\|\mathbf{D}_p^k\|_{\textrm{F}}^2=0. \end{equation} Combining with Lemma \ref{lem2}, it follows that every limit point of $\{\mathbf{H}_p^k\}$ is a stationary point of \eqref{stcca-l}.
   \end{IEEEproof}

\subsection{Complexity Analysis}   	
For our proposed Algorithm~\ref{al1}, the overall computational complexity is 
$O(Tm(r^m + d_{a} N + d_{a}^2 r))$, 
where $T$ is the number of outer iterations, $m$ is the number of views, and $d_{a} = \max_p d_p$. 
The main computational cost arises from constructing the covariance tensor, solving the SSN subproblem, evaluating the objective function, and performing retraction onto the Stiefel manifold. The runtime comparison will be provided in the following section.

\begin{table}[t]
    \renewcommand\arraystretch{1.2} 
    \caption{Statistics of all selected datasets.}\label{data} 
    \centering
    \setlength{\tabcolsep}{1.0mm} 
    \vspace{-0.1cm}
    \begin{tabular}{|c|c|c|c|c|c|}
      \hline
     ~~Sizes~~& ~~~{Datasets}~~~ & ~{{Instances}}~ & ~{{Clusters}}~ &~~~~{{Views}}~~~~& ~~{{Dim}}~~  \\
      \hline \hline
      \multirow{10}{*}{{Small}}&\multirow{3}{*}{{3Sources}} &\multirow{3}{*}{169}&\multirow{3}{*}{6}& Reuters&3068\\ \cline{5-6}
			~&~&~&~& BBC&3560 \\ \cline{5-6}
			~&~&~&~& Guardian&3631  \\
             \cline{2-6}
      ~&\multirow{5}{*}{{ MSRC}} &\multirow{5}{*}{210}&\multirow{5}{*}{7}& CN&24\\ \cline{5-6}
			~&~&~&~& HOG&576 \\ \cline{5-6}
			~&~&~&~& GIST&512\\ \cline{5-6}
            ~&~&~&~& LBP&256 \\ \cline{5-6}
			~&~&~&~& CENT&256\\
            \cline{2-6}
       ~&\multirow{2}{*}{{ BBCsport}} &\multirow{2}{*}{544}&\multirow{2}{*}{5}& View1&3183\\ \cline{5-6}
			~&~&~&~& View2&3203 \\
           \hline
         \multirow{14}{*}{{Medium}}& \multirow{5}{*}{{ Reuters}} &\multirow{5}{*}{1200}&\multirow{5}{*}{6}& English&2000\\ \cline{5-6}
			~&~&~&~& French&2000\\ \cline{5-6}
			~&~&~&~& German&2000\\ \cline{5-6}
            ~&~&~&~& Spanish&2000\\ \cline{5-6}
             ~&~&~&~& Italian&2000\\
             \cline{2-6}
         ~&\multirow{4}{*}{{Caltech101}}&\multirow{4}{*}{1474}&\multirow{4}{*}{7}&Gabor &48\\ \cline{5-6}
			~&~&~&~&WM&40	 \\ \cline{5-6}
			~&~&~&~&CENT&254\\ \cline{5-6}
			~&~&~&~&HOG&1984 \\
           \cline{2-6}
         ~&\multirow{5}{*}{{Handwritten}} &\multirow{5}{*}{2000}&\multirow{5}{*}{10}& FOU&76\\ \cline{5-6}
			~&~&~&~& FAC&216 \\ \cline{5-6}
			~&~&~&~& KAR&64  \\ \cline{5-6}
            ~&~&~&~& PIX&240 \\ \cline{5-6}
			~&~&~&~& ZER&47  \\
            \hline
        \multirow{6}{*}{{Large}}~&\multirow{2}{*}{{ MNIST}} &\multirow{2}{*}{10000}&\multirow{2}{*}{10}& ISO&30\\ \cline{5-6}
			~&~&~&~& NPE&30  \\
            \cline{2-6}
			~&\multirow{4}{*}{{Animal}} &\multirow{4}{*}{11673}&\multirow{4}{*}{20}& View 1&2688\\  \cline{5-6}
			~&~&~&~& View 2&2000\\ \cline{5-6}
			~&~&~&~& View 3&2001  \\ \cline{5-6}
            ~&~&~&~& View 4&2000 \\
      \hline
    \end{tabular}
\end{table}

	\section{Numerical Experiments}\label{sec4}
	In this section, we evaluate the performance of our proposed STCCA-L with various competitive methods on eight well-known multi-view datasets, covering 3Sources\footnote{http://mlg.ucd.ie/datasets/3sources.html}, MSRC\footnote{https://github.com/zhudafa/Multi-view-datasets\label{MLdates}}, BBCsport\footnote{http://mlg.ucd.ie/datasets/bbc.html}, Reusters\textsuperscript{\ref{MLdates}},  Caltech101\footnote{https://data.caltech.edu/records/mzrjq-6wc02}, Handwritten\footnote{https://github.com/cvdfoundation/mnist}, MNIST\footnote{https://tensorflow.google.cn/datasets/catalog/mnist}, and Animal\textsuperscript{\ref{MLdates}}. These datesets can be divided into three groups based on the sample size, as shown in Table \ref{data}, where small, medium, and large respectively represent sample sizes less than 1000, greater than 1000 and less than 10000, and greater than 10000.

Section \ref{A} gives the experimental settings, Section \ref{B} provides the experimental results, Section \ref{C} presents the ablation studies, including the contribution of each group, the influence of graph construction, and the effect of algorithm initialization, and Section \ref{E} discusses the robustness, parameter analysis, stability, and efficiency.

\subsection{Experimental Settings}\label{A}

\subsubsection{Implementation Details}
First, it applies the principal component analysis to all the datasets, reducing the dimension from 2 to 20 at intervals of 2. Then, for $p\in [m]$, compute the $p$th view projection matrix as $\mathbf{Z}_p = \mathbf{X}_p^\top \mathbf{H}_p$. The final representation is obtained by concatenating all view projections, \textit{i.e.},  $\mathbf{Z}=[\mathbf{Z}_1,\cdots,\mathbf{Z}_m] \in \mathbb{R}^{N \times mr}$. Finally, the K-nearest neighbor (KNN) classifier is used in our experiments to measure classification performance. For fair comparison, the number of neighbors K is adjusted according to dataset characteristics, while the same K is used across all competing methods on each dataset. Moreover, it selects the adaptive neighbor graph method initial weight matrix $\mathbf{W}_p$. Different from the traditional KNN graph with fixed neighbor weights, this method learns the neighbor weights of each sample by minimizing the adaptive weight distribution under the constraint of reconstruction error, thereby automatically determining the optimal local structure. Each penalty parameter is determined using cross-validation techniques, and the test ratio is set to $0.3$. The mean accuracy values and related standard deviations are also recorded after each experiment is randomly repeated 10 times.

\begin{table*}[t]
    \renewcommand\arraystretch{1.2} 
    \caption{Classification accuracy ($\%$) of all compared methods under the best dimensions. }\label{acc} 
    \centering
    	 \setlength{\tabcolsep}{1.0mm}
    \vspace{-0.1cm}
        \begin{threeparttable}
    \begin{tabular}{|c|c|c|c|c|c|c|c|c|}
      \hline
     ~~~{{Methods}}~~~&~~{{3Sources}}~~& ~~{{MSRC}}~~  & ~~{{BBCsport}}~~ & ~~{{Reusters}}~~ & ~~{Caltech101}~~ &~~{{Handwritten}}~~ &~~{{MNIST}}~~& ~~{{Animal}}~~  \\
      \hline \hline
     {KNN} &~$82.00{(\pm 6.46)}$~&~$71.74{(\pm 0.36)}$~&~$93.25{(\pm 1.37)}$~&~$70.83{(\pm 1.57)}$~&~$85.47{(\pm 1.38)}$~&~$85.40{(\pm 4.99)}$~&~$92.87{(\pm 0.31)}$~&~$27.23{(\pm 0.50)}$~\\ \hline
      {CCA} &$86.20{(\pm 5.99)}$&$73.17{(\pm 0.54)}$&$95.71{(\pm 1.45)}$&$71.67{(\pm 5.11)}$&$87.73{(\pm 0.22)}$&$87.43{(\pm 1.18)}$&$92.40{(\pm 0.55)}$&$27.38{(\pm 0.36)}$\\ \hline
      {SCCA } & $63.50{(\pm 8.99)}$&$63.65{(\pm 8.84)}$&$61.76{(\pm 9.23)}$&$41.25{(\pm 12.37)}$&$82.86{(\pm 9.26)}$&$71.37{(\pm 3.54)}$&$48.39{(\pm 3.20)}$&$17.40{(\pm 2.34)}$\\ \hline
      {SGCCA} & $63.50{(\pm 8.99)}$&$63.65{(\pm 8.84)}$&$61.76{(\pm 9.23)}$&$41.25{(\pm 12.37)}$&$82.86{(\pm 9.26)}$&$71.37{(\pm 3.54)}$&$48.39{(\pm 3.20)}$&$17.40{(\pm 2.34)}$\\ \hline
      {TCCA }  &$83.00{(\pm 6.23)}$&$\underline{85.24}{(\pm \underline{4.30})}$&$91.00{(\pm 1.49)}$&$52.08{(\pm 1.76)}$& $89.98{(\pm 1.51)}$&$\underline{94.52}{(\pm \underline{1.46})}$&$83.53{(\pm 0.94)}$&$\underline{27.77}{(\pm \underline{0.64})}$\\ \hline
      {TCCA-O}  &$\underline{90.50}{(\pm \underline{1.91})}$&$73.02{(\pm 6.03)}$&$\underline{96.32}{(\pm \underline{1.07})}$&$\underline{72.91}{(\pm \underline{1.37})}$&$89.37{(\pm 1.33)}$&$80.37{(\pm 2.68)}$&$\underline{93.13}{(\pm \underline{0.13})}$&$21.89{(\pm 0.86)}$\\ \hline
      {TCCA-OS} &$83.50{(\pm 5.97)}$&$74.13{(\pm 4.97)}$&$95.19{(\pm 1.66)}$&$72.67{(\pm 1.37)}$& $\underline{90.61}{(\pm \underline{1.24})}$&$78.10{(\pm 1.92)}$&$92.49{(\pm 0.28)}$&$22.31{(\pm 0.63)}$\\ \hline
      {TMCCA}  & $64.60{(\pm 4.34)}$&$53.97{(\pm 7.18)}$&$94.58{(\pm 1.75)}$&$56.67{(\pm 1.17)}$& $84.73{(\pm 3.10)}$&$87.45{(\pm 4.73)}$&$59.50{(\pm 6.83)}$&$18.22{(\pm 4.85)}$\\ \hline
      RTSL  & $68.00{(\pm 5.06)}$&$63.49{(\pm 2.47)}$&$90.49{(\pm 1.18)}$&$71.94{(\pm 1.37)}$& $85.75{(\pm 0.96)}$&$94.50{(\pm 1.65)}$&-&-\\ \hline
      CDPML  & $79.20{(\pm 5.67 )}$&$66.03{(\pm 7.63)}$&$92.14{(\pm 2.27)}$&$58.58{(\pm 4.44)}$& $90.27{(\pm 2.05 )}$&$93.07{(\pm 1.16)}$&$87.77{(\pm 0.66)}$&$19.70{(\pm 1.21)}$\\ \hline
     STCCA-L (Our) &$\mathbf{95.00}{(\pm \mathbf{4.24})}$&$\mathbf{93.65}{(\pm \mathbf{3.42})}$ &$\mathbf{98.01}{(\pm \mathbf{0.90})}$&$\mathbf{76.11}{(\pm \mathbf{0.23})}$&$\mathbf{94.29}{(\pm \mathbf{0.39})}$&$\mathbf{98.45}{(\pm \mathbf{0.63})}$&$\mathbf{94.48}{(\pm \mathbf{0.40})}$&$\mathbf{33.06}{(\pm \mathbf{0.77})}$ \\
      \hline
    \end{tabular}
    \begin{tablenotes}    
        \footnotesize               
        \item[1] If there is insufficient memory, it will not be shown.         
      \end{tablenotes}  
       \end{threeparttable}
\end{table*}
\begin{table*}[t]
    \renewcommand\arraystretch{1.2} 
    \caption{F1 Score ($\%$) of all compared methods under the best dimensions.}\label{f1score} 
    \centering
    \setlength{\tabcolsep}{1.0mm}
    \vspace{-0.1cm}
        \begin{threeparttable}
    \begin{tabular}{|c|c|c|c|c|c|c|c|c|}
      \hline
     ~~~{{Methods}}~~~ &~{{3Sources}}~& ~{{MSRC}}~ & ~{{BBCsport}}~ & ~{{Reusters}}~ & ~{Caltech101}~&~{{Handwritten}}~ &~{{MNIST}}~& ~{{Animal}}~   \\
      \hline \hline
      {KNN} &~$75.61{(\pm 4.91)}$~&~$84.82{(\pm 3.54)}$~&~$93.54{(\pm 1.45)}$~&~$68.52{(\pm 3.02)}$~&~$58.82{(\pm 3.41)}$~&~$77.82{(\pm 1.06)}$~&~$92.36{(\pm 0.71)}$~&~$23.53{(\pm 0.32)}$~\\ \hline
      {CCA } &$83.75{(\pm 3.60)}$&${74.01}{(\pm 1.57)}$&$95.82{(\pm 1.98)}$&$71.55{(\pm 2.84)}$& $56.39{(\pm 4.48)}$&$78.07{(\pm 0.44)}$&$\underline{92.37}{(\pm \underline{0.39})}$&$23.26{(\pm 0.45)}$\\ \hline
      {SCCA } &$80.13{(\pm 6.38)}$&$\underline{88.17}{(\pm \underline{1.51})}$&$96.62{(\pm 0.79)}$&$68.20{(\pm 1.01)}$&$57.44{(\pm 3.57)}$&$79.09{(\pm 0.21)}$&$92.19{(\pm 0.22)}$&$\underline{23.95}{(\pm \underline{0.58})}$ \\ \hline
      {SGCCA} &$51.45{(\pm 12.01)}$&$57.30{(\pm 1.93)}$&$54.66{(\pm 10.31)}$&$46.68{(\pm 5.43)}$&$45.39{(\pm 6.46)}$&$74.28{(\pm 4.26)}$&$25.05{(\pm 1.77)}$&$15.24{(\pm 2.01)}$\\ \hline
      {TCCA }  &$83.41{(\pm 6.75)}$&$85.72{(\pm 2.18)}$&$90.23{(\pm 6.07)}$&$61.83{(\pm 1.00)}$& $58.69{(\pm 2.13)}$&$\underline{96.49}{(\pm \underline{0.82})}$&$77.12{(\pm 2.33)}$&$23.65{(\pm 0.61)}$\\ \hline
      {TCCA-O} &$\underline{87.51}{(\pm \underline{4.66})}$&$72.44{(\pm 2.25)}$&$\underline{97.20}{(\pm \underline{2.32})}$&${71.89}{(\pm {3.01})}$ & $\underline{71.34}{(\pm \underline{1.43})}$&$92.05{(\pm 0.82)}$&$91.13{(\pm 0.24)}$&$18.10{(\pm 1.02)}$\\ \hline
     {TCCA-OS} &$80.24{(\pm 3.92)}$&$68.72{(\pm 4.89)}$&$95.64{(\pm 1.44)}$&$67.07{(\pm 0.22)}$& $71.13{(\pm 1.92)}$&$91.39{(\pm 1.26)}$&$91.78{(\pm 0.33)}$&$18.22{(\pm 0.98)}$\\ \hline
     {TMCCA}  &$87.38{(\pm 4.95)}$&$79.69{(\pm 3.27)}$&$94.08{(\pm 0.51)}$&$62.77{(\pm 2.28)}$& $52.23{(\pm 6.93)}$&$86.52{(\pm 5.78)}$&$59.75{(\pm 4.95)}$&$16.22{(\pm 3.57)}$\\ \hline
     RTSL  & $37.75{(\pm 4.90)}$&$61.42{(\pm 3.06)}$&$88.26{(\pm 1.69)}$&$\underline{72.06}{(\pm \underline{1.46})}$& $48.83{(\pm 3.69)}$&$94.55{(\pm 1.65)}$&-&-\\ \hline
      CDPML & $75.06{(\pm 7.84)}$&$65.15{(\pm 7.21)}$&$92.65{(\pm 2.10)}$&$58.76{(\pm 4.27)}$& $70.97{(\pm 5.59) }$&$93.09{(\pm 1.20)}$&$87.55{(\pm 0.71)}$&$15.81{(\pm 1.04)}$\\ \hline
      STCCA-L (Our) &$\mathbf{95.16}{(\pm \mathbf{4.55})}$&$\mathbf{93.54}{(\pm \mathbf{2.08})}$&$\mathbf{98.63}{(\pm \mathbf{0.47})}$&$\mathbf{75.05}{(\pm \mathbf{2.67})}$&$\mathbf{77.25}{(\pm \mathbf{1.56})}$&$\mathbf{98.05}{(\pm \mathbf{0.70})}$&$\mathbf{94.44}{(\pm \mathbf{0.36})}$&$\mathbf{27.36}{(\pm \mathbf{0.44})}$\\
      \hline
    \end{tabular}
     \begin{tablenotes}    
        \footnotesize               
        \item[1] If there is insufficient memory, it will not be shown.         
      \end{tablenotes}  
       \end{threeparttable}
\end{table*}
 \begin{figure*}[!h]
    \centering
    \subfigcapskip=-1pt
     \subfigure[KNN]{
        \centering
        \includegraphics[width=4.3cm]{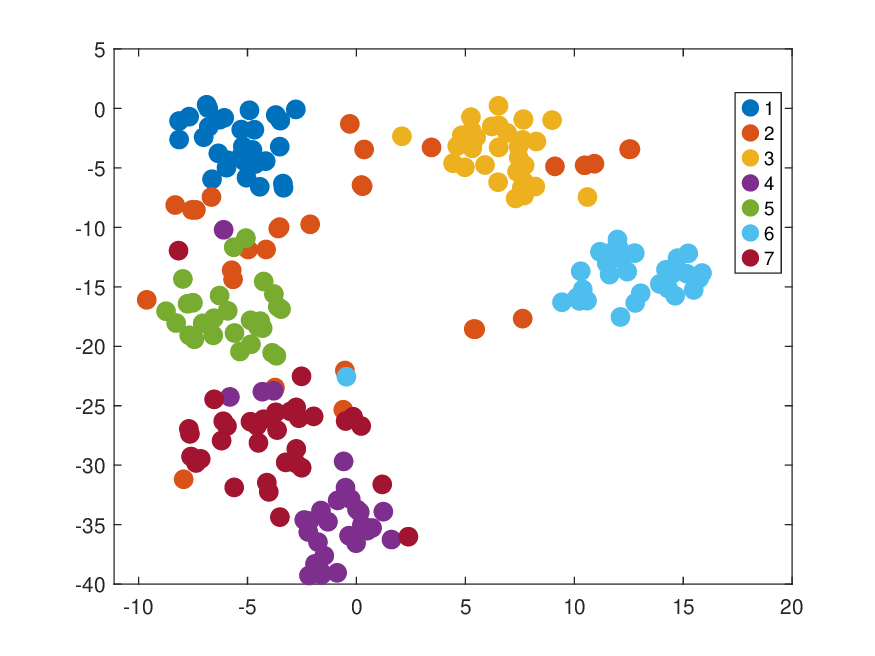}
    }\hspace{-2mm}
    \subfigcapskip=-1pt
    \subfigure[CCA]{
        \centering
        \includegraphics[width=4.3cm]{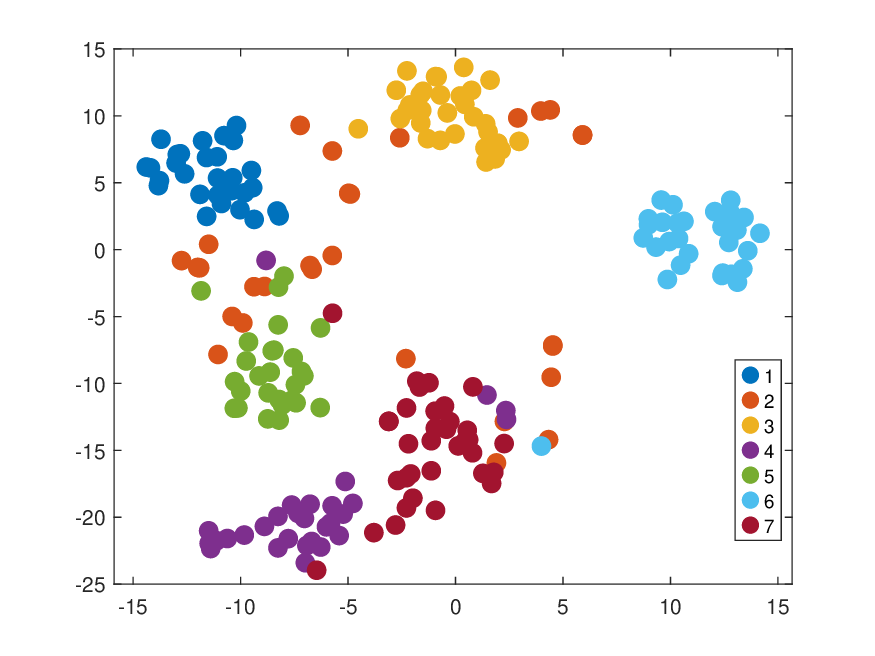}
    }\hspace{-2mm}
    \subfigcapskip=-1pt
    \subfigure[SCCA]{
        \centering
        \includegraphics[width=4.3cm]{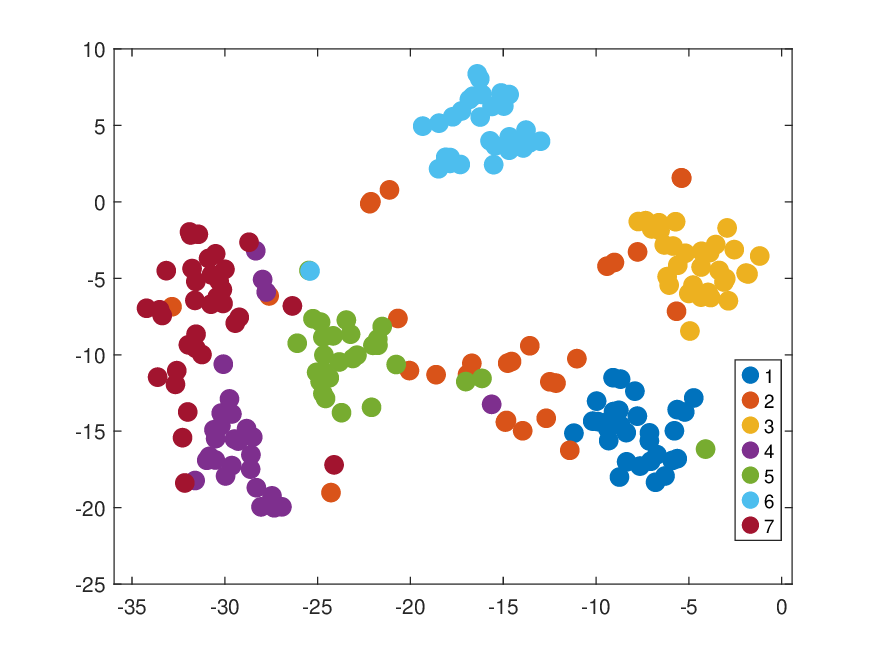}
    }\hspace{-2mm}
    \subfigcapskip=-1pt
    \subfigure[SGCCA]{
        \centering
        \includegraphics[width=4.3cm]{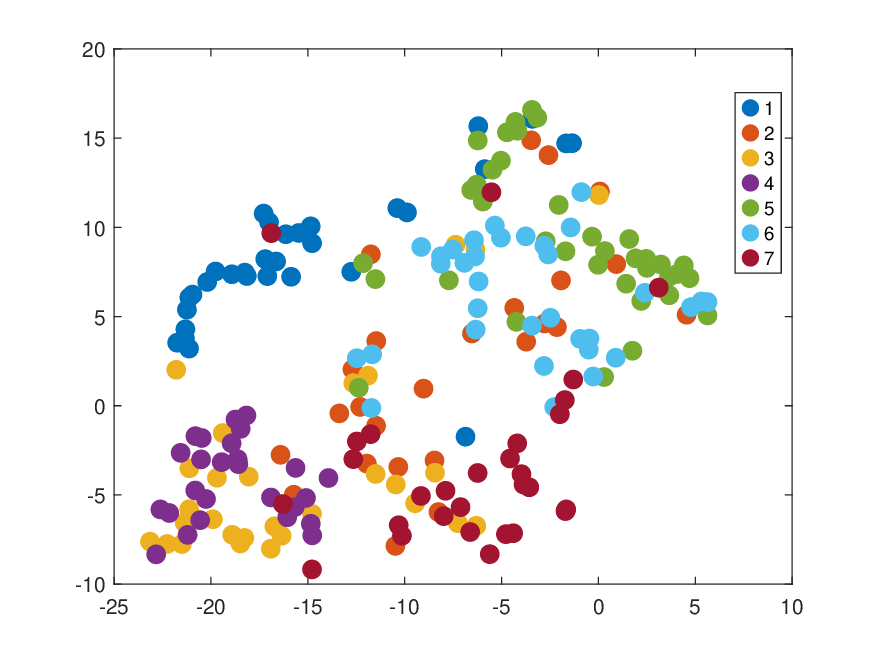}
    }\hspace{-2mm}

 \subfigure[TCCA]{
        \centering
        \includegraphics[width=4.3cm]{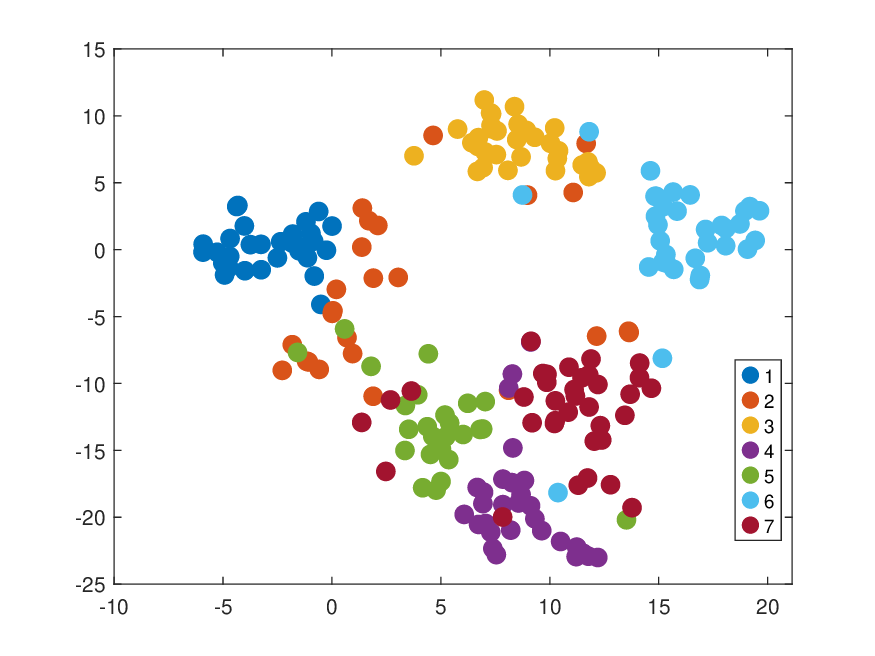}
    }\hspace{-2mm}
    \subfigcapskip=-1pt
    \subfigure[TCCA-O]{
        \centering
        \includegraphics[width=4.3cm]{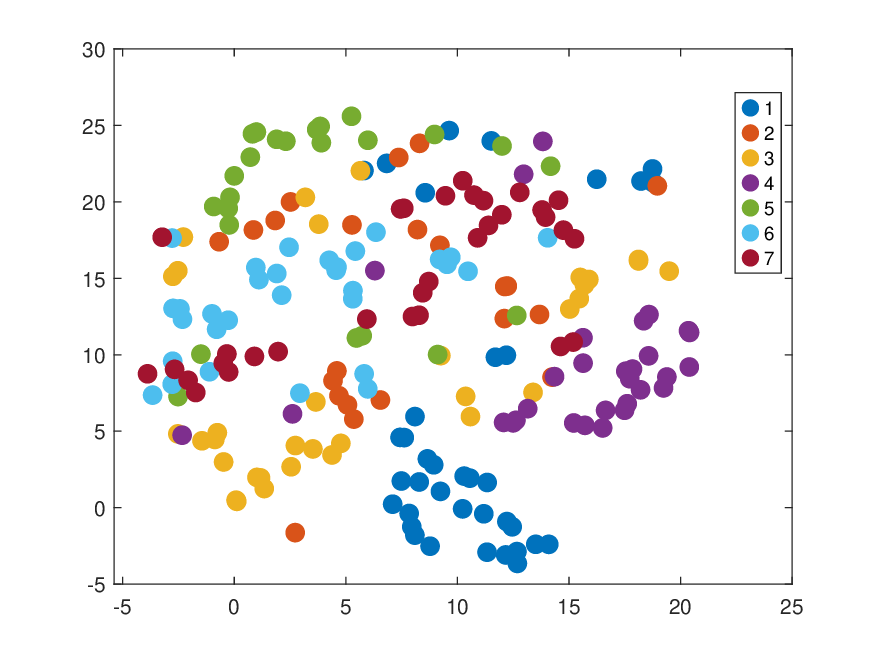}
    }\hspace{-2mm}
    \subfigcapskip=-1pt
    \subfigure[TCCA-OS]{
        \centering
        \includegraphics[width=4.3cm]{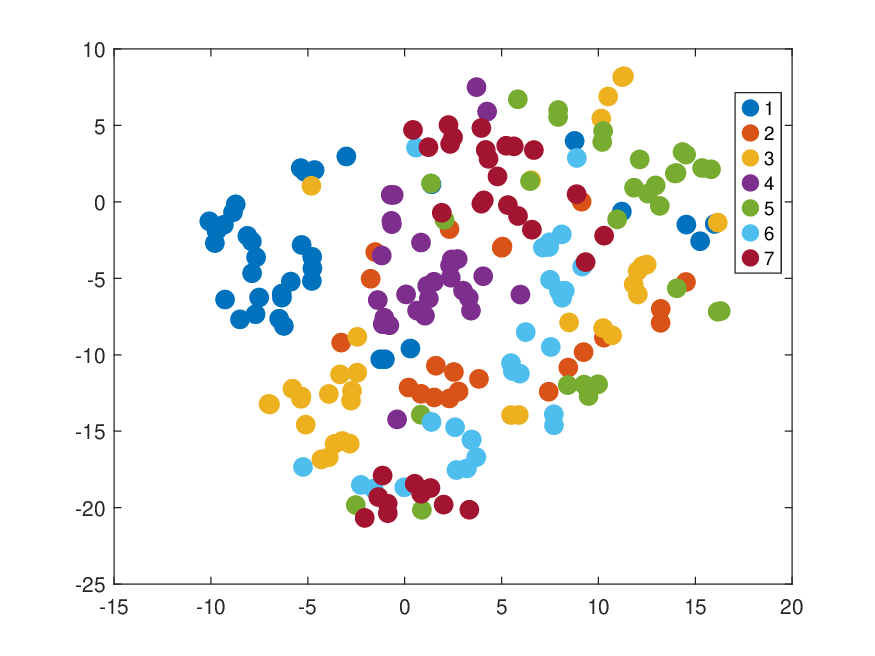}
    }\hspace{-2mm}
    \subfigcapskip=-1pt
    \subfigure[TMCCA]{
        \centering
        \includegraphics[width=4.3cm]{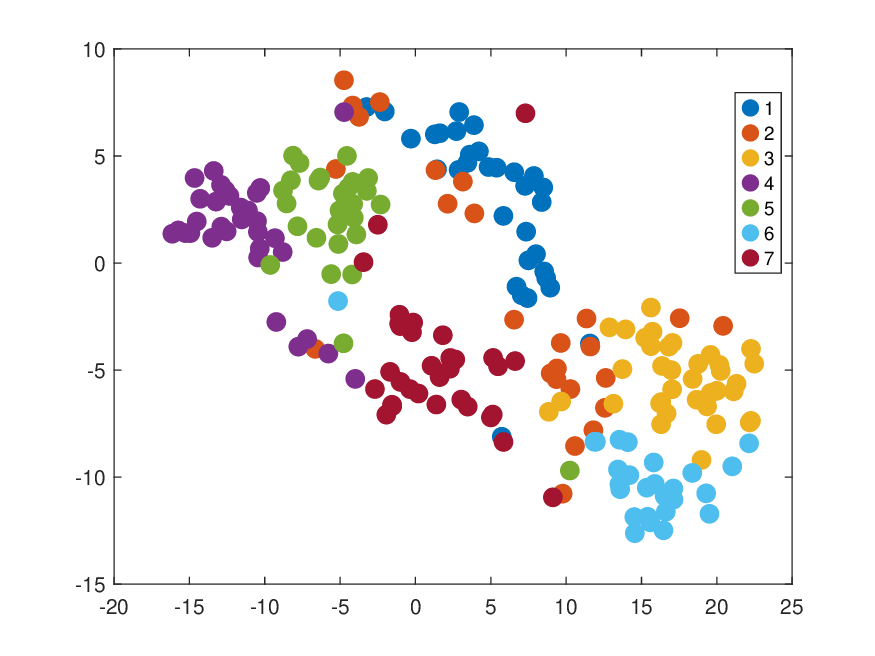}
    }\hspace{-2mm}
    
    \subfigure[RTSL]{
        \centering
        \includegraphics[width=4.3cm]{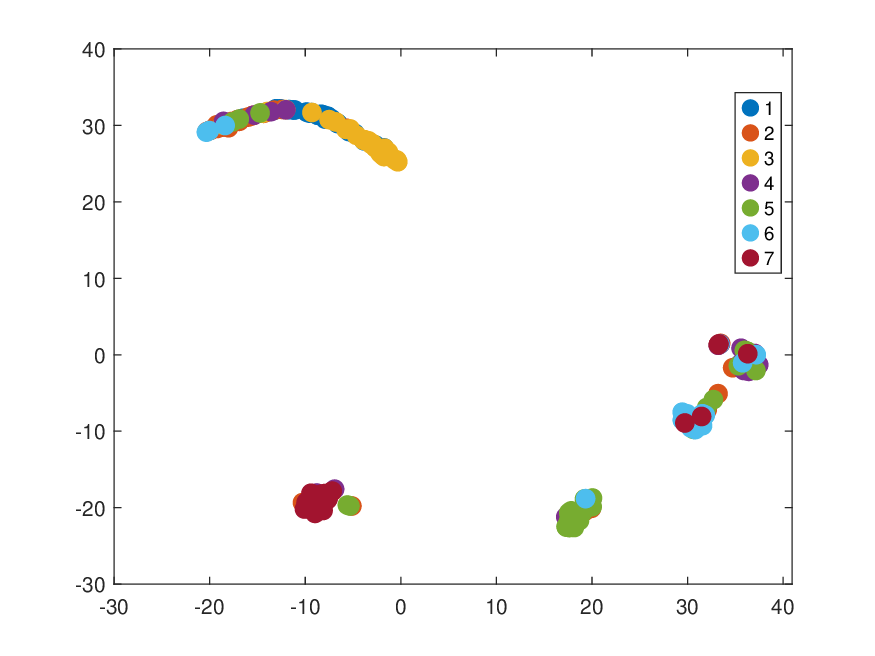}
    }\hspace{-2mm}
    \subfigcapskip=-1pt
    \subfigure[CDPMVL]{
        \centering
        \includegraphics[width=4.3cm]{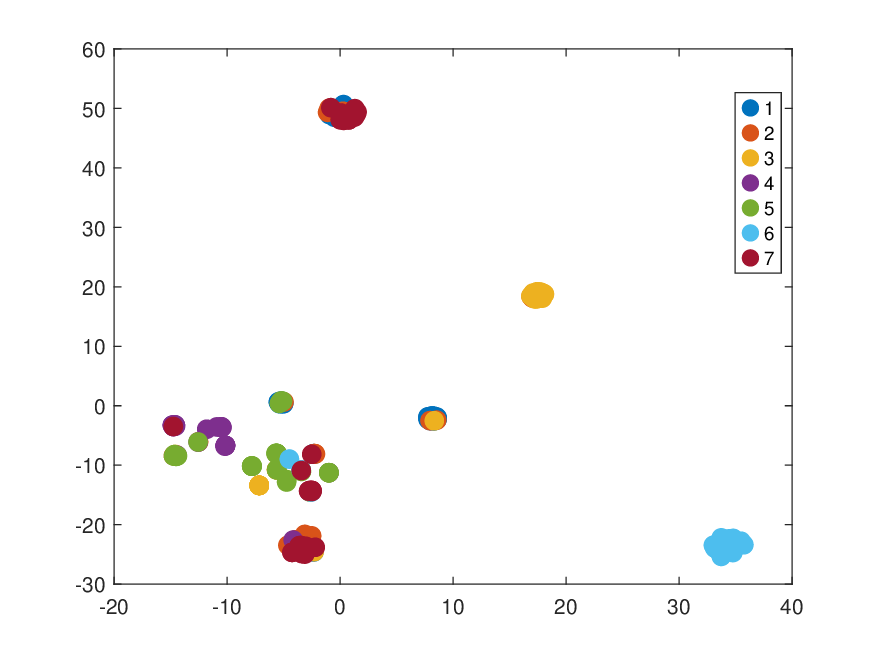}
    }\hspace{-2mm}
    \subfigcapskip=-1pt
    \subfigure[STCCA-L]{
        \centering
        \includegraphics[width=4.3cm]{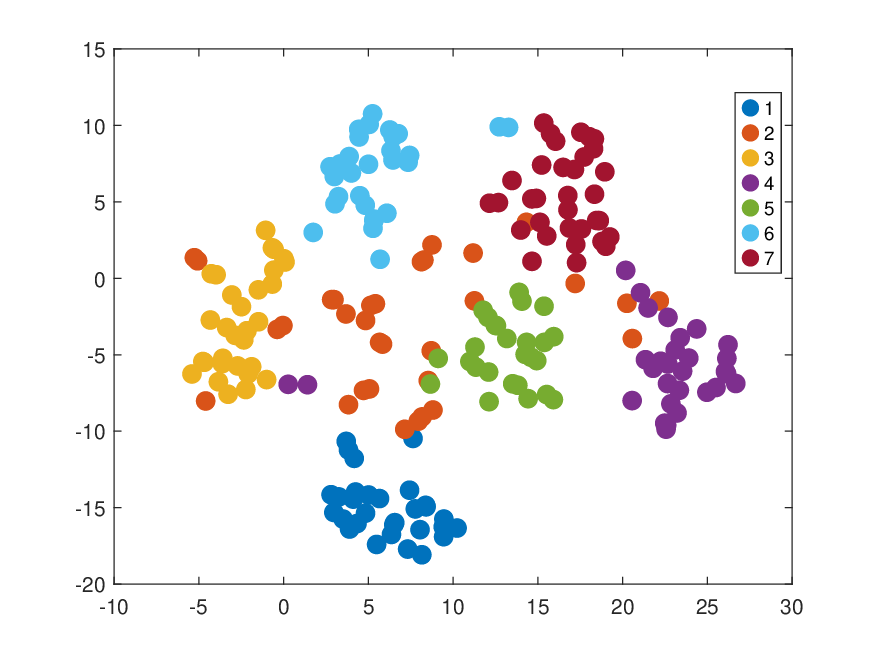}
    }\hspace{-2mm}
    \subfigcapskip=-1pt
    \subfigure[Truth]{
        \centering
        \includegraphics[width=4.3cm]{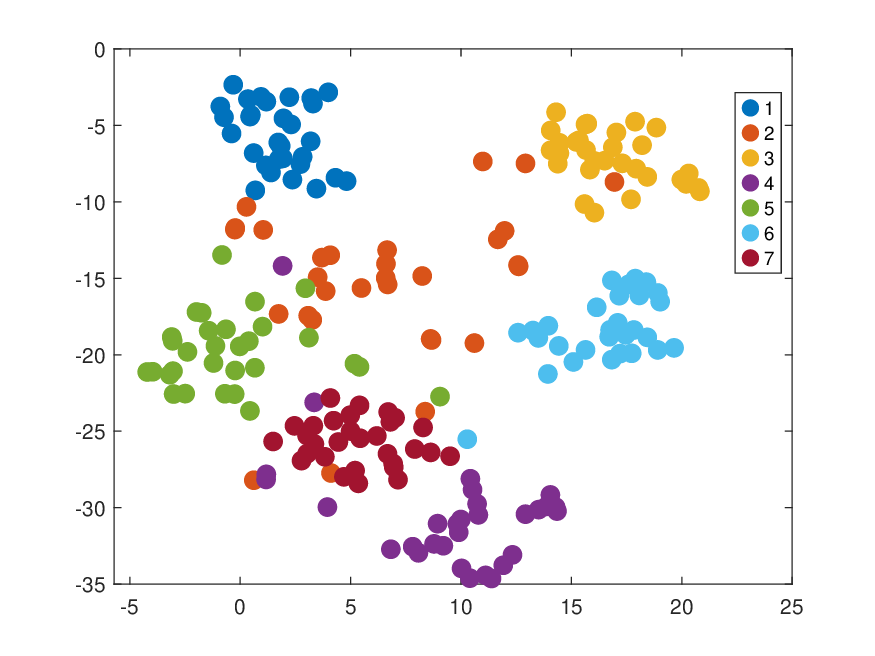}
    }\hspace{-2mm}
      \vspace{-0.1cm}
\caption{Visualization of t-SNE on the MSRC dataset, where (a)-(k) are the results of compared methods and (l) is the truth.}\label{tnse}
\end{figure*}

\begin{figure*}[!h]
    \centering
    \subfigcapskip=-1pt
     \subfigure[3Sources]{
        \centering
        \includegraphics[width=4.3cm]{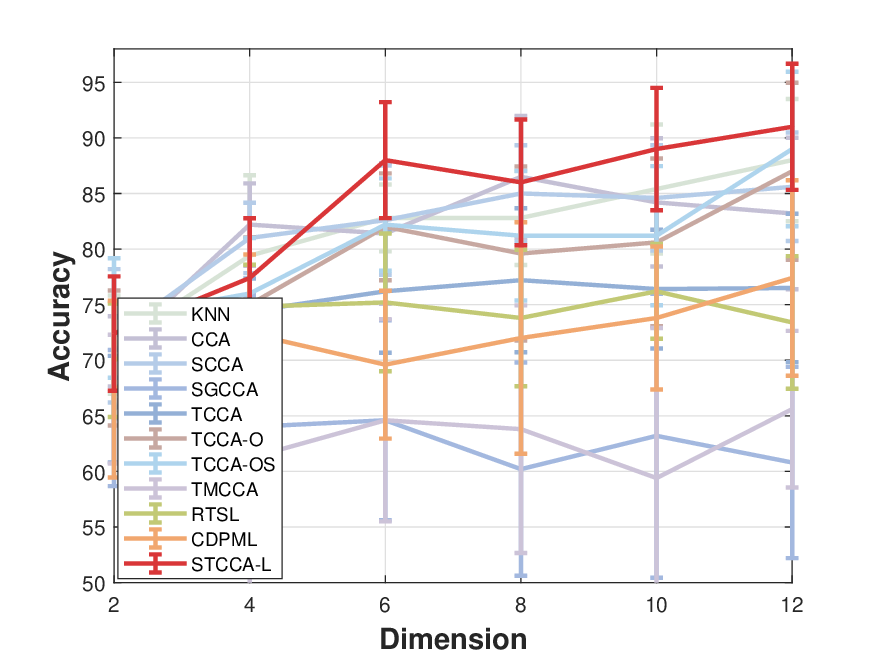}
    }\hspace{-2mm}
    \subfigcapskip=-1pt
    \subfigure[MSRC]{
        \centering
        \includegraphics[width=4.3cm]{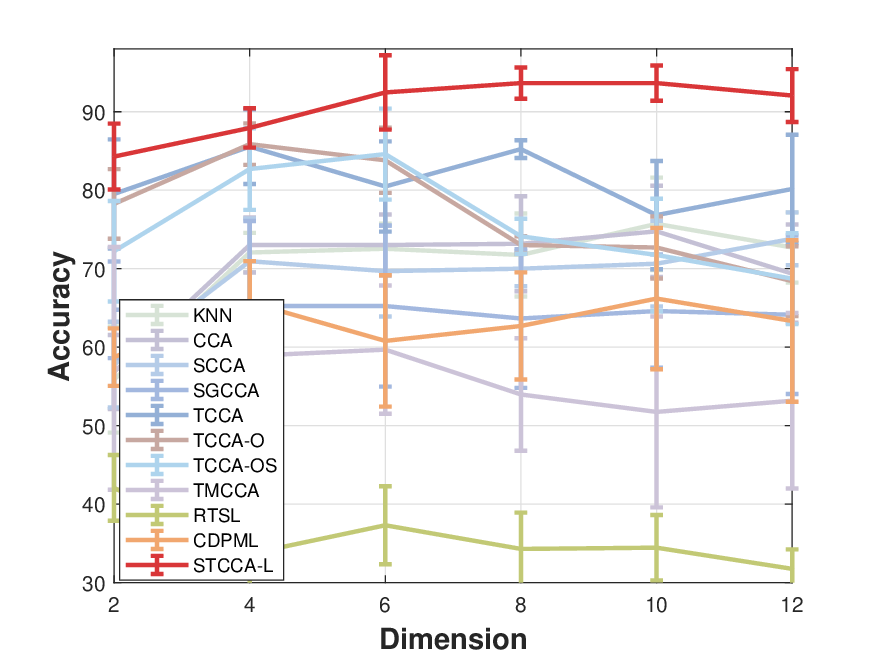}
    }\hspace{-2mm}
    \subfigcapskip=-1pt
    \subfigure[BBCsport]{
        \centering
        \includegraphics[width=4.3cm]{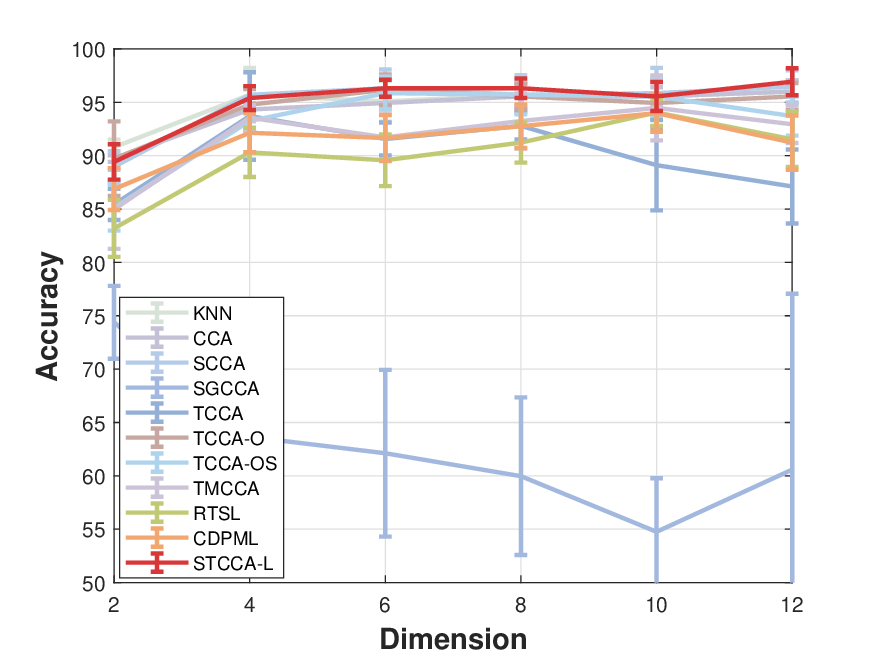}
    }\hspace{-2mm}
    \subfigcapskip=-1pt
    \subfigure[Reusters]{
        \centering
        \includegraphics[width=4.3cm]{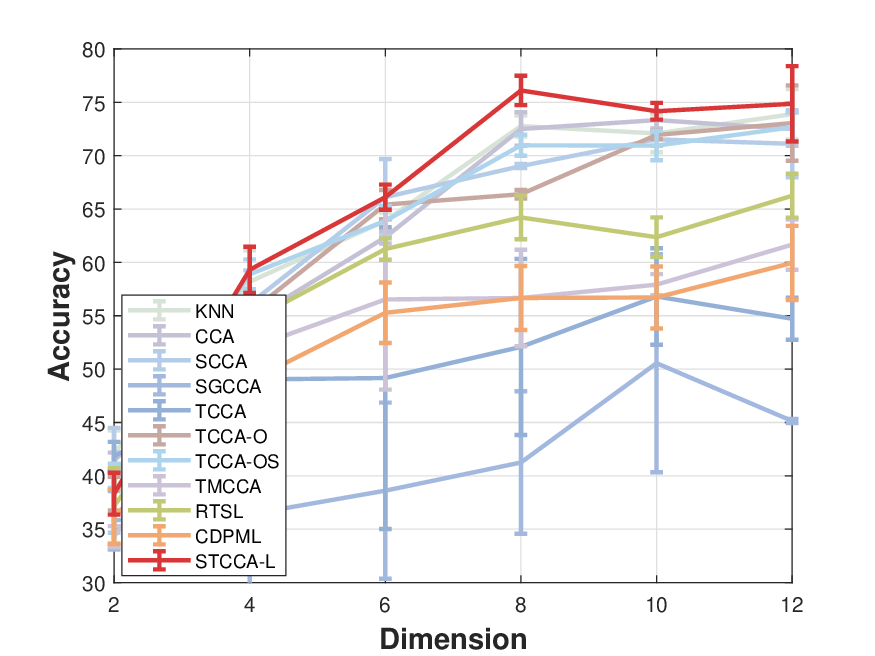}
    }\hspace{-2mm}

    \subfigure[Caltech101]{
        \centering
        \includegraphics[width=4.3cm]{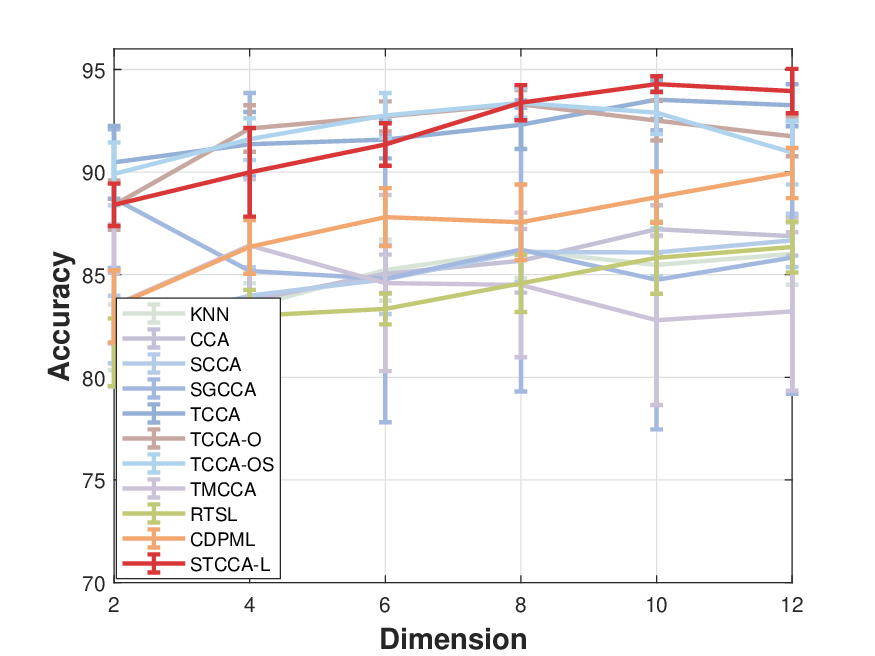}
    }\hspace{-2mm}
    \subfigcapskip=-1pt
    \subfigure[Handwritten]{
        \centering
        \includegraphics[width=4.3cm]{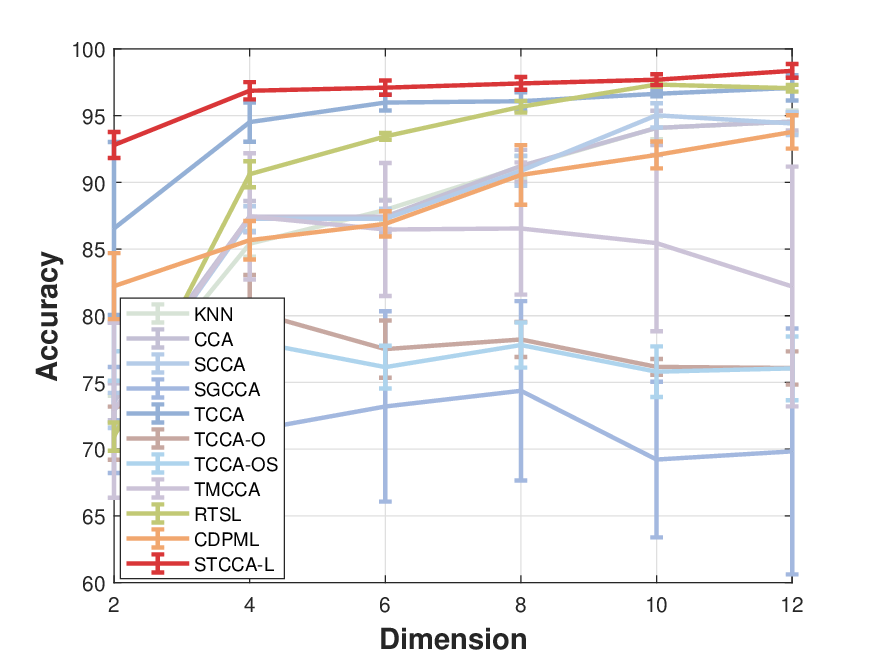}
    }\hspace{-2mm}
    \subfigcapskip=-1pt
    \subfigure[MNIST]{
        \centering
        \includegraphics[width=4.3cm]{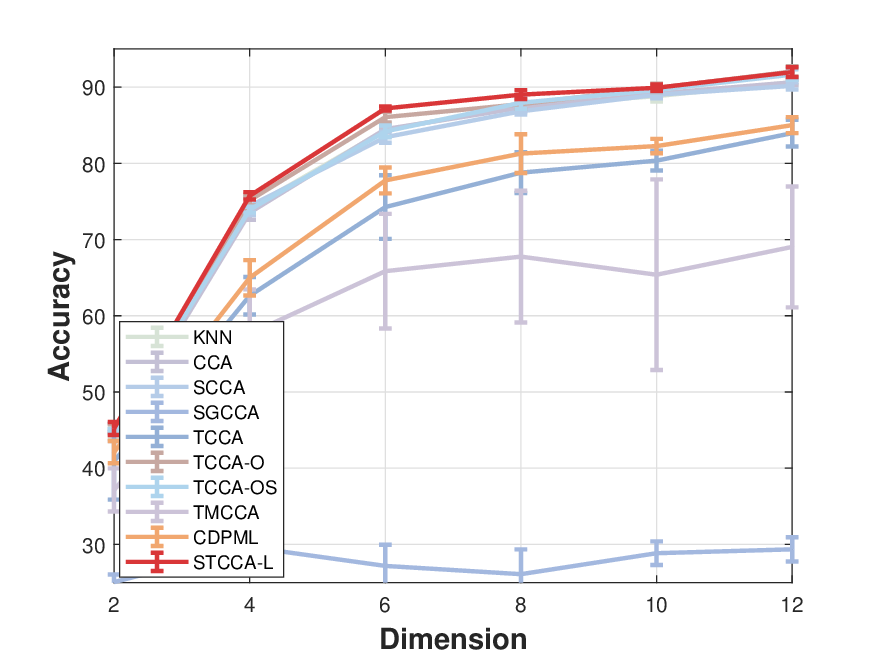}
    }\hspace{-2mm}
    \subfigcapskip=-1pt
    \subfigure[Animal]{
        \centering
        \includegraphics[width=4.3cm]{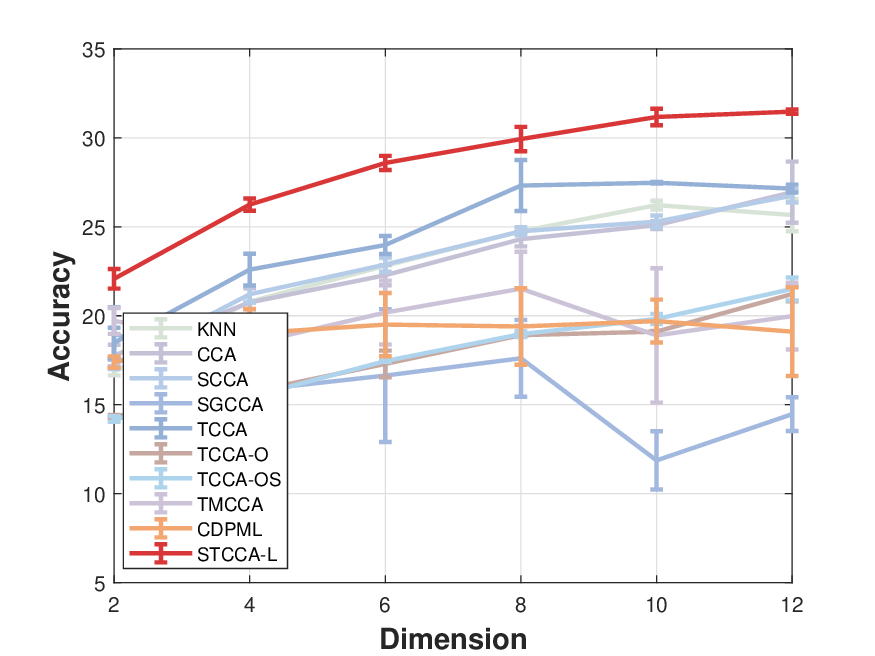}
    }\hspace{-2mm}
      \vspace{-0.1cm}
\caption{Classification accuracy of all compared methods under different dimensions.}\label{dim}
\end{figure*}

\subsubsection{Comparison Methods}
To evaluate its effectiveness, the proposed STCCA-L is compared against the classical KNN classifier and a range of state-of-the-art CCA methods. These benchmarks include matrix CCA methods such as CCA\footnote{https://github.com/tmarino2/scca} (2009), SCCA\footnote{https://github.com/htpusa/scanoncorr} (2014), and SGCCA\footnote{https://github.com/kelenlv/SGCCA2023} (2024), as well as tensor CCA methods including TCCA\footnote{https://github.com/yluopku/TCCA} (2015), TCCA-O\footnote{ https://github.com/xianchaoxiu/TCCA \label{TCCAO}} (2023), TCCA-OS\textsuperscript{\ref{TCCAO}} (2023), and TMCCA (2025). Specifically, it is compared with two state-of-the-art multi-view learning methods, robust tensor subspace learning (RTSL)\footnote{https://github.com/suxiao1824308603/Multi-view-Learning}(2024) and consensus and diversity-fusion partial-view-shared multi-view learnin (CDPML)\footnote{https://github.com/zzf495/CDPMVL} (2025).

\subsubsection{Evaluation Measures}
In this paper, it employs classification accuracy and F1 Score as standard metrics.

\textbf{Accuracy} is a key metric that measures how accurate the classification model produces.   Accuracy is defined as
\begin{equation*}
  \textrm{Accuracy}=\frac{\sum_{i=1}^{C}(\textrm{TP}_{i}+\textrm{TN}_{i})}{\sum_{i=1}^{C}(\textrm{TP}_{i}+\textrm{FP}_{i}+\textrm{TN}_{i}+\textrm{ FN}_{i})},
\end{equation*}
where $C$ is the number of types, $\textrm{TP}_{i}$ is the number of type $i$ samples that are successfully predicted, $\textrm{TN}_{i}$ is the number of other types samples that are successfully predicted, $\textrm{FP}_{i}$ is the number of samples that wrongly predict other types of samples as type $i$, $\textrm{FN}_{i}$ is the number of of type $i$ samples that are wrongly predicted as those of other types.

\textbf{F1 Score} is calculated by combining the precision and recall of the model. The F1 score can be particularly useful when the class distribution is unbalanced and the user is seeking a trade-off between precision and recall. F1 score is defined as
\begin{equation*}
  \textrm{F1 score}=\frac{1}{C}\sum_{i=1}^{C}\frac{2\textrm{TP}_{i}}{2\textrm{TP}_{i}+\textrm{FP}_{i}+\textrm{FN}_{i}}.
\end{equation*}

A higher accuracy and F1 score value indicates better classification performance.

\begin{table*}[t]
    \renewcommand\arraystretch{1.2} 
    \caption{Ablation studies of our proposed method. }\label{ablation} 
    \centering
  \setlength{\tabcolsep}{1.0mm}
    \vspace{-0.1cm}
    \begin{tabular}{|c|c|c|c|cc|cc|cc|}
      \hline
      ~~~\multirow{2}{*}{{Cases}}~~~ & 
     \multirow{2}{*}{{Orthogonality}}&\multirow{2}{*}{{Sparse}}& \multirow{2}{*}{{Laplacian}}&\multicolumn{2}{c|}{{3Sources}}&\multicolumn{2}{c|}{{MSRC}}&\multicolumn{2}{c|}{{Caltech101}}\\
        {~}& {~}&{~}&{~}&{Accuracy} &{F1 Score} &{Accuracy} &{F1 Score}&{Accuracy} &{F1 Score}  \\
      \hline \hline
     {Case I}&  $\times$&${\checkmark}$ &${\checkmark}$  & ~$40.50{(\pm 2.52)}$~ & ~$33.90{(\pm 5.54)}$~& ~$25.00{(\pm 7.14)}$~ & ~$24.18{(\pm 6.25)}$~ & ~$52.26{(\pm 1.15)}$~ & ~$23.18{(\pm 3.92)}$~   \\ \hline
     {Case II}&  ${\checkmark}$&${\times}$ & ${\checkmark}$ & $\underline{85.50}{(\pm \underline{7.00})}$ & $\underline{78.95}{(\pm \underline{8.71})}$&$  68.25{(\pm 3.43)}$ & $66.00{(\pm 3.63)}$ &$  87.40{(\pm 1.23)}$ & $55.05{(\pm 3.54)}$\\ \hline
    {Case III}&  ${\checkmark}$&${\checkmark}$ &${\times}$   & $82.00{(\pm 5.03)}$  & $77.00{(\pm 6.57)}$& $\underline{84.92}{(\pm \underline{4.76})}$ & $\underline{83.85}{(\pm \underline{4.89})}$ & $\underline{90.72}{(\pm \underline{1.09})}$ & $\underline{73.25}{(\pm \underline{3.26})}$\\ \hline
   Case IV&   ${\checkmark}$&${\checkmark}$ &${\checkmark}$&$\mathbf{ 91.50}{(\pm \mathbf{4.12})}$ & $\mathbf{89.67}{(\pm \mathbf{4.94})}$ &$\mathbf{89.92}{(\pm \mathbf{2.71})}$ & $\mathbf{83.84}{(\pm \mathbf{2.82})}$  & $\mathbf{92.08}{(\pm \mathbf{1.23})}$ &$\mathbf{78.42}{(\pm \mathbf{3.23})}$   \\
      \hline
    \end{tabular}
\end{table*}
\begin{figure*}[t]
    \centering
    \subfigcapskip=-1pt
     \subfigure[Case I (3Sources)]{
        \centering
        \includegraphics[width=4.5cm]{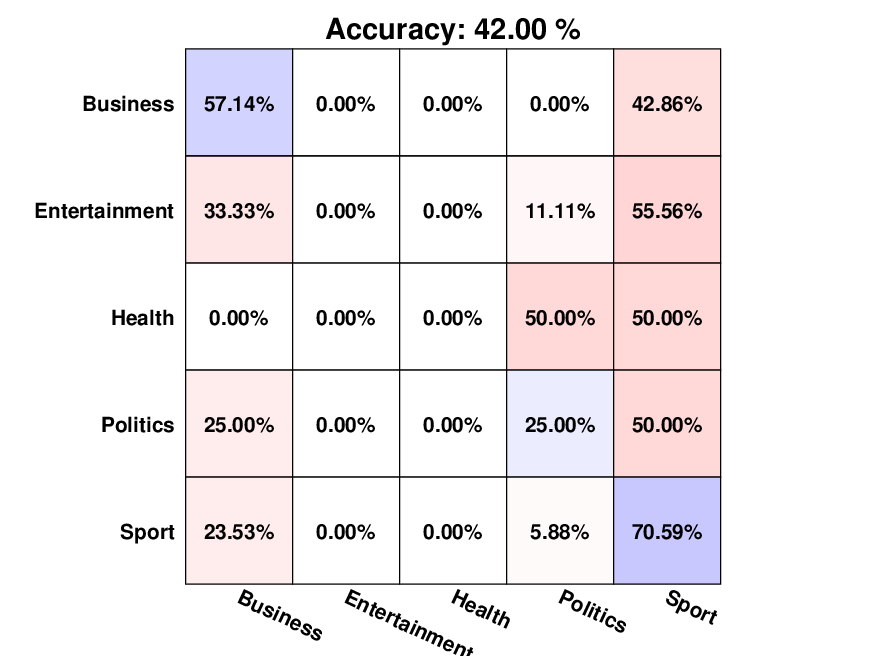}
    }\hspace{-4mm}
    \subfigcapskip=-1pt
    \subfigure[Case II (3Sources)]{
        \centering
        \includegraphics[width=4.5cm]{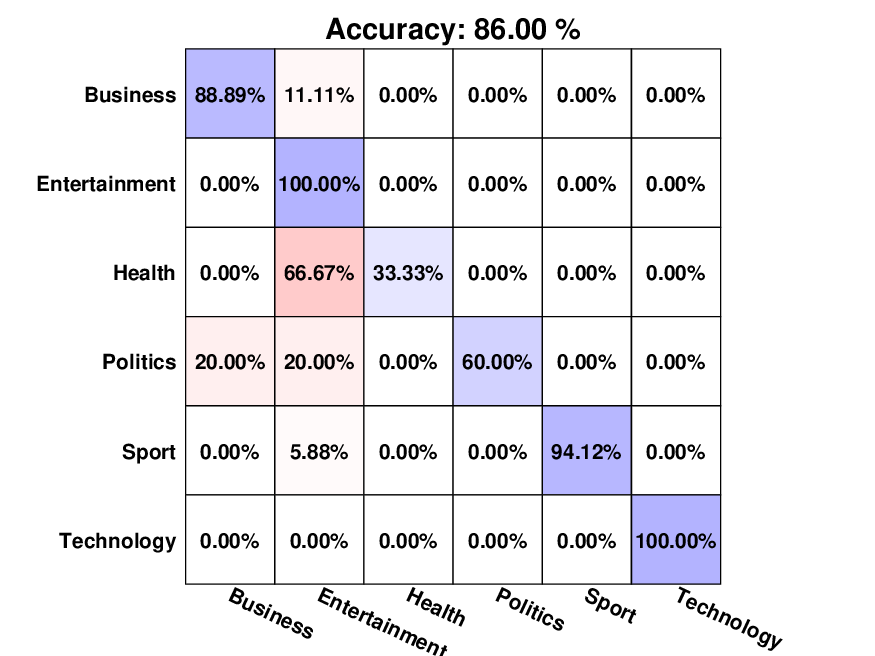}
    }\hspace{-4mm}
    \subfigcapskip=-1pt
    \subfigure[Case III (3Sources)]{
        \centering
        \includegraphics[width=4.5cm]{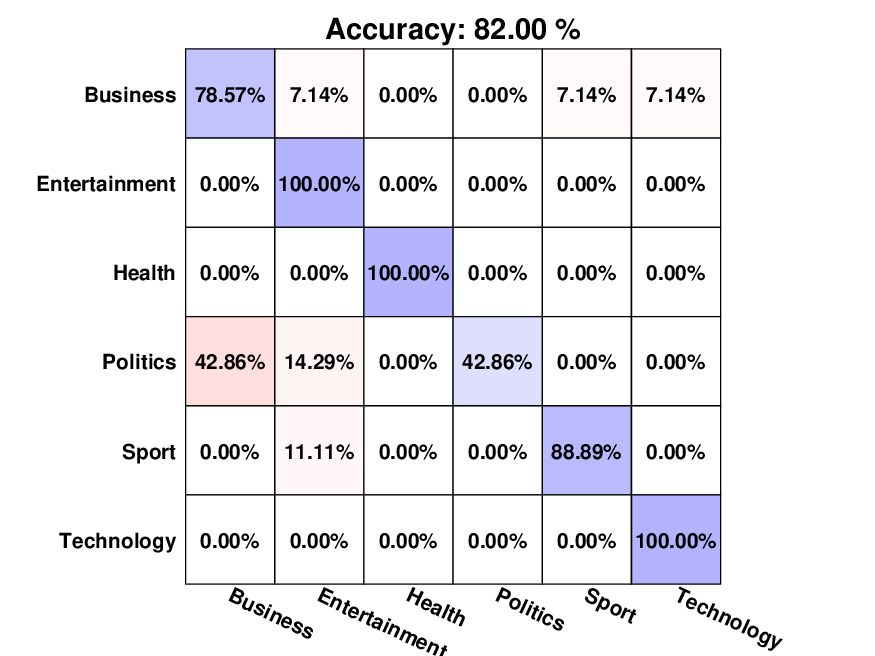}
    }\hspace{-4mm}
    \subfigcapskip=-1pt
    \subfigure[Case IV (3Sources)]{
        \centering
        \includegraphics[width=4.5cm]{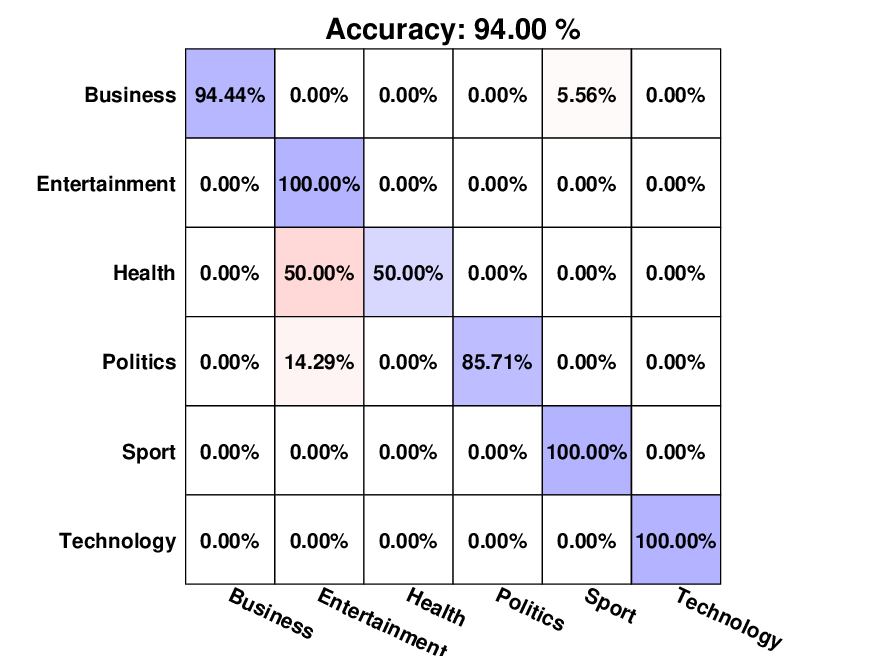}
    }\hspace{-4mm}

    \subfigure[Case I (MSRC)]{
        \centering
        \includegraphics[width=4.5cm]{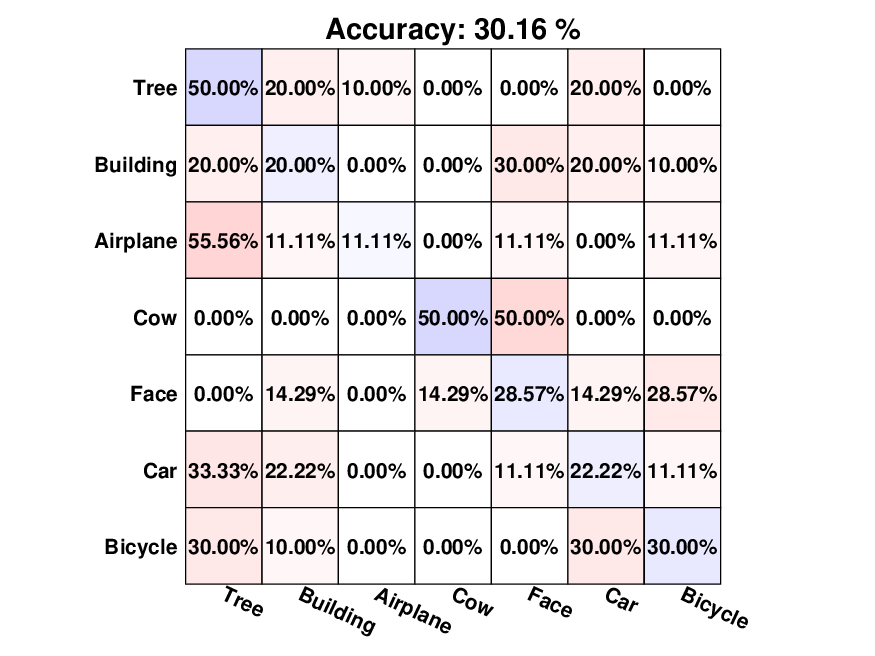}
    }\hspace{-4mm}
    \subfigcapskip=-1pt
    \subfigure[Case II (MSRC)]{
        \centering
        \includegraphics[width=4.5cm]{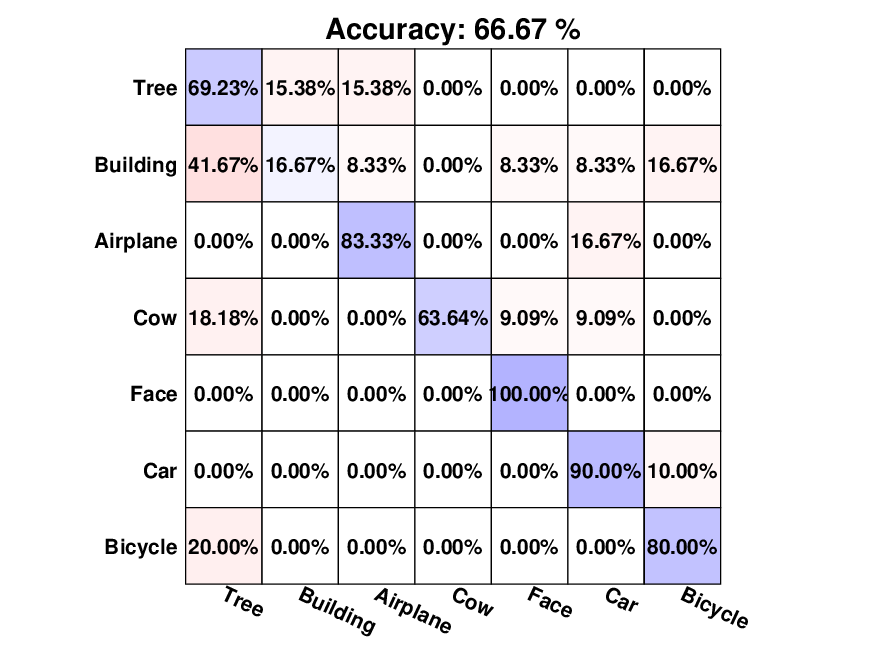}
    }\hspace{-4mm}
    \subfigcapskip=-1pt
    \subfigure[Case III (MSRC)]{
        \centering
        \includegraphics[width=4.5cm]{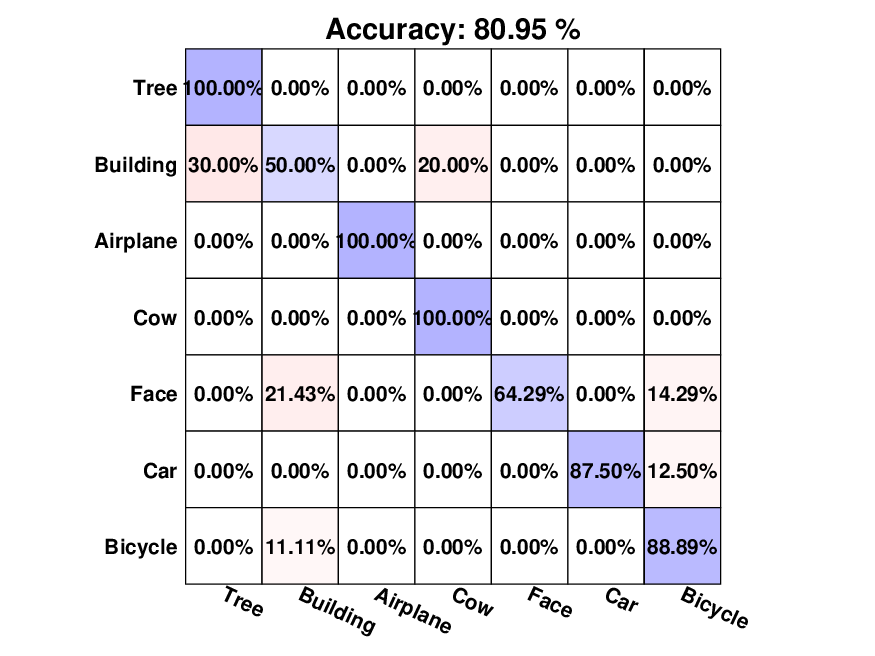}
    }\hspace{-4mm}
    \subfigcapskip=-1pt
    \subfigure[Case IV (MSRC)]{
        \centering
        \includegraphics[width=4.5cm]{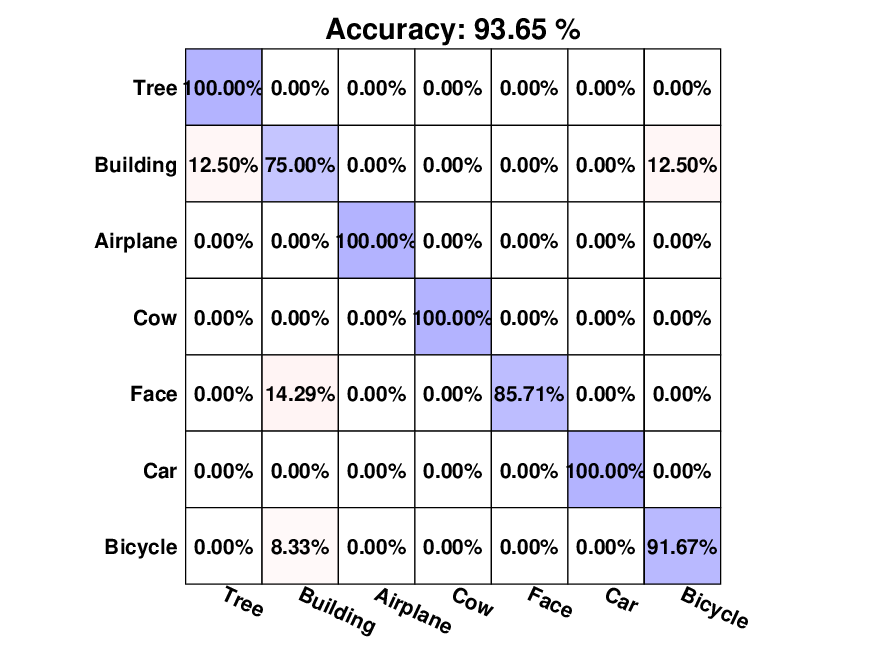}
    }\hspace{-4mm}
    
    \subfigure[Case I (Caltech101)]{
        \centering
        \includegraphics[width=4.5cm]{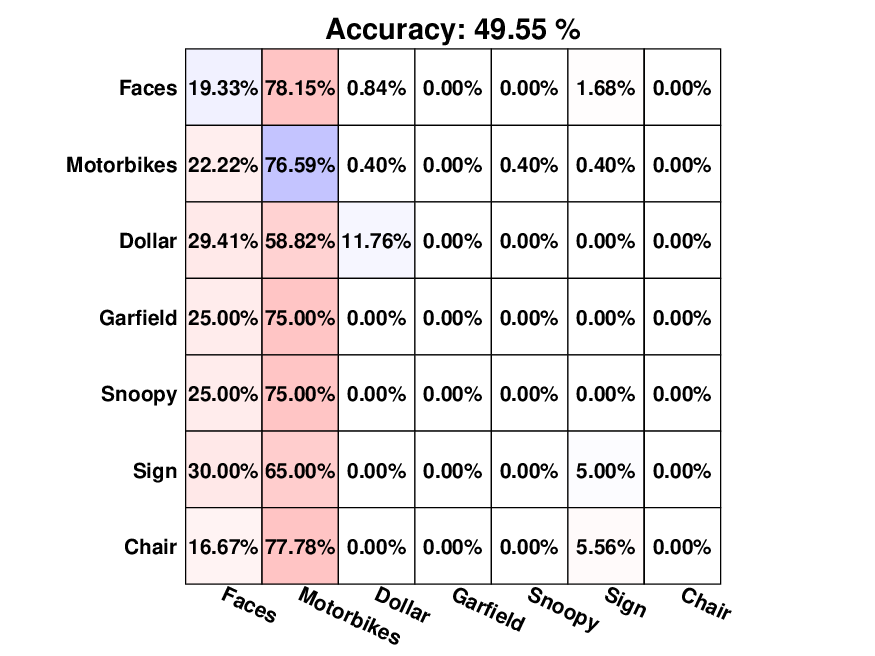}
    }\hspace{-5mm}
    \subfigcapskip=-1pt
    \subfigure[Case II (Caltech101)]{
        \centering
        \includegraphics[width=4.5cm]{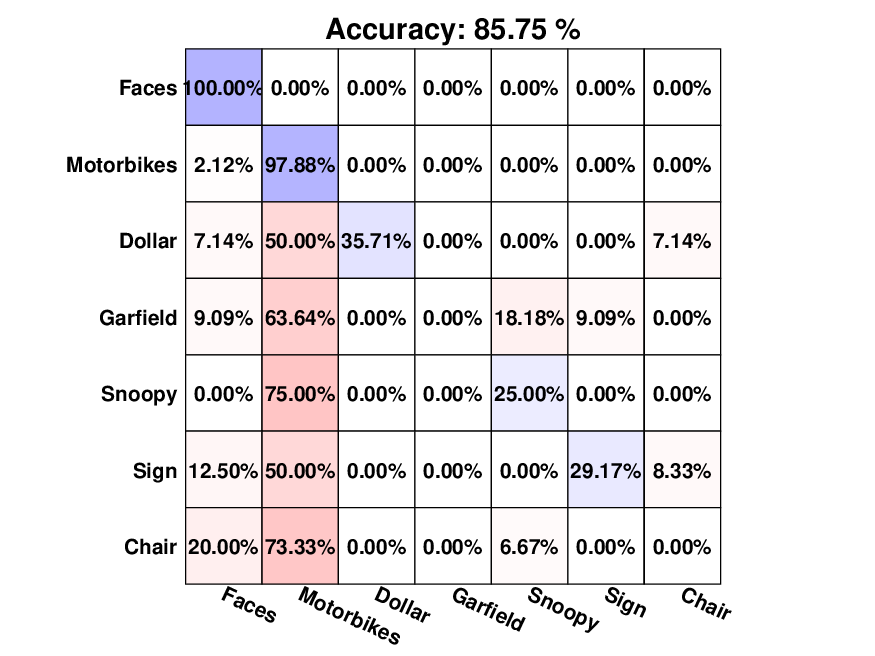}
    }\hspace{-5mm}
    \subfigcapskip=-1pt
    \subfigure[Case III (Caltech101)]{
        \centering
        \includegraphics[width=4.5cm]{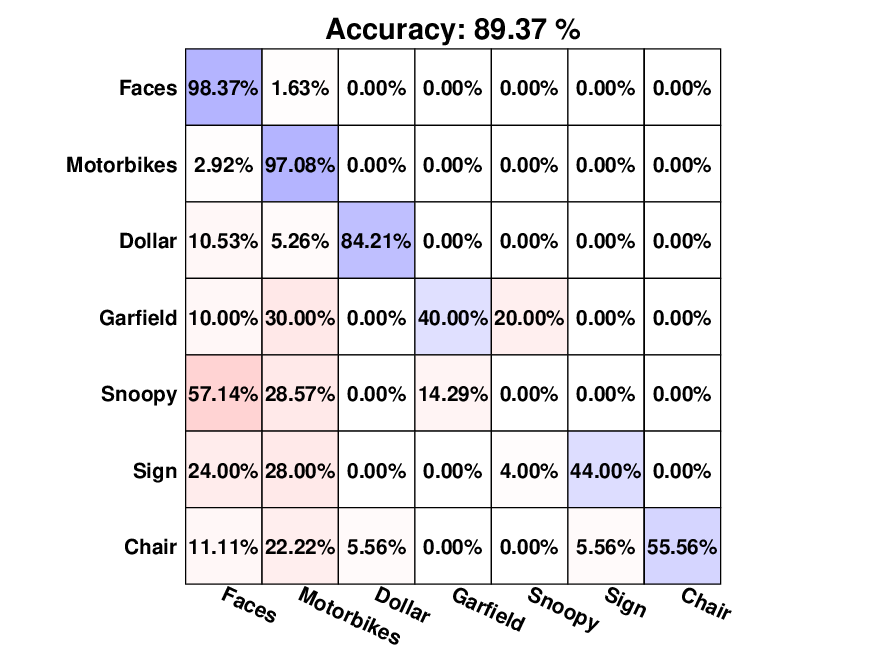}
    }\hspace{-5mm}
    \subfigcapskip=-1pt
    \subfigure[Case IV (Caltech101)]{
        \centering
        \includegraphics[width=4.5cm]{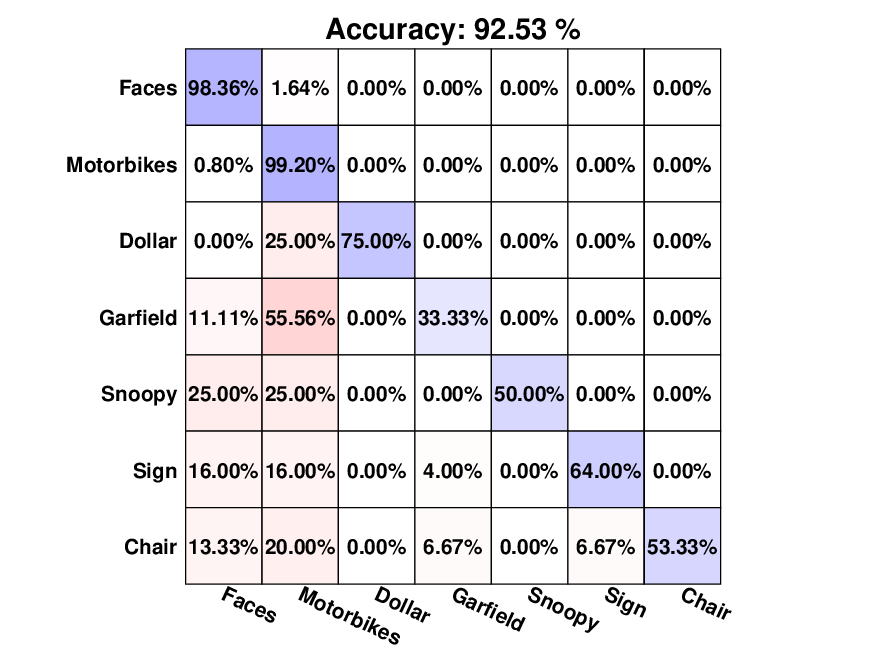}
    }\hspace{-5mm}
      \vspace{-0.1cm}
\caption{Visualization of the confusion matrix on the  3Sources, MSRC, and Caltech101 datasets.}\label{Ablation}
\end{figure*}

\begin{table*}[h] 
	\renewcommand\arraystretch{1.2} 
	\caption{Ablation studies of the graph method. }\label{abgraph} 
	\centering
	\setlength{\tabcolsep}{2.0mm}
	\vspace{-0.1cm}
			\begin{tabular}{|c|cc|cc|cc|}
				\hline
				~~~\multirow{2}{*}{{Methods}}~~~ &\multicolumn{2}{c|}{{3Sources}}&\multicolumn{2}{c|}{{MSRC}}&\multicolumn{2}{c|}{{Caltech101}}\\
				{~}&{Accuracy} &{F1 Score} &{Accuracy} &{F1 Score}&{Accuracy} &{F1 Score}  \\
				\hline \hline
				{Gaussian Graph}&  ~$78.60{(\pm 7.12)}$~ & ~$69.74{(\pm 9.01)}$~& ~$88.09{(\pm 4.38)}$~ & ~$87.25{(\pm 5.15)}$~ & ~$93.25{(\pm 1.21)}$~ & ~$77.15{(\pm 4.11)}$~  \\ \hline
				{KNN Graph}& $\underline{77.80}{(\pm 4.47)}$ & $\underline{69.75}{(\pm 7.28)}$&$90.47{(\pm 4.85)}$ & $89.51{(\pm 5.83)}$ &$ 92.89{(\pm 1.35)}$ & $74.95{(\pm 3.06)}$\\ \hline
				{Cosine Graph}& $78.40{(\pm 6.98)}$  & $68.88{(\pm 8.97)}$& ${89.36}{(\pm \underline{2.70})}$ & $88.84{(\pm 3.09)}$ & $\underline{93.34}{(\pm \underline{0.89})}$ & $\underline{76.81}{(\pm \underline{2.37})}$\\ \hline
				{Sparse Graph}& $75.40{(\pm 5.96)}$  & $64.98{(\pm 6.54)}$& $87.94{(\pm 5.03)}$ & $87.21{(\pm 5.89)}$ & $\underline{93.19}{(\pm \underline{0.87})}$ & $\underline{76.63}{(\pm \underline{2.93})}$\\ \hline
				Adaptive Neighbor Graph& $\mathbf{ 91.50}{(\pm \mathbf{4.12})}$ & $\mathbf{89.67}{(\pm \mathbf{4.94})}$ &$\mathbf{91.27}{(\pm \mathbf{2.72})}$ & $\mathbf{91.02}{(\pm \mathbf{3.16})}$  & $\mathbf{93.62}{(\pm \mathbf{0.99})}$ &$\mathbf{78.62}{(\pm \mathbf{4.41})}$   \\
				\hline
			\end{tabular}
		\end{table*}

\begin{table*}[h]
	\renewcommand\arraystretch{1.2} 
	\caption{Ablation studies of the initialization strategy. }\label{abIni} 
	\centering
	\setlength{\tabcolsep}{2.0mm}
	\vspace{-0.1cm}
			\begin{tabular}{|c|cc|cc|cc|}
				\hline
				~~~\multirow{2}{*}{{Methods}}~~~ &\multicolumn{2}{c|}{{3Sources}}&\multicolumn{2}{c|}{{MSRC}}&\multicolumn{2}{c|}{{Caltech101}}\\
				{~}&{Accuracy} &{F1 Score} &{Accuracy} &{F1 Score}&{Accuracy} &{F1 Score}  \\
				\hline \hline
				{SVD Initialization}&  ~$89.60{(\pm 3.37)}$~ & ~$85.53{(\pm 8.33)}$~& ~$87.78{(\pm 3.35)}$~ & ~$87.31{(\pm 3.29)}$~ & ~$91.04{(\pm 1.07)}$~ & ~$74.61{(\pm 1.36)}$~  \\ \hline
				{Identity Initialization}& $\underline{90.00}{(\pm \underline{4.42})}$ & $\underline{87.41}{(\pm \underline{5.47})}$&$ 87.30{(\pm 3.17)}$ & $86.99{(\pm 3.17)}$ &$ 91.47{(\pm 1.08)}$ & $74.61{(\pm 1.36)}$\\ \hline
				{Orthogonal Initialization}& $89.60{(\pm 3.09)}$  & $87.49{(\pm 4.72)}$& $\underline{83.81}{(\pm \underline{7.32})}$ & $\underline{76.04}{(\pm \underline{3.54})}$ & $\underline{91.71}{(\pm \underline{1.07})}$ & $\underline{76.04}{(\pm \underline{3.54})}$\\ \hline
				Random Initialization& $\mathbf{ 91.50}{(\pm \mathbf{4.12})}$ & $\mathbf{89.67}{(\pm \mathbf{4.94})}$ &$\mathbf{89.92}{(\pm \mathbf{2.71})}$ & $\mathbf{83.84}{(\pm \mathbf{2.82})}$  & $\mathbf{92.08}{(\pm \mathbf{1.23})}$ &$\mathbf{78.42}{(\pm \mathbf{3.23})}$   \\
				\hline
			\end{tabular}
		\end{table*}

\begin{figure*}[t]
    \centering
    \subfigcapskip=-1pt
     \subfigure[3Sources]{
        \centering
        \includegraphics[width=4.3cm]{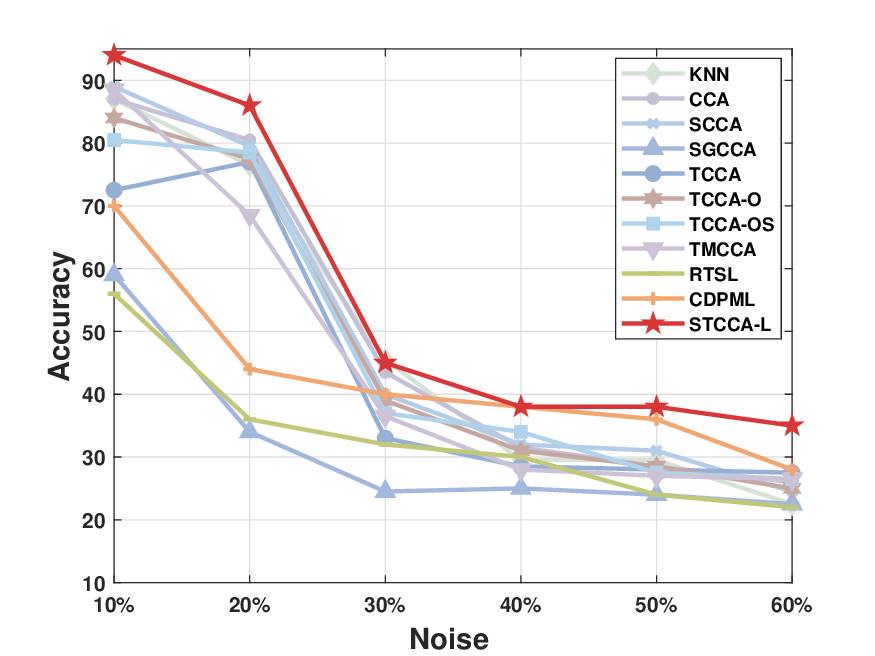}
    }\hspace{-2mm}
    \subfigcapskip=-1pt
    \subfigure[MSRC]{
        \centering
        \includegraphics[width=4.3cm]{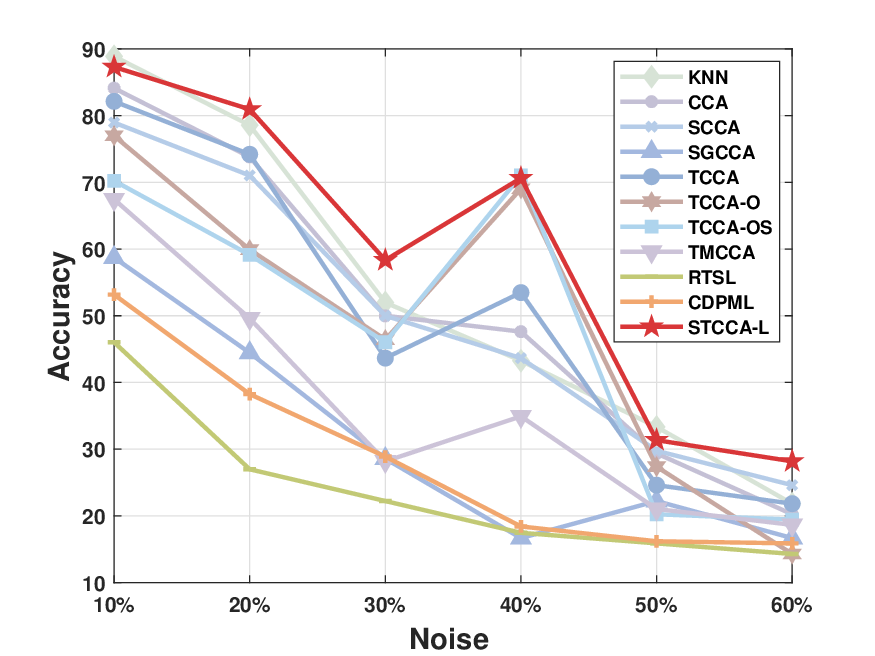}
    }\hspace{-2mm}
    \subfigcapskip=-1pt
    \subfigure[BBCsport]{
        \centering
        \includegraphics[width=4.3cm]{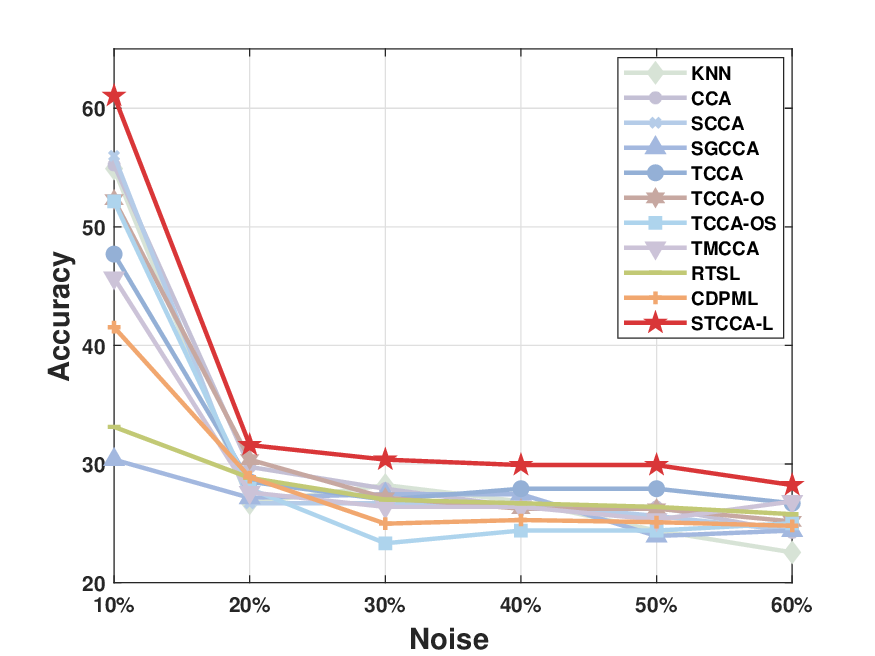}
    }\hspace{-2mm}
    \subfigcapskip=-1pt
    \subfigure[Reusters]{
        \centering
        \includegraphics[width=4.3cm]{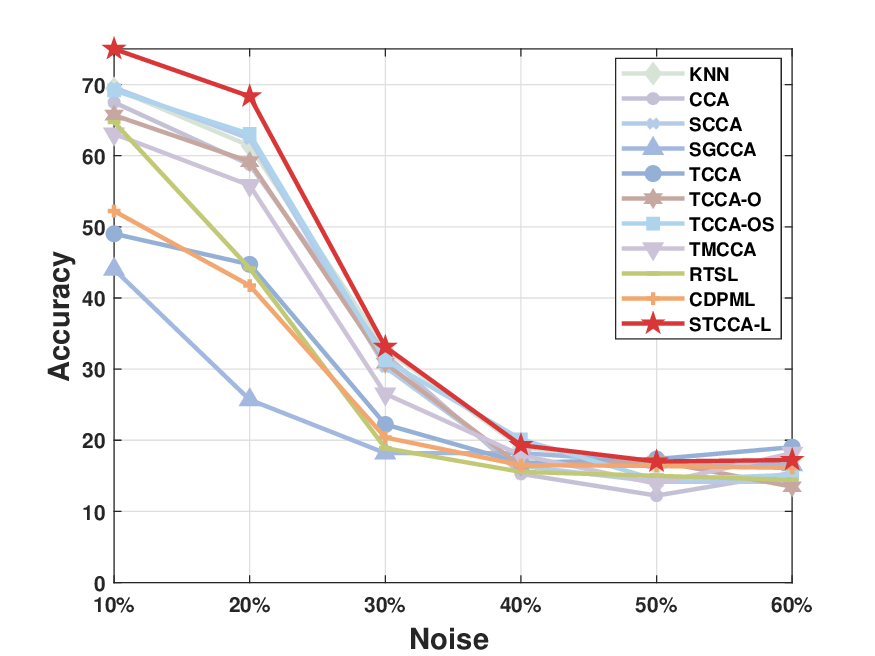}
    }\hspace{-2mm}

    \subfigure[Caltech101]{
        \centering
        \includegraphics[width=4.3cm]{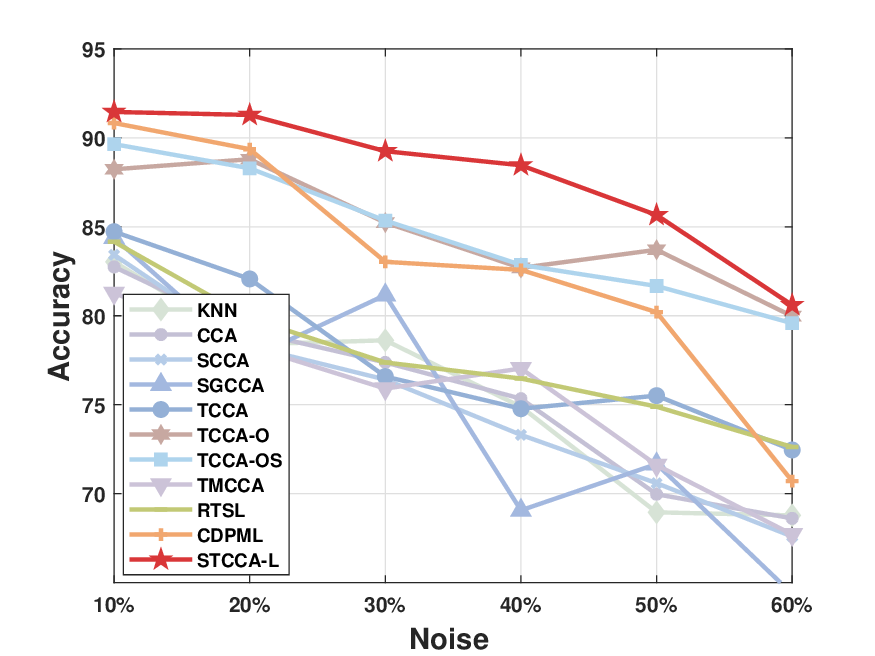}
    }\hspace{-2mm}
    \subfigcapskip=-1pt
    \subfigure[Handwritten]{
        \centering
        \includegraphics[width=4.3cm]{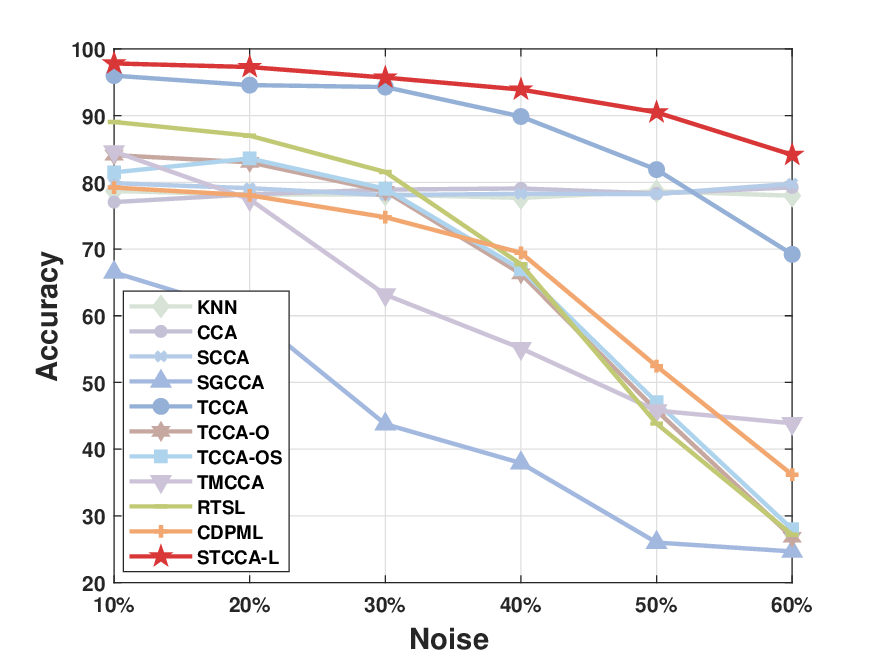}
    }\hspace{-2mm}
    \subfigcapskip=-1pt
    \subfigure[MNIST]{
        \centering
        \includegraphics[width=4.3cm]{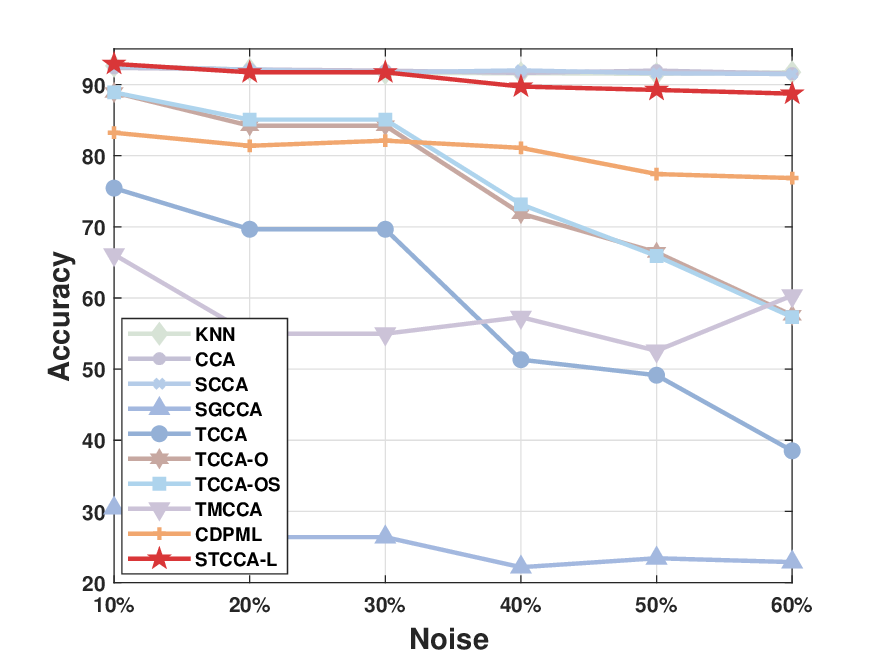}
    }\hspace{-2mm}
    \subfigcapskip=-1pt
    \subfigure[Animal]{
        \centering
        \includegraphics[width=4.3cm]{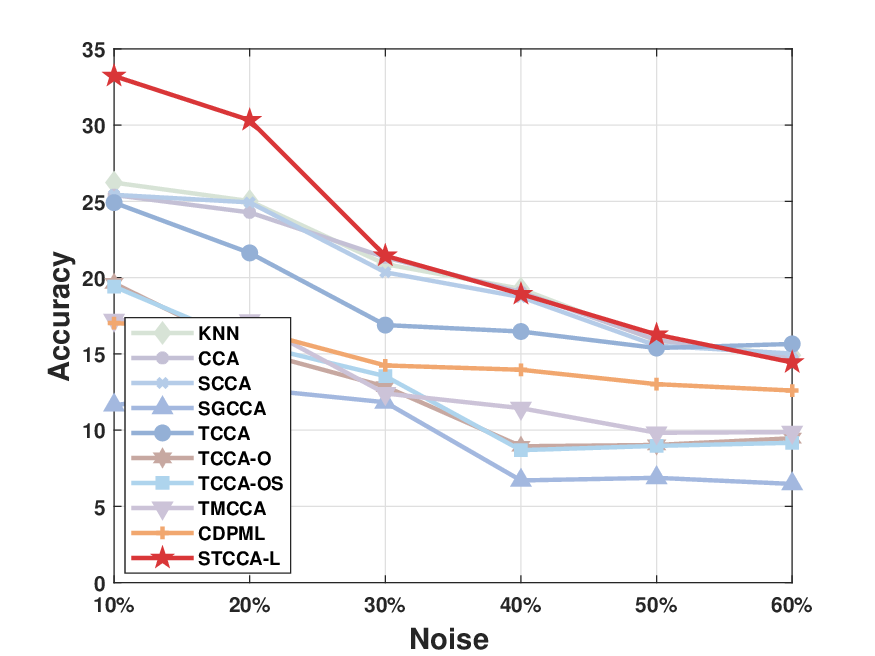}
    }\hspace{-2mm}
      \vspace{-0.1cm}
\caption{Classification accuracy of all compared methods on different datasets with varying proportions of Gaussian noise. }\label{noise}
\end{figure*}
	
\begin{figure*}[t]
    \centering
    \subfigcapskip=-1pt
     \subfigure[3Sources]{
        \centering
        \includegraphics[width=4.3cm]{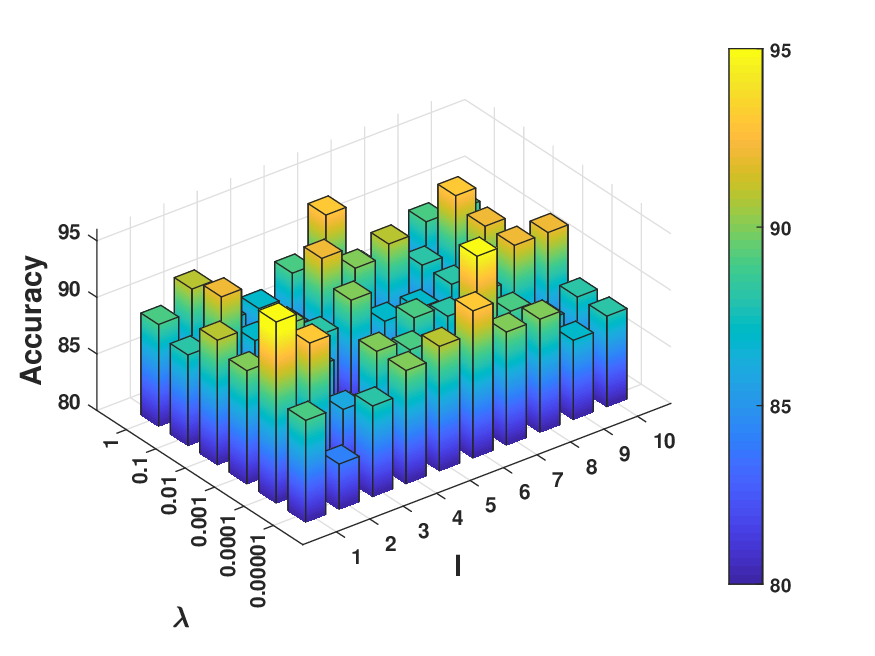}
    }\hspace{-2mm}
    \subfigcapskip=-1pt
    \subfigure[MSRC]{
        \centering
        \includegraphics[width=4.3cm]{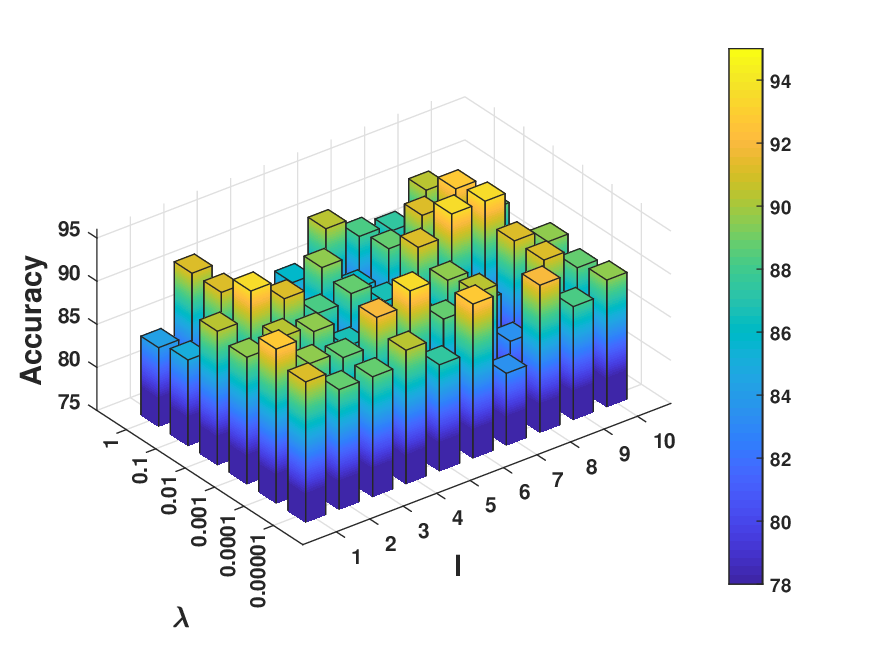}
    }\hspace{-2mm}
    \subfigcapskip=-1pt
    \subfigure[BBCsport]{
        \centering
        \includegraphics[width=4.3cm]{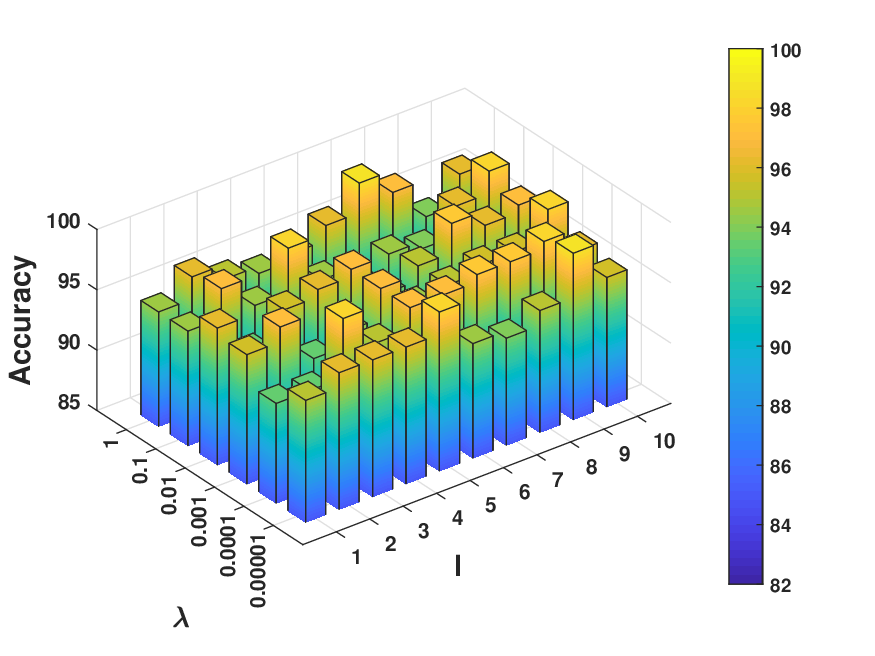}
    }\hspace{-2mm}
    \subfigcapskip=-1pt
    \subfigure[Reusters]{
        \centering
        \includegraphics[width=4.3cm]{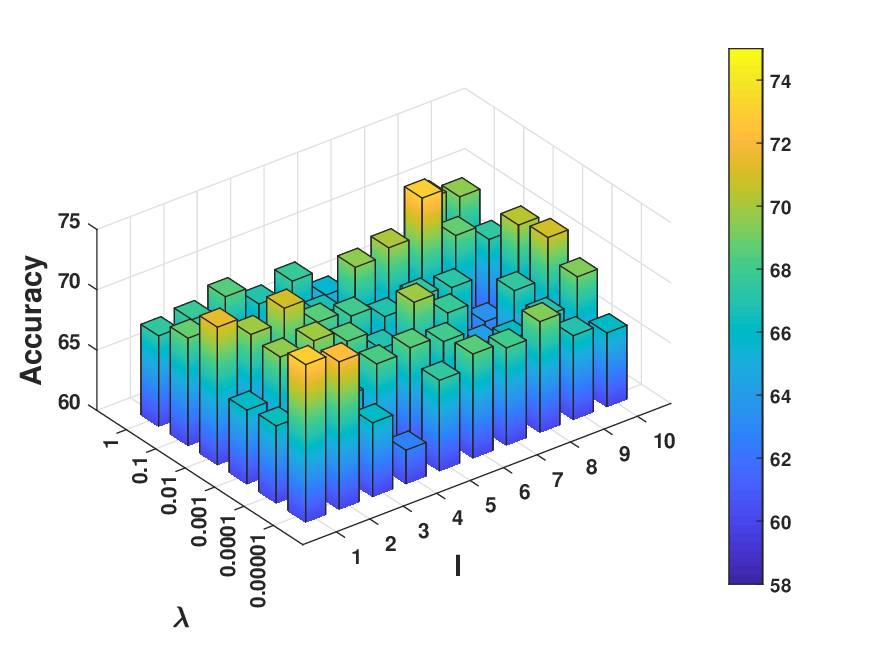}
    }\hspace{-2mm}

    \subfigure[Caltech101]{
        \centering
        \includegraphics[width=4.3cm]{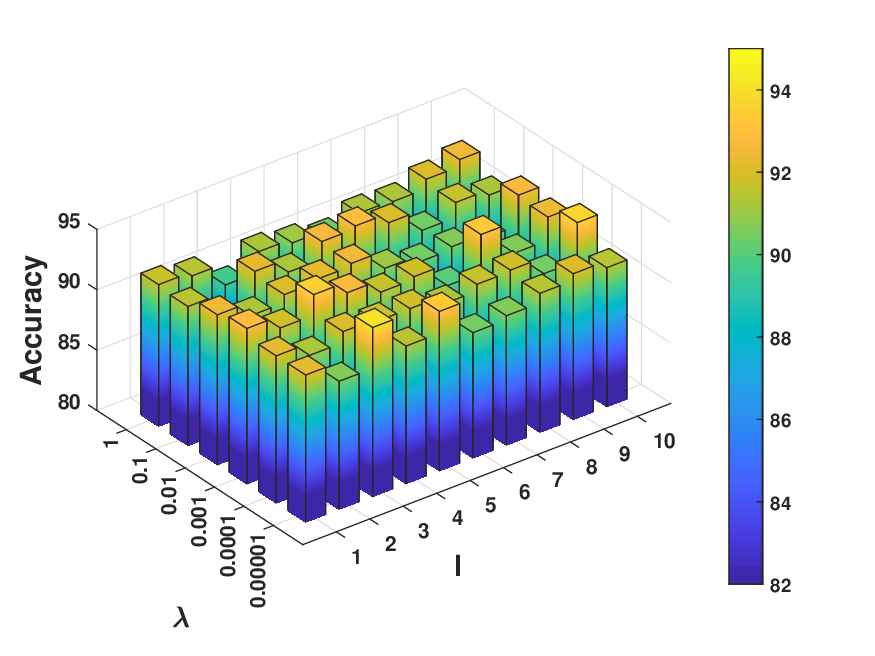}
    }\hspace{-2mm}
    \subfigcapskip=-1pt
    \subfigure[Handwritten]{
        \centering
        \includegraphics[width=4.3cm]{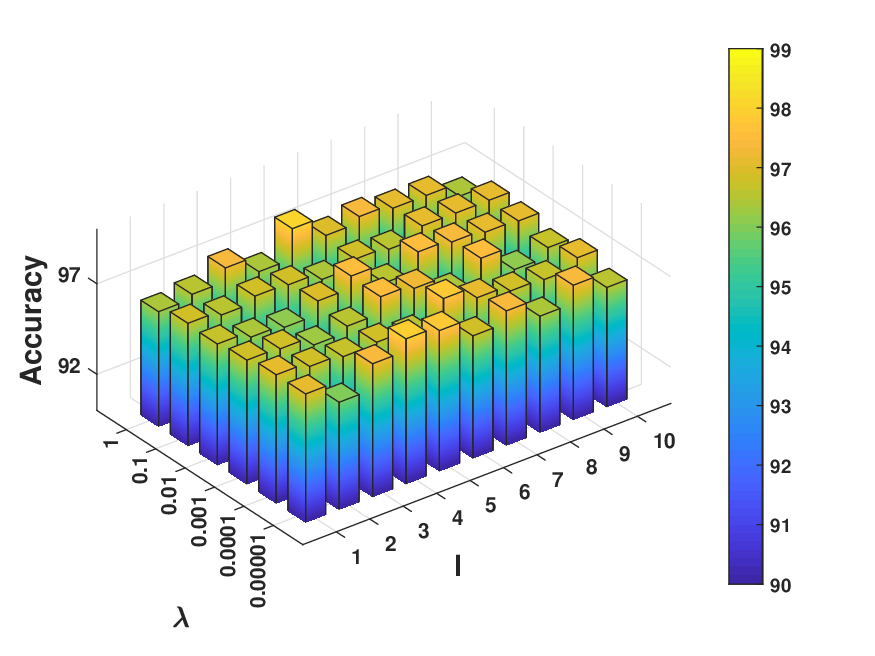}
    }\hspace{-2mm}
    \subfigcapskip=-1pt
    \subfigure[MNIST]{
        \centering
        \includegraphics[width=4.3cm]{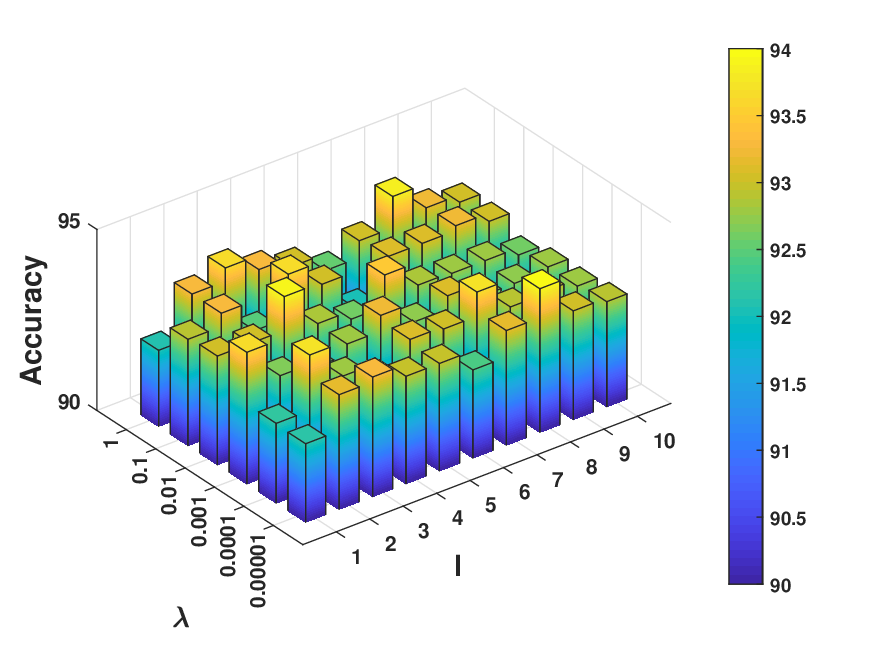}
    }\hspace{-2mm}
    \subfigcapskip=-1pt
    \subfigure[Animal]{
        \centering
        \includegraphics[width=4.3cm]{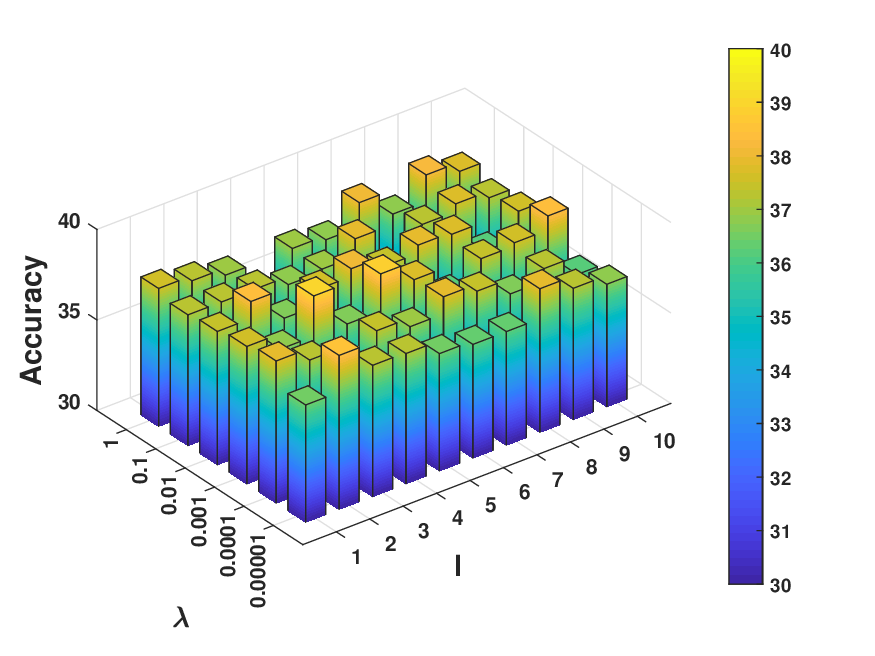}
    }\hspace{-2mm}
      \vspace{-0.1cm}
\caption{Impact of different parameters on different datasets.}\label{para}
\end{figure*}

\begin{figure*}[t]
    \centering
    \subfigcapskip=-1pt
     \subfigure[3Sources]{
        \centering
        \includegraphics[width=4.3cm]{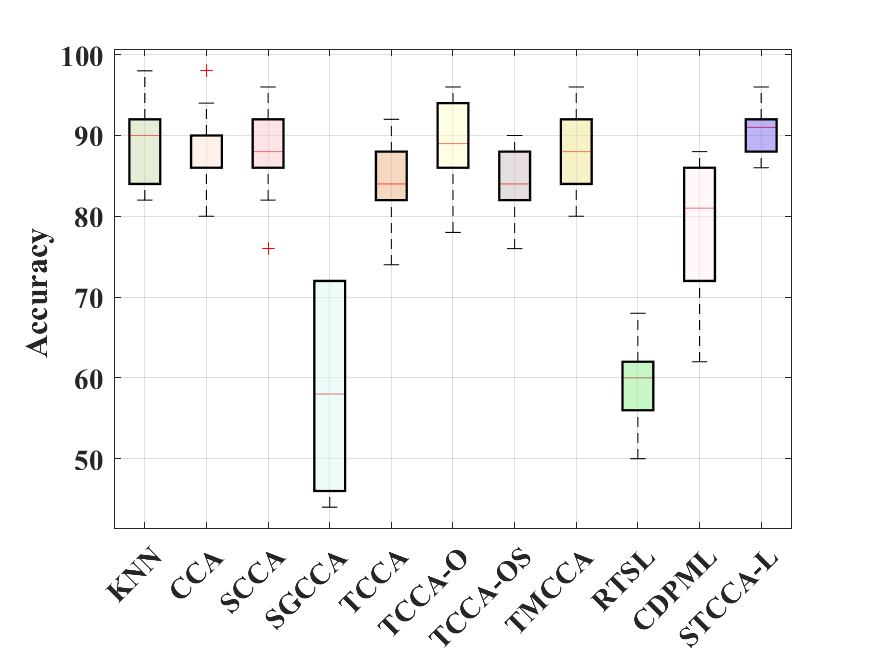}
    }\hspace{-2mm}
    \subfigcapskip=-1pt
    \subfigure[MSRC]{
        \centering
        \includegraphics[width=4.3cm]{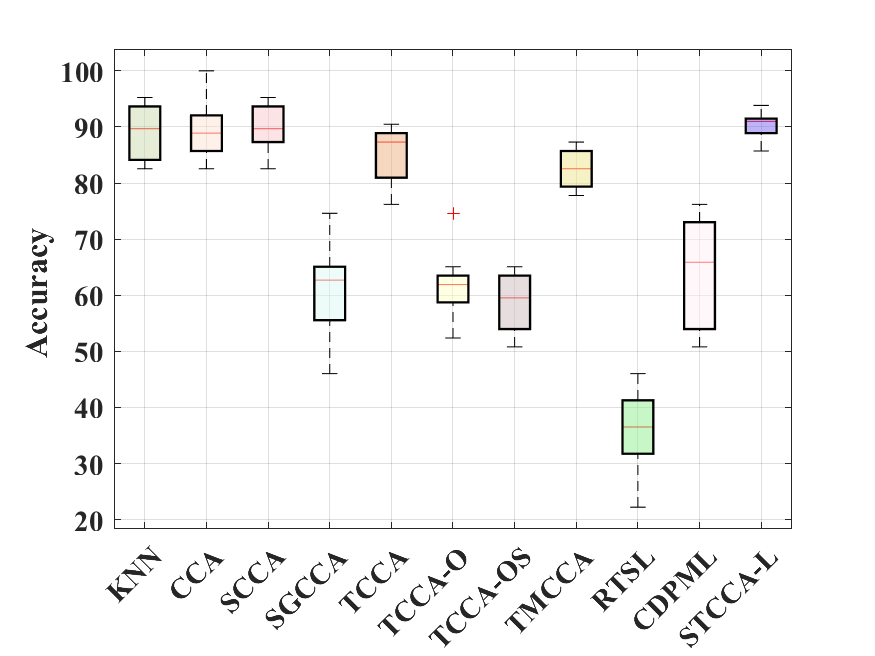}
    }\hspace{-2mm}
    \subfigcapskip=-1pt
    \subfigure[BBCsport]{
        \centering
        \includegraphics[width=4.3cm]{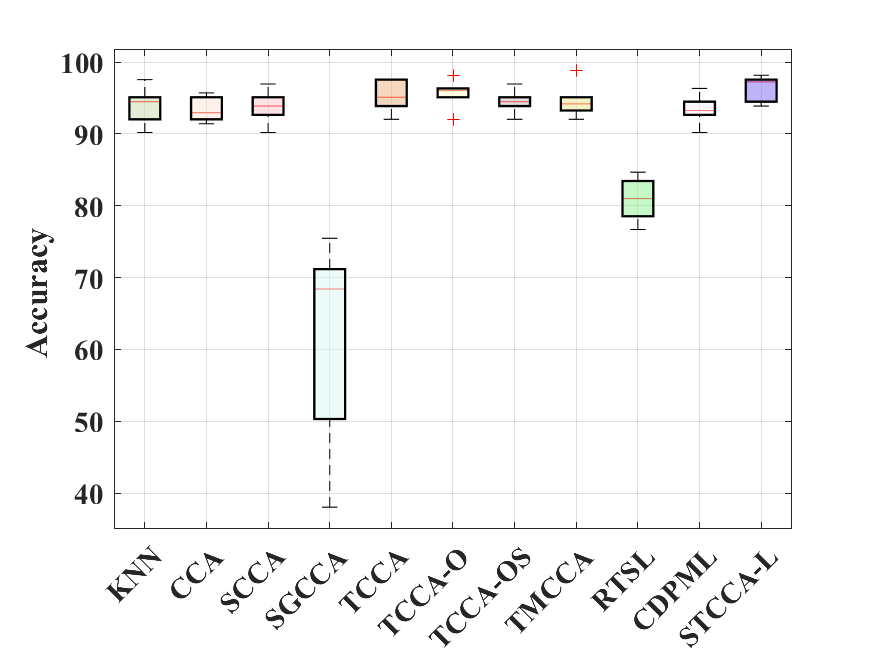}
    }\hspace{-2mm}
    \subfigcapskip=-1pt
    \subfigure[Reusters]{
        \centering
        \includegraphics[width=4.3cm]{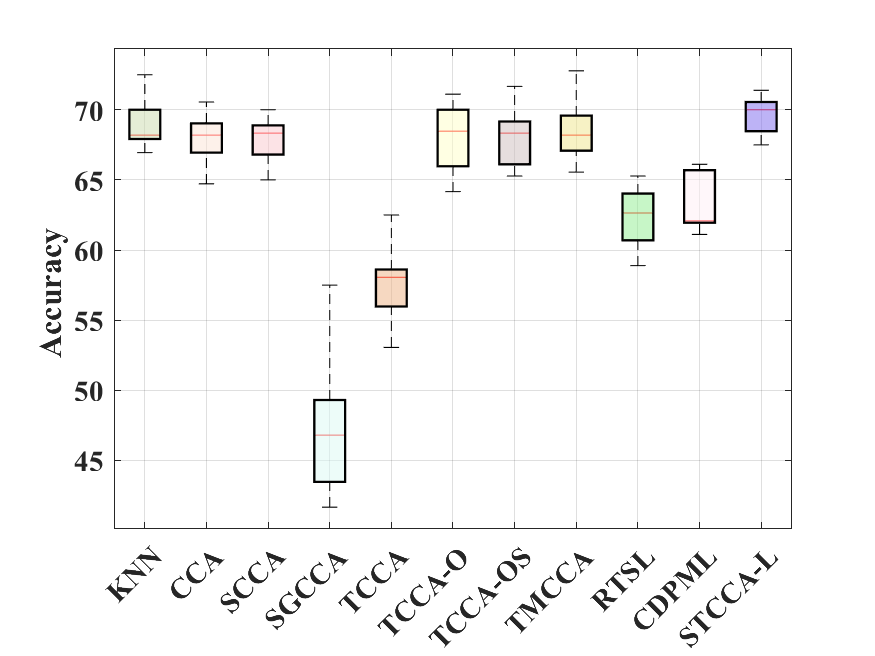}
    }\hspace{-2mm}

    \subfigure[Caltech101]{
        \centering
        \includegraphics[width=4.3cm]{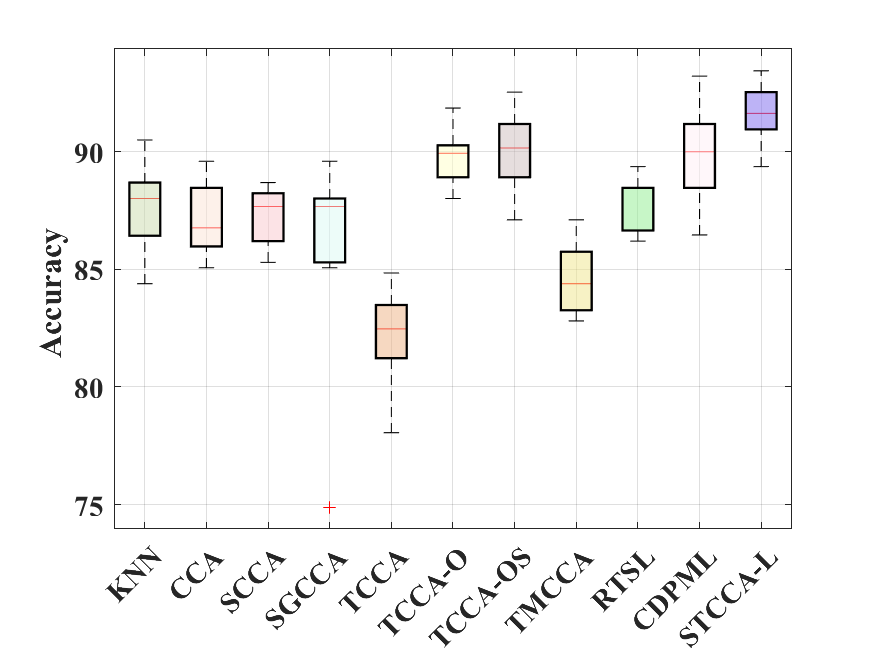}
    }\hspace{-2mm}
    \subfigcapskip=-1pt
    \subfigure[Handwritten]{
        \centering
        \includegraphics[width=4.3cm]{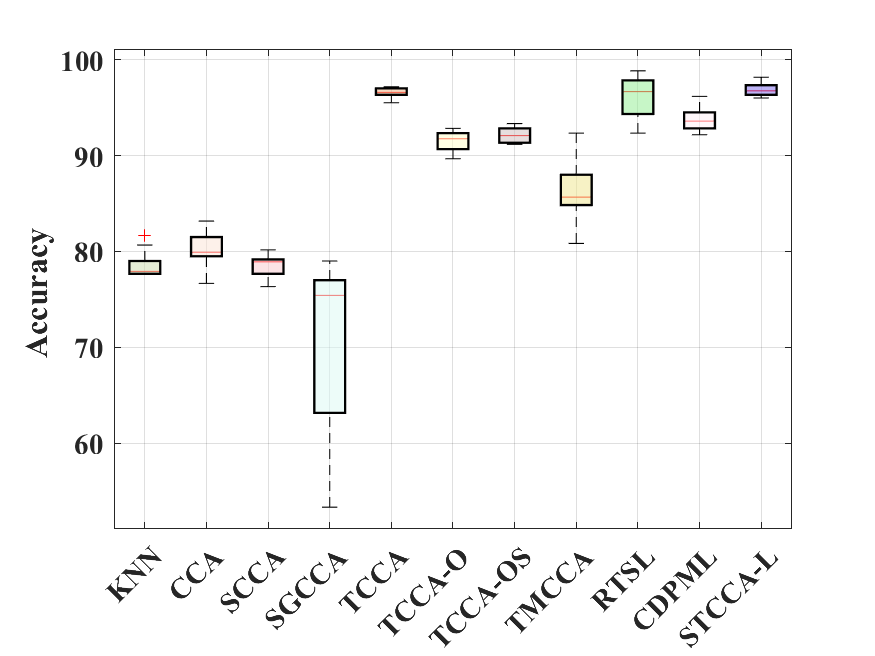}
    }\hspace{-2mm}
    \subfigcapskip=-1pt
    \subfigure[MNIST]{
        \centering
        \includegraphics[width=4.3cm]{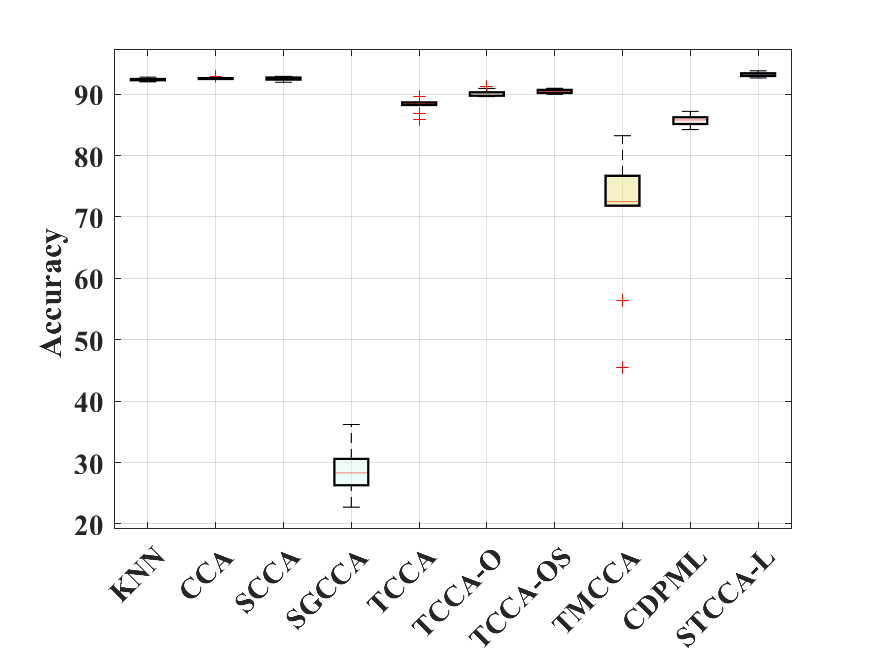}
    }\hspace{-2mm}
    \subfigcapskip=-1pt
    \subfigure[Animal]{
        \centering
        \includegraphics[width=4.3cm]{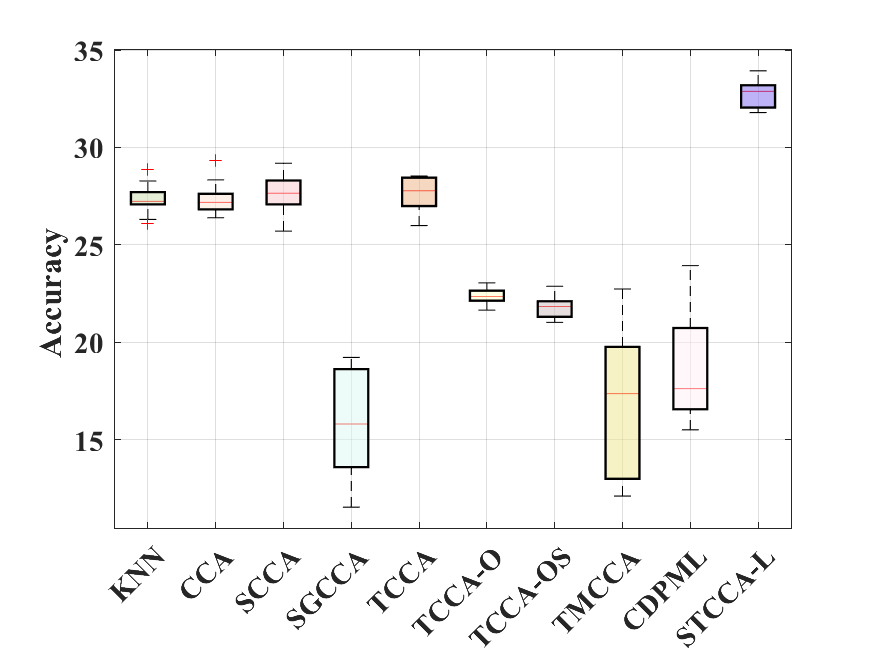}
    }\hspace{-2mm}
     \vspace{-0.1cm}
\caption{Visualization of model stability analysis on different datasets.}\label{box}
\end{figure*}

\begin{table*}[t]
    \renewcommand\arraystretch{1.2} 
    \caption{Time consuming (s) of all tensor CCA methods. }\label{time} 
    \centering
 \setlength{\tabcolsep}{2.0mm}
    \vspace{-0.1cm}
    \begin{tabular}{|c|c|c|c|c|c|c|c|c|}
      \hline
     ~~~{{Methods}}~~~ &~{3Sources}~& ~{MSRC}~  & ~{BBCsport}~& ~{Reusters}~ & ~{Caltech101}~&~{Handwritten}~ &~{MNIST}~ & ~{Animal}~ \\
      \hline \hline
     { TCCA}&  0.12&$\mathbf{0.89}$ &\underline{0.18}&\underline{7.56}& 1.90&\underline{6.58}& $\mathbf{3.18}$& \underline{10.75}\\ \hline
    { TCCA-O}&\underline{0.11}&\underline{1.01}&0.18&7.77& \underline{1.35}&6.79&\underline{3.31}&10.95\\ \hline
   {	TCCA-OS}& 0.69& 5.71& 0.53&10.95& 4.03&8.42& 3.82& 25.83\\ \hline
    {  TMCCA}& 0.34& 15.98& 1.32&45.64& 2.41&26.53& 54.39& 97.90\\ \hline
    	STCCA-L (Our)&$\mathbf{0.10}$&1.08&$\mathbf{0.17}$&$\mathbf{7.47}$& $\mathbf{1.27}$&$\mathbf{6.43}$& 3.59&$\mathbf{10.60}$\\
      \hline
    \end{tabular}
\end{table*}
\subsection{Experimental Results}\label{B}
Table \ref{acc} and Table \ref{f1score} list the classification accuracy and F1 scores, respectively, of the proposed STCCA-L compared to other state-of-the-art methods on eight multi-view datasets. The best and second-best results are highlighted in bold and underlined, respectively. It can be observed that,

\begin{itemize}
    \item The proposed STCCA-L outperforms other state-of-the-art methods in terms of both classification accuracy and F1 score on most datasets, demonstrating its effectiveness and superiority. For example, on the Animals, MSRC, and 3Sources datasets, STCCA-L achieves accuracy improvements of 5.29$\%$, 4.76$\%$, and 4.50$\%$, respectively, over the second-best method, along with F1 score gains of 3.41$\%$, 5.37$\%$, and 6.77$\%$. Moreover, t-SNE is a dimensionality reduction technique primarily used for visualizing high-dimensional data in a lower-dimensional space. Fig. \ref{tnse} visualizes the t-SNE results of all methods on the MSRC dataset. It can be observed that among these methods, STCCA-L has the most points of the same color, and its classification results are the most distinct.
    \item Compared with the baseline KNN, the classification performance of the multi-view subspace method has basically improved. Compared with other multi-view subspace methods such as RTSL and CDPMVL, the method based on CCA has more robust performance. It is worth noting that RTSL shows an error of insufficient memory on large datasets, indicating that this method is not applicable to large datasets. Among the matrix CCA methods, which include CCA, SCCA, and SGCCA, SCCA achieves the highest performance, which can be attributed to its incorporation of sparse regularization. Compared with matrix CCA methods, the proposed STCCA-L demonstrates significant advantages. For instance, on the Handwritten  dataset, STCCA-L improves the accuracy and F1 score by 11.02$\%$ and 18.96$\%$, respectively, over the best matrix CCA method. The main reason is that the covariance tensor of STCCA-L captures a more comprehensive multi-view relationship.
    \item The tensor CCA methods, which include TCCA, TCCA-O, TCCA-OS, TMCCA, and STCCA-L, outperform the matrix CCA methods in terms of classification performance. However, the original TCCA does not provide satisfactory results on datasets such as BBCSport, Reuters, and MNIST, primarily due to the absence of regularization mechanisms like orthogonality and sparsity constraints. In contrast, the proposed STCCA-L integrates orthogonal regularization, sparse regularization, and Laplacian regularization, allowing it to more effectively capture rich and discriminative information from each view. As a result, STCCA-L consistently achieves strong classification performance across all eight datasets. Notably, on the MSRC  dataset with 5 views, STCCA-L improves classification accuracy and F1 score by 8.41$\%$ and 7.82$\%$, respectively, compared with the best-performing alternative tensor CCA method.
	\end{itemize}
    
Fig. \ref{dim} is a line graph of accuracy with error bars, reflecting the results of classification accuracy in different dimensions after being processed by different methods. The classification accuracy of STCCA-L is significantly higher than that of other methods on eight datasets. In terms of the trend, STCCA-L has a more stable trend with the increase in the number of extracted features. For example, on the Caltech101 dataset, TCCA-O and TCCA-OS have a decreasing trend when the number of extracted features is greater than 8, and there may be feature redundancy. Our proposed STCCA-L effectively avoids feature redundancy because it can make good use of graph information. From the analysis of error bar size, STCCA-L is also significantly smaller than other methods, and the classification results are more accurate.

\subsection{Ablation Studies}\label{C}
\subsubsection{Contribution of Each Group}
The proposed STCCA-L integrates the orthogonal constraint, the structural sparse regularization, and the Laplacian regularization in a unified framework. To demonstrate their effectiveness, it conducts ablation studies on the 3Sources, MSRC, and Caltech101 datasets. 
Table \ref{ablation} displays the control group set and the performance of these groups.
 It is evident that the proposed STCCA-L without the sparse regularization or the Laplacian regularization has worse classification results than STCCA-L. Furthermore, the classification results of the proposed  STCCA-L without orthogonality constraints are quite different from those of STCCA-L. Therefore, the orthogonality constraint is a crucial part of the proposed STCCA-L.

 Fig. \ref{Ablation} visualizes the classification confusion matrices of STCCA-L and its three degradation models on the 3Sources, MSRC, and Caltech101 datasets. The confusion matrix is a situation analysis table for the prediction results of a classification model. It summarizes the records in the dataset in matrix form based on two criteria: the true category and the predicted category. The rows of the matrix represent the true values, and the columns of the matrix represent the predicted values. The diagonal structure of the confusion matrix indicates that the prediction results of the classification model are close to the true values. As shown in Fig. \ref{Ablation}, our proposed  STCCA-L presents nearly perfect diagonal structures on three datasets, demonstrating its excellent classification performance and indicating the necessity of orthogonal constraints, Laplacian regularization, and sparse regularization.
 
\subsubsection{Influence of Graph Construction} 
To validate the effectiveness of the graph construction strategy, the adaptive neighbor graph is adopted. The initial graph $\mathbf{W}_p$ is constructed and evaluated using different graph methods, including Gaussian kernel, KNN, cosine similarity, and sparse representation. As shown in Table \ref{abgraph}, the proposed method demonstrates consistent performance across different graph methods, indicating its robustness to the selection of the initial weight matrix. It is worth noting that the adaptive neighbor graph achieved the best results on all datasets, highlighting its ability to mitigate the impact of noise and bias in the initial graph by dynamically adjusting the neighbor weights.

\subsubsection{Effect of Algorithm Initialization} 
Furthermore, to verify the advantages of using the random matrix initialization, it conducts an ablation study using four initialization strategies: SVD, identity, orthogonal, and random. The results are summarized in Table \ref{abIni}. It can be seen that all strategies have achieved comparable performance under minor fluctuations between initializations, indicating that the developed algorithm is robust to initialization. It is worth noting that the results of random initialization on most datasets are slightly better, and thus it is adopted as the default initialization strategy in this paper.

  \subsection{Discussion}\label{E}
\subsubsection{Robustness Verification}\label{D}

This section presents experiments on noisy datasets, an aspect often overlooked by most CCA classification methods. However, evaluating performance under noisy conditions is crucial and deserves further investigation. To assess the robustness of the proposed STCCA-L, varying proportions (10$\%$-60$\%$) of Gaussian noise are added to the eight original multi-view datasets, resulting in eight noisy multi-view datasets. Fig. \ref{noise} shows the classification accuracy of all compared methods on these noisy datasets.

 The classification performance of all methods on noisy datasets has declined to varying degrees. Although the performance of the proposed STCCA-L has also declined, compared with other methods, the extent of its decline is relatively small. For instance, on the MNIST dataset, the classification performance of our method is almost unaffected by noise. On the Caltech101 dataset with 20$\%$, 30$\%$, and 40$\%$ Gaussian noise, the proposed STCCA-L improves the classification accuracy by at least 2.49\%, 3.90\%, and 5.60\%, respectively.

       \begin{table}[t]
	\renewcommand\arraystretch{1.2} 
	\caption{Ablation studies of the SSN method. }\label{ablationsub} 
	\centering
	\setlength{\tabcolsep}{1.0mm}
	\vspace{-0.1cm}
			\begin{tabular}{|c|cc|cc|}
				\hline
				~~~\multirow{2}{*}{{Datasets}}~~~ &\multicolumn{2}{c|}{{ADMM}}&\multicolumn{2}{c|}{{SSN}}\\
				{~}&{Total } &{Subproblem } &{Total } &{Subproblem }\\
				\hline \hline
				{3Sources}& $8.85\times 10^{-3}$  & $2.21\times 10^{-3}$&$\mathbf{6.09\times 10^{-3}}$ & $\mathbf{0.26\times 10^{-3}}$  \\ \hline
				{MSRC}& $1.19\times 10^{-2}$ & $0.35\times 10^{-3}$&$\mathbf{6.65 \times 10^{-3}}$ & $\mathbf{0.15\times 10^{-3}}$    \\
				\hline
                {Caltech101}& $4.91\times 10^{-2}$ & $2.12\times 10^{-2}$ & $\mathbf{3.93\times 10^{-2}}$ &$\mathbf{1.09\times 10^{-2}}$   \\
				\hline
			\end{tabular}
		\end{table}

\subsubsection{Parameter Analysis}
This section evaluates the parameter sensitivity of the proposed STCCA-L and selects the best parameters. Our method has two parameters, \textit{i.e.}, $l$ and $\lambda$, which must be chosen carefully. $l$ represents the maximum order of the multi-order graph, and $\lambda$ represents the importance of the sparse structure. It first sets a range empirically and then chooses a set of parameter values with the best classification performance from this range. Their variation ranges are $\lambda=\{0.00001,0.0001,\cdots,1\}$ and $l=\{1,2,\cdots,10\}$, respectively. Fig. \ref{para} shows the classification accuracy results under different parameters on the eight datasets.

As shown in Fig.~\ref{para}, decreasing $\lambda$ tends to enhance accuracy for the majority of datasets. For example, on the BBCSport dataset, the accuracy reaches 98.20$\%$ at $l = 3 $ and $\lambda= 0.0001$, significantly outperforming 94.50$\%$ obtained when $\lambda= 1$. The parameter $l$ shows no consistent trend, but optimal performance tends to occur when $l$ is in the range of $3$ to $7$. In contrast, extreme values such as $l=1$ or $l=10$ often result in suboptimal performance.

\subsubsection{Stability Analysis}\label{F}
The stability of our model is analyzed using box plots on the eight datasets. In terms of model stability, tensor CCA methods, \textit{i.e.}, TCCA-O and STCCA-L, are significantly superior to other methods. Compared with TCCA-O, our proposed STCCA-L has higher classification accuracy. Therefore, STCCA-L has stable classification results compared with other competing models.

\subsubsection{Time Consuming}\label{H}

Table \ref{time} presents the average CPU time consumption of all tensor CCA methods on eight datasets. It can be seen that the proposed STCCA-L achieves competitive time costs on most datasets. While TCCA and TCCA-O generally have lower runtimes, their classification performance is poor. Note that TMCCA has the slowest runtime on MNIST and Animal, indicating that its computational overhead is unsuitable for practical applications on large-scale datasets. Therefore, our method has good computational efficiency while taking into account performance.

Next, the ablation study of the SSN method for solving the subproblem is added to quantitatively evaluate the optimization efficiency of different methods. Specifically, the total computing time and subproblem resolution time of the SSN method and the ADMM on three datasets are compared. As shown in Table \ref{ablationsub}, the SSN method always has a lower computational cost than ADMM, demonstrating its superior efficiency. These results verify the contribution of the SSN method to the overall performance of our alternating manifold proximal gradient algorithm.

\begin{figure}[t]
	\centering
	\subfigcapskip=-1pt
	\subfigure[MNIST]{
		\centering
		\includegraphics[width=4cm]{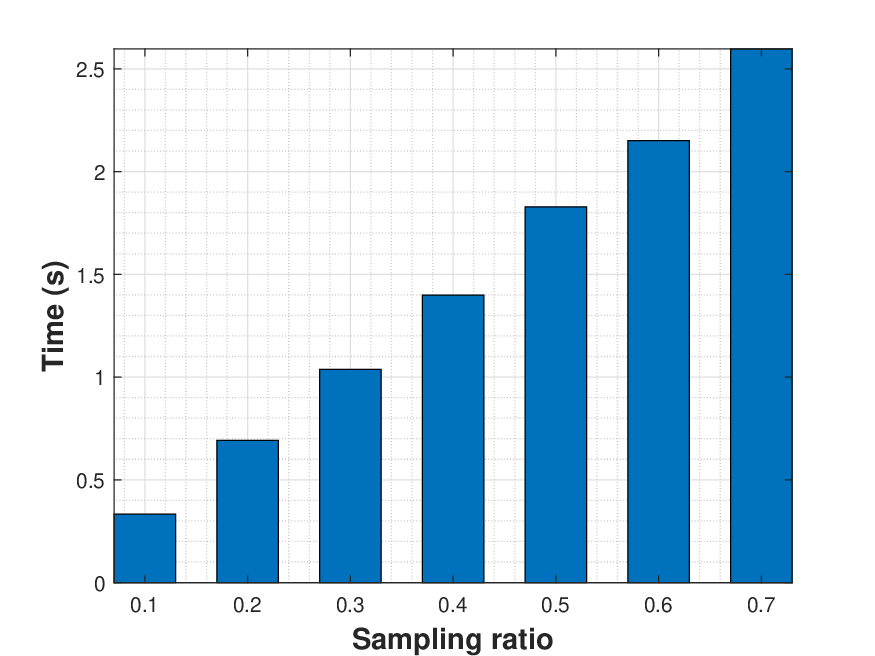}
	}\hspace{-2mm}
	\subfigcapskip=-1pt
	\subfigure[Animal]{
		\centering
		\includegraphics[width=4cm]{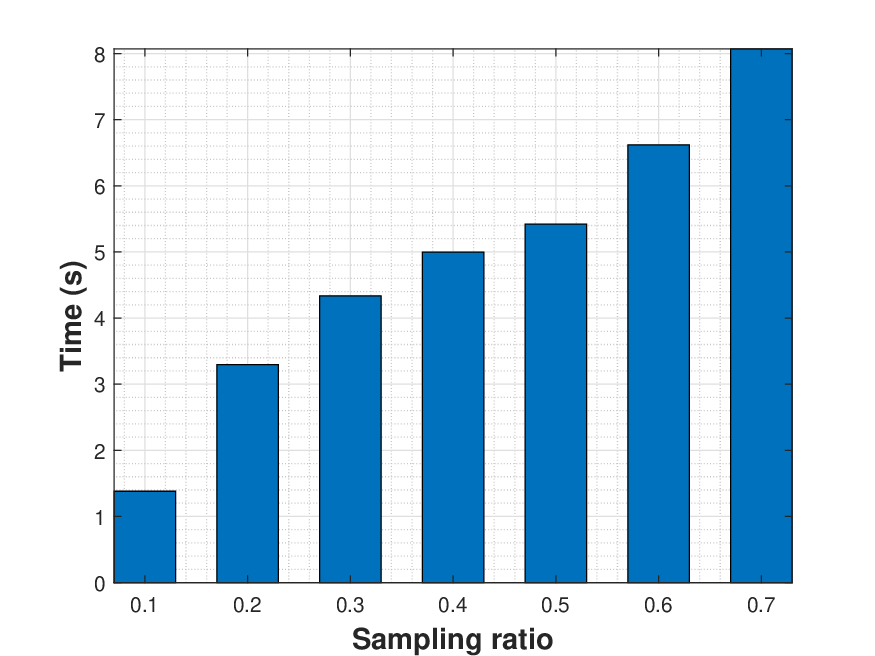}
	}\hspace{-2mm}
	\vspace{-0.1cm}
	\caption{Runtime analysis under different sampling ratios.}\label{sample}
\end{figure}

In addition, on large datasets (\textit{i.e.}, MNIST and Animal), by changing the sampling ratio, the sample size $N$ is effectively changed while other parameters remained unchanged. The results of the running time are shown in Fig.~\ref{sample}. The bar chart indicates that on the MNIST and Animal datasets, the total running time increases approximately linearly with $N$. This empirical trend is consistent with the theoretical complexity of the algorithm when other parameters are fixed.


\section{Conclusion}\label{sec5}
In this paper, we address the issues of feature redundancy and the neglect of individual view information in existing TCCA methods by proposing STCCA-L, a novel method that incorporates sparse regularization on canonical matrices and Laplacian regularization of multi-order graphs. To solve the resulting optimization problem, we develop an alternating manifold proximal gradient algorithm, further accelerated with the SSN method. We theoretically prove that the sequence generated by our algorithm converges to a stationary point. Experimental results on real-world datasets demonstrate the superiority of the proposed method.

In the future, we are interested in extending the proposed method to distributed settings \cite{guo2024adaptive} to accommodate scenarios where multi-view data may come from independent sources. Additionally, developing efficient optimization algorithms based on deep unfolding networks \cite{Deng2025} to enable automatic parameter learning is also an area worth exploring.

\section*{Acknowledgment}
 
The authors would like to thank the Associate Editor and anonymous reviewers for their numerous constructive comments. They would also like to thank Prof. Cheng Liang, Prof. Siddartha Reddy, and Prof. Chenglong Zhang for sharing the code of related papers.

	\begin{IEEEbiography}
	[{\includegraphics[width=1in,height=1.25in,clip,keepaspectratio]{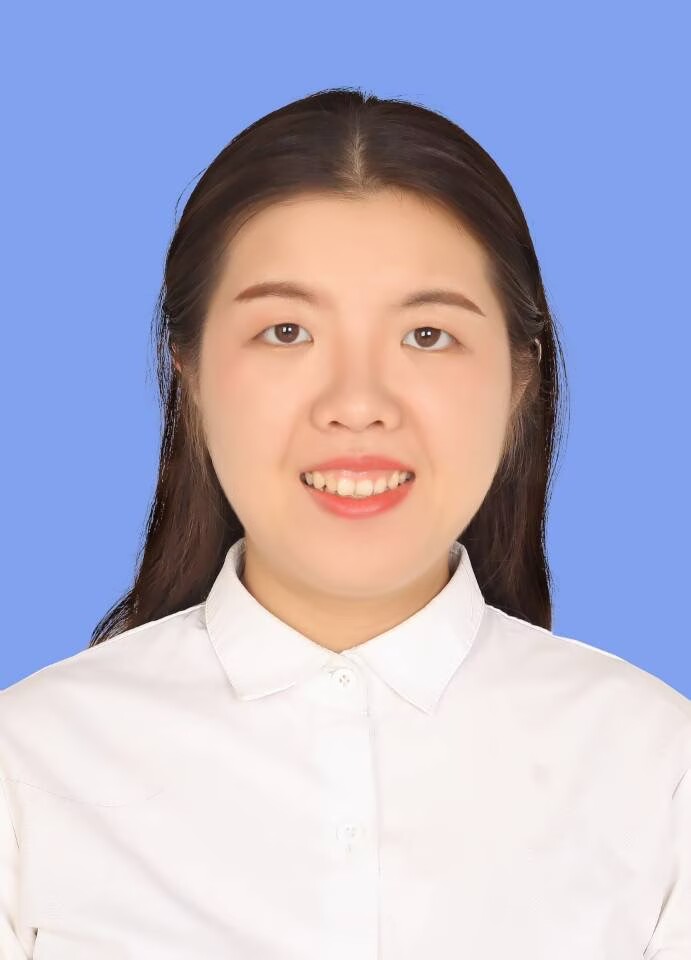}}] 
	{Yanjiao Zhu} received the Ph.D. degree in Statistics from Qufu Normal University, China, in 2024. She is a Postdoctoral Researcher at the School of Intelligent Systems Engineering, Sun Yat-sen University, Guangzhou, China.  \\
	Her current research interests include machine learning and pattern recognition.
\end{IEEEbiography} 

\begin{IEEEbiography}	[{\includegraphics[width=1in,height=1.25in,clip,keepaspectratio]{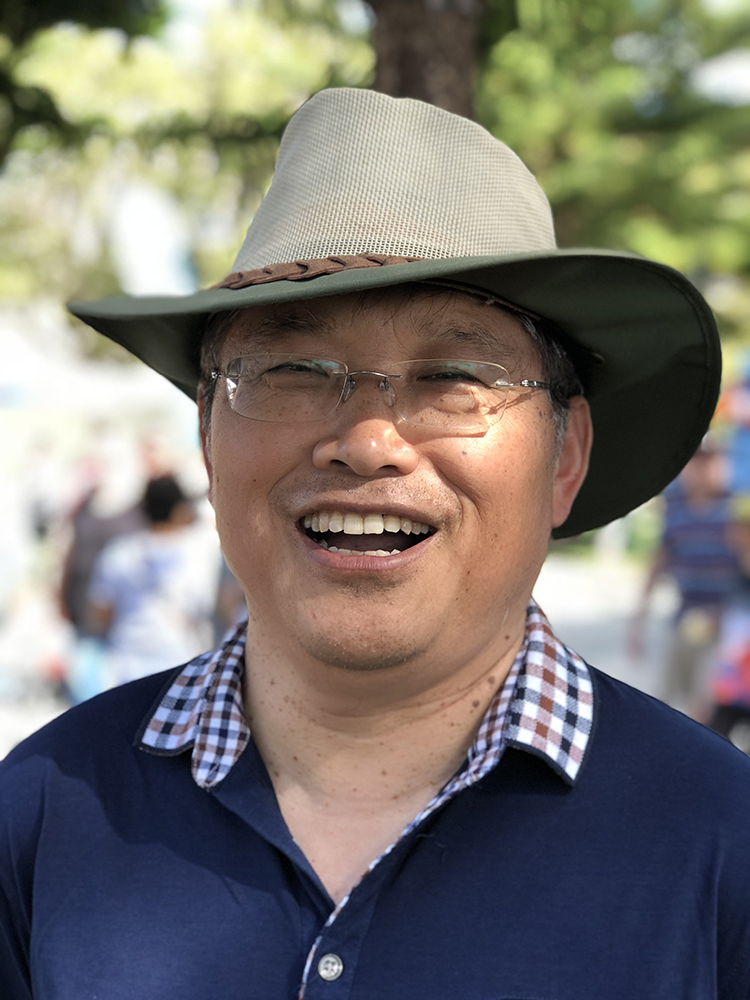}}] 
	{Wanquan Liu} received the B.S. degree in Applied Mathematics from Qufu Normal University, China, in 1985, the M.S. degree in Control Theory and Operation Research from Chinese Academy of Science in 1988, and the Ph.D. degree in Electrical Engineering from Shanghai Jiaotong University, in 1993. He once held the ARC Fellowship, U2000 Fellowship and JSPS Fellowship and attracted research funds from different resources over 2.4 million dollars. He is currently a Full Professor at the School of Intelligent Systems Engineering, Sun Yat-sen University, Guangzhou, China. 
	\\
	His current research interests include large-scale pattern recognition, signal processing, machine learning, and control systems.
\end{IEEEbiography}	

\begin{IEEEbiography}	[{\includegraphics[width=1in,height=1.25in,clip,keepaspectratio]{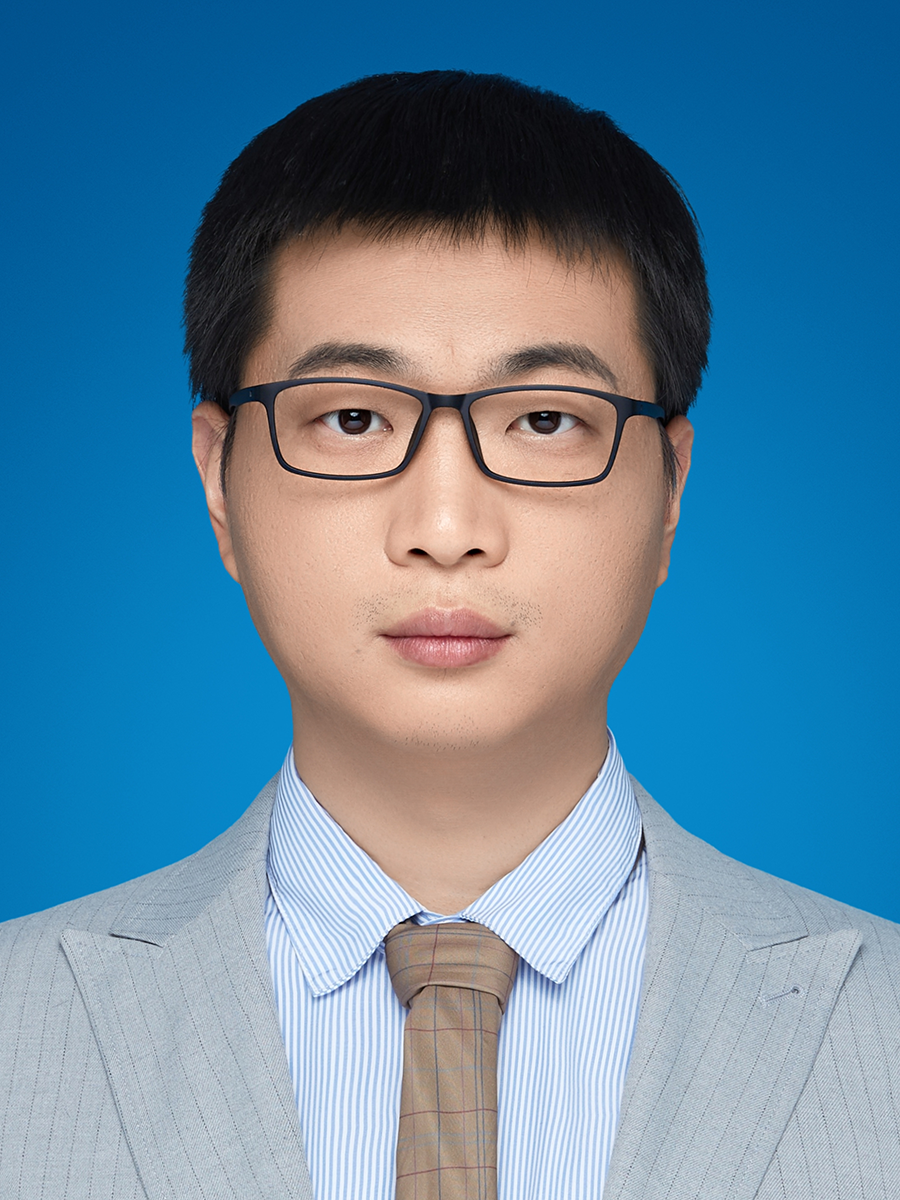}}] 
	{Xianchao Xiu} received the Ph.D. degree in Operations Research from Beijing Jiaotong University, China, in 2019. From June 2019 to May 2021, he worked as a Postdoctoral Researcher at Peking University, China. He is an Associate Professor at the School of Mechatronic Engineering and Automation, Shanghai University, China. 
	\\
	His current research interests include sparse optimization, signal processing, deep learning, and large language models. 
\end{IEEEbiography}

\begin{IEEEbiography}
[{\includegraphics[width=1in,height=1.25in,clip,keepaspectratio]{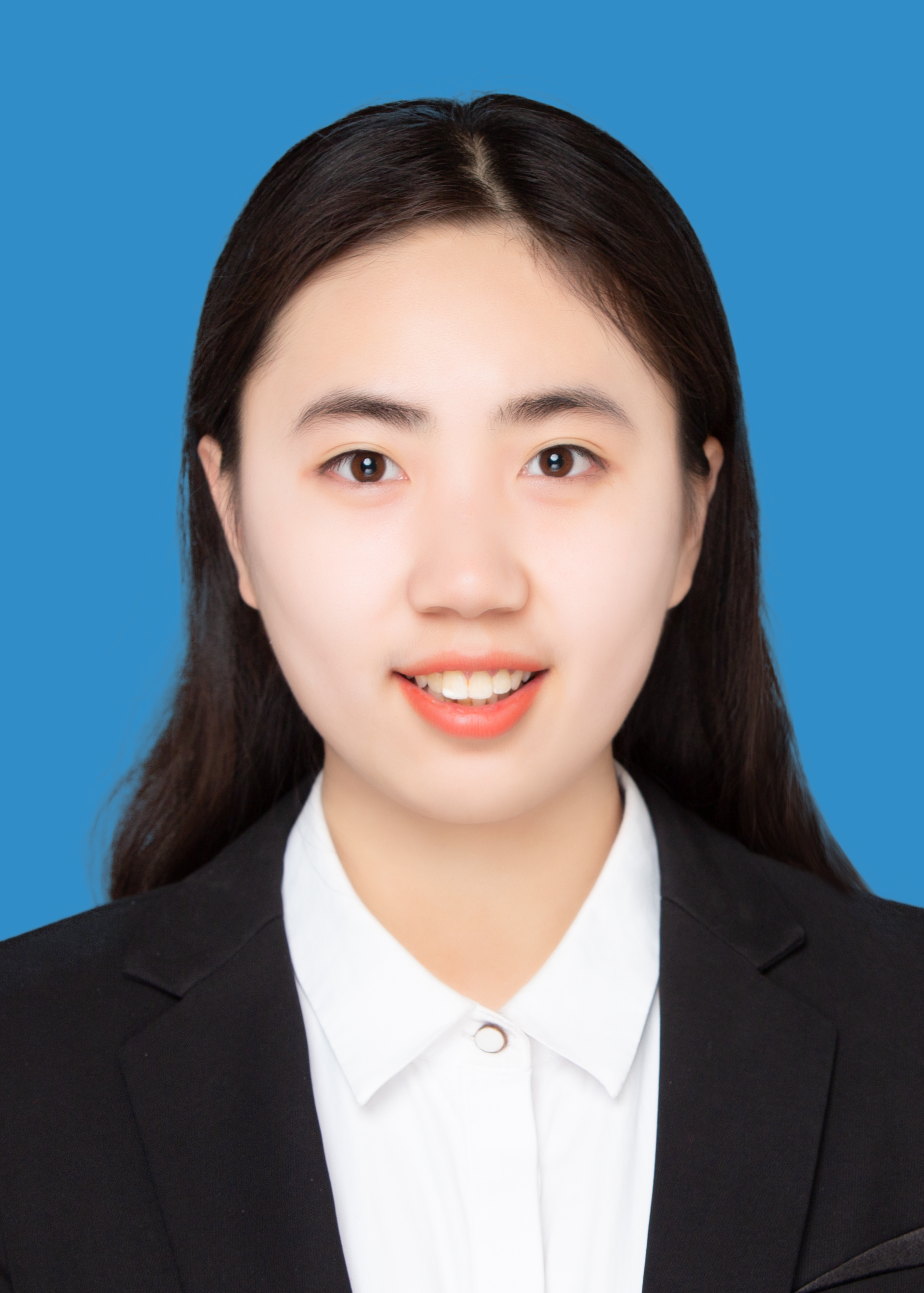}}] 
	{Jianqin Sun} is currently pursuing the M.S. degree with the School of Mathematics and Statistics, Beijing Jiaotong University, Beijing, China. 
	\\
	Her current research interests include tensor representation and sparse optimization.
\end{IEEEbiography}
\end{document}